\documentstyle[amscd,amssymb,verbatim,diagrams,12pt]{amsart}
\newcommand{\Perf}{\operatorname{Perf}}
\newcommand{\ev}{{\operatorname{ev}}}

\renewcommand{\Bar}{\operatorname{Bar}}

\newcommand{\Si}{\Sigma}

\renewcommand{\mod}{\operatorname{mod}}

\newcommand{\Cone}{\operatorname{Cone}}

\newcommand{\gr}{\operatorname{gr}}
\newcommand{\cusp}{{\operatorname{cusp}}}

\newcommand{\Tor}{\operatorname{Tor}}
\newcommand{\und}{\underline}

\newcommand{\OO}{{\cal O}}
\newcommand{\RR}{{\cal R}}

\newcommand{\qgr}{\operatorname{qgr}}

\newcommand{\be}{{\bf e}}

\newcommand{\II}{{\cal I}}

\newcommand{\BB}{{\cal B}}

\newcommand{\Tw}{\operatorname{Tw}}
\newcommand{\tors}{\operatorname{tors}}

\newcommand{\G}{{\Bbb G}}

\newcommand{\mg}{{\frak m}}
\newcommand{\hra}{\hookrightarrow}
\newcommand{\lan}{\langle}
\newcommand{\ran}{\rangle}
\newcommand{\Coh}{\operatorname{Coh}}
\newcommand{\GG}{{\cal G}}
\newcommand{\CC}{{\cal C}}

\newcommand{\tr}{\operatorname{tr}}

\newcommand{\Spec}{\operatorname{Spec}}

\newcommand{\Mat}{\operatorname{Mat}}

\newcommand{\Proj}{\operatorname{Proj}}
\renewcommand{\P}{{\Bbb P}}

\newcommand{\si}{\sigma}

\newcommand{\ga}{\gamma}
\newcommand{\de}{\delta}

\renewcommand{\ker}{\operatorname{ker}}
\newcommand{\im}{\operatorname{im}}

\newcommand{\A}{{\Bbb A}}

\numberwithin{equation}{subsection}
\newcommand{\HHom}{{\cal H}om}

\newcommand{\GL}{\operatorname{GL}}

\newtheorem{thm}{Theorem}[subsection]
\newtheorem{prop}[thm]{Proposition}
\newtheorem{lem}[thm]{Lemma}
\newtheorem{cor}[thm]{Corollary}
{  \theoremstyle{definition}
\newtheorem{defi}[thm]{Definition}
\newtheorem{ex}[thm]{Example}

\newtheorem{rem}[thm]{Remark}

}

\newcommand{\Pf}{\noindent {\it Proof}}
\newcommand{\id}{\operatorname{id}}

\newcommand{\ov}{\overline}

\renewcommand{\AA}{{\cal A}}

\newcommand{\FF}{{\cal F}}
\newcommand{\EE}{{\cal E}}
\newcommand{\MM}{{\cal M}}

\newcommand{\XX}{{\cal X}}
\newcommand{\HH}{{\cal H}}
\newcommand{\PP}{{\cal P}}
\newcommand{\QQ}{{\cal Q}}
\newcommand{\VV}{{\cal V}}
\newcommand{\SS}{{\cal S}}
\newcommand{\LL}{{\cal L}}

\newcommand{\Hom}{\operatorname{Hom}}

\newcommand{\Ext}{\operatorname{Ext}}
\newcommand{\End}{\operatorname{End}}

\newcommand{\Aut}{\operatorname{Aut}}

\renewcommand{\a}{\alpha}
\renewcommand{\b}{\beta}
\newcommand{\om}{\omega}

\newcommand{\la}{\lambda}
\newcommand{\th}{\theta}
\newcommand{\C}{{\Bbb C}}

\newcommand{\Z}{{\Bbb Z}}

\newcommand{\wt}{\widetilde}
\newcommand{\ot}{\otimes}

\newcommand{\ad}{\operatorname{ad}}
\newcommand{\sub}{\subset}
\newcommand{\ed}{\qed\vspace{3mm}}

\newcommand{\PGL}{\operatorname{PGL}}

\newcommand{\cha}{\operatorname{char}}

\newcommand{\sph}{\operatorname{sph}}
\newcommand{\Ad}{\operatorname{Ad}}

\newcommand{\grmod}{\operatorname{grmod}}

\title{$A_\infty$-structures associated with pairs of $1$-spherical objects and noncommutative orders over curves}
\author{Alexander Polishchuk}
\address{University of Oregon, National Research University Higher School of Economics, and
Korea Institute for Advanced Study}
\thanks{Supported in part by the NSF grant DMS-1700642 and by the Russian Academic Excellence Project `5-100'}

\begin{document}

\begin{abstract} We show that pairs $(X,Y)$ of $1$-spherical objects in $A_\infty$-categories, such that the
morphism space $\Hom(X,Y)$ is concentrated in degree $0$, can be described by certain noncommutative orders
over (possibly stacky) curves. In fact, we establish a more precise correspondence at the level of isomorphism 
of moduli spaces which we show to be affine schemes of finite type over $\Z$.
\end{abstract}

\maketitle

\section*{Introduction}

\subsection{Main results}

The study of $A_\infty$-categories has become an important part of the study of derived categories in algebraic
geometry, especially in connection with the homological mirror symmetry. In \cite{P-ainf,P-more-pts}
we started to develop a systematic approach to
the moduli spaces of minimal $A_\infty$-structures on a given graded vector space.
In \cite{P-ainf,LP-m1n,P-more-pts} we related certain moduli spaces of $A_\infty$-stuctures to appropriate moduli
spaces of curves, and in \cite{LP-hmstori} this was applied to proving an arithmetic version of homological mirror
symmetry for $n$-punctured tori.

Philosophically, when replacing a geometric object by the corresponding $A_\infty$-category, one
enters the world of noncommutative geometry. Thus, it is natural that objects of noncommutative geometry should
appear in descriptions of more general moduli spaces of $A_\infty$-structures. 

In the present paper we 
consider examples of such moduli spaces parametrizing $A_\infty$-categories generated by pairs of
$1$-spherical objects $(X,Y)$ such that morphism space $\Hom(X,Y)$ is concentrated in degree $0$.
Note that examples of such pairs come from considering simple vector bundles on Calabi-Yau curves, as well as from
Fukaya categories of punctured surfaces. 

Recall (see \cite{ST}, \cite[I.5]{Seidel-book})
that an object $X$ of a $k$-linear $A_\infty$-category $\CC$, where $k$ be a field, 
is called {\it $n$-spherical} if $\Hom^i(X,X)=0$ for $i\neq 0,n$,
$\Hom^0(X,X)=\Hom^n(X,X)=k$, and for any object $Y$ of $\CC$ the pairing between the morphism spaces in the cohomology category $H^*\CC$,
$$\Hom^{n-i}(Y,X)\ot\Hom^i(X,Y)\to \Hom^n(X,X),$$
induced by $m_2$, is perfect.

Note that if we have a pair of $1$-spherical objects $(X,Y)$ in a minimal $A_\infty$-category, 
such that $\Hom(X,Y)$ is concentrated in degree $0$,
then the only interesting double products invovling $X$ and $Y$ are the perfect pairings
$$\Hom^1(Y,X)\ot \Hom^0(X,Y)\to \Hom^1(X,X)\simeq k, \ \ 
\Hom^0(X,Y)\ot \Hom^1(Y,X)\to \Hom^1(Y,Y)\simeq k.$$
Thus, up to an isomorphism, the graded associative algebra $\Hom(X\oplus Y, X\oplus Y)$ is determined by a
linear automorphism $g$ of the finite-dimensional space $\Hom^0(X,Y)$, which measures the difference between
the above two pairings. More precisely, fixing trivializations of $\Hom^1(X,X)$ and $\Hom^1(Y,Y)$ and a basis 
$\a_1,\ldots,\a_n$ in $\Hom^0(X,Y)$,
we get an identification of graded associative algebras
$$\Hom(X\oplus Y,X\oplus Y)\simeq \SS(k^n,g),$$
where $\SS(g)=\SS(k^n,g)$ is a certain $(2n+4)$-dimensional algebra depending on $g\in\GL_n(k)$ 
(see Sec.\ \ref{moduli-setup-sec}).
Furthermore, it is easy to see that replacing $g$ by $\la\cdot g$, where $\la\in k^*$, leads to an isomorphic algebra.

Since $X$ and $Y$ were objects of a minimal $A_\infty$-category, we get
a minimal $A_\infty$-structure on $\SS(g)$ extending the given $m_2$.
Thus, the problem of describing pairs of $1$-spherical objects $(X,Y)$ as above with $\Hom^0(X,Y)$ of dimension $n$
(in this case we refer to $(X,Y)$ as an {\it $n$-pair})
fits into the framework of \cite[Sec.\ 2.2]{P-more-pts}. As in \cite{P-more-pts} we consider
the moduli space of all minimal $A_\infty$-structures on the family of algebras $\SS(\cdot)$
over $\PGL_n$ (extending the given $m_2$).
The corresponding moduli functor $\MM_\infty(\sph,n)$ associates to a commutative ring $R$
the set of gauge equivalence classes
of minimal $R$-linear $A_\infty$-structures on an algebra of the form $\SS(R^n,g)$, where $g\in\PGL_n(R)$
(see Sec.\ \ref{moduli-setup-sec} for the precise definition).


Note that by definition, we have a natural projection
$$\MM_\infty(\sph,n)\to \PGL_n,$$ 
where we identify the affine scheme $\PGL_n$ with the corresponding functor on commutative rings.


In the case $n=2$ we need to restrict possible elements $g$, so
we consider the principal open subscheme $\PGL_2[\tr^{-1}]\sub \PGL_2$ where $\tr(g)$ is invertible.
We denote by
$$\MM_\infty(\sph,2)[\tr^{-1}]\sub \MM_\infty(\sph,2)$$ 
the corresponding subfunctor.

Our first main result, Theorem A below, relates $A_\infty$-structures in $\MM_\infty(\sph,n)$ for $n\ge 3$
(resp., $\MM_\infty(\sph,2)[\tr^{-1}]$ for $n=2$)
to certain noncommutative projective schemes in the sense of \cite{AZ}. 
Recall that for a Noetherian graded algebra $\RR$ one considers the quotient-category $\qgr \RR=\grmod-\RR/\tors$
of finitely generated graded $\RR$-modules by the subcategory of torsion modules (it should be viewed
as the category of coherent sheaves on the corresponding noncommutative scheme). 

\begin{defi}
Let $R$ be a commutative ring, $V$ a finitely generated projective $R$-module, and $\LL$ an invertible $R$-module.
For $g\in\End(V)\ot \LL$ we denote by $\End_g(V)\sub\End(V)$ the $R$-submodule of transformations $a$ such that
$\tr(ga)=0$. We define the subalgebra in $\End(V)[z]$ by
\begin{equation}\label{E-V-eq}
\EE(V,g):=\{a_0+a_1z+\ldots\in \End(V)[z] \ |\ a_0\in R\cdot \id, a_1\in \End_g(V)\}.
\end{equation}
We view $\EE(V,g)$ as a graded $R$-algebra, where $\deg(z)=1$.
\end{defi}


\medskip

\noindent
{\bf Theorem A}. {\it For $n\ge 2$,
let us consider the functor $\MM_{filt}(n)$ associating with a commutative ring $R$ the following data:
a morphism $g:\Spec(R)\to \PGL_n$ and an isomorphism class of filtered $R$-algebras $(A,F_\bullet)$ equipped with 
an isomorphism 
$$\gr^F A\simeq \EE(R^n,g)^{op},$$ 
where we denote also by $g$ the pull-back under $g$ of the universal matrix in $H^0(\PGL_n,\Mat_n(\OO)\ot \OO(1))$.
Then for $n\ge 3$, there is an isomorphism of functors
$$\MM_\infty(\sph,n)\simeq \MM_{filt}(n)$$
and each of these functors is representable by an affine scheme of finite type over $\Z$.

In the case $n=2$ we have an isomorphism of modified functors
$$\MM_\infty(\sph,2)[\tr^{-1}]\simeq \MM_{filt}(n)[\tr^{-1}]$$
where we impose the condition that $\tr(g)$ is invertible. These functors are still representable by an affine scheme of
finite type over $\Z$.

In either case, if $(A,F_\bullet)$ is the filtered $S$-algebra corresponding to an $n$-pair $(E,F)$ of $1$-spherical objects
over a Noetherian commutative ring $S$, 
then there exists an $A_\infty$-functor from $\lan E,F\ran$ to the $A_\infty$-enhancement of the derived category of 
$\qgr \RR(A)$, where $\RR(A)$ is the Rees algebra of $A$, inducing a quasi-equivalence with its image.
}

\medskip

Note that in the case $n=1$ the equivalence of Theorem A still holds if we restrict to working over $\Z[1/6]$. Indeed, this
follows from the results of \cite{LPer}, where the moduli scheme $\MM_\infty(\sph,1)$ is identified with $\A^2$. The corresponding 
pairs of spherical objects are realized geometrically as $(\OO_C,\OO_p)$, where $C$ is an irreducible curve of 
arithmetic genus $1$ and $p\in C$ is a smooth point.

Another case that has a nice geometric realization is that of $n=2$, $g=\id$. Namely, working over $\C$,
Seidel showed in \cite[Sec.\ 2]{Seidel-flux} that 
pairs of spherical objects $(X,Y)$ with $\dim \Hom^0(X,Y)=2$ can be realized as $(\OO_C,\pi^*\OO_{\P^1}(1))$,
where $\pi:C\to \P^1$ is a (possibly singular) double cover of $\P^1$.



It seems that in the case $n=\dim\Hom^0(X,Y)>2$ 
one cannot avoid using some noncommutative geometry to realize pairs of spherical
objects $(X,Y)$. However, it is still of a sufficiently simple kind. 
Namely, working over a field $k$, using the classical result of Small-Warfield \cite{Small-Warfield} on algebras of GK dimension one, we deduce that any filtered algebra appearing in Theorem A is finite over its center. 
Using this we establish a natural bijection between $\MM_\infty(\sph,n)$ and
certain orders over integral stacky curves. Here by an {\it order} over an integral stacky curve $C$ we mean a coherent
sheaf of $\OO_C$-algebras, torsion free as $\OO_C$-module, whose stalk at the generic point of $C$ is a simple $k(C)$-algebra.

We make the following definition concerning the types of stacky curves and orders we consider.

\begin{defi}\label{neat-stacky-curve-def}
A {\it neat pointed stacky curve} over $k$ is an integral $1$-dimensional proper stack $C$ over $k$ with a smooth stacky
point of the form $p=\Spec(k)/\mu_d$, such that $C\setminus \{p\}$ is an affine (non-stacky) curve. In addition we assume
that the coarse moduli space $\ov{C}$ is a projective curve satisfying $H^0(\ov{C},\OO)=k$, and
there exists an
\'etale morphism of the form $f: U/\mu_d\to C$, where $U$ is a smooth affine curve with a $\mu_d$-action and $k$-point
$q$, such that $\mu_d$ acts faithfully on the tangent space to $q$ and $f(q)=p$.
\end{defi}

Note that a neat stacky curve either has a unique stacky point $p$ (when $d>1$) or is a usual curve (when $d=1$).

\begin{defi}\label{sph-order-def} 
Let $\AA$ be an order over an integral proper stacky curve $C$, such that $h^0(C,\OO)=k$. We say that $\AA$ is 
{\it spherical} if $\AA$ is a $1$-spherical object in the perfect derived category of right $\AA$-modules, $\Perf(\AA^{op})$, 
We say that $\AA$ is
{\it weakly spherical} if $h^0(C,\AA)=h^1(C,\AA)=1$.
\end{defi}

It is easy to see that if $\AA$ is spherical then it is weakly spherical. We will prove that an order is spherical if and only if
$h^0(C,\AA)=1$ and there exists a morphism of coherent sheaves $\tau:\AA\to\om_C$ such that the pairing
$$\AA\ot \AA\to \om: x\ot y\mapsto \tau(xy)$$
is perfect (see Sec.\ \ref{spherical-dual}). 
We say that a spherical order $\AA$ is {\it symmetric} if the above pairing (which is uniquely
defined up to a scalar) is symmetric.
The importance of spherical orders is due to the fact that they give rise to cyclic $A_\infty$-structures
(see Corollary C below).


\medskip

\noindent
{\bf Theorem B}. {\it Fix a field $k$ and a vector space $V$ over $k$ of dimension $n\ge 2$. Let us consider the following
two groupoids:

\noindent
(1) filtered algebras $(A,F_\bullet)$ with a fixed isomorphism $\gr^F A\simeq \EE(V,g)^{op}$ for some $g\in\P\End(V)$
(here morphisms exist only when the corresponding elements $g\in\P\End(V)$ are equal);

\noindent
(2) data $(C,p,v,\AA,\tau,\phi)$, where $C$ is a neat pointed stacky curve, $p\simeq \Spec(k)/\mu_d$ 
a unique (smooth) stacky point on $C$, such that
$\AA$ is a weakly spherical order over $C$ with the center $\OO_C$,
such that $h^1(C,\AA(p))=0$; $v$ is a nonzero tangent vector at $p$; 
and $\phi:\AA|_p\simeq \rho_*\End(V)^{op}$ is an isomorphism of algebras,
where $\rho:\Spec(k)\to p$ is the projection.

Then the map associating to $(C,p,\AA)$ the algebra $A=H^0(C\setminus p,\AA)$ with its natural filtration
extends to an equivalence of groupoids
(1) and (2). 

\noindent
Furthermore, the element $g$ in (1) is invertible if and only the corresponding order $\AA$ is spherical.
In this case, assuming in addition that either $n\ge 3$ or $\tr(g)\neq 0$, we have an equivalence between 
the perfect derived category $\Perf(\AA^{op})$ and
the triangulated envelope of the $A_\infty$-category with an $n$-spherical pair $(X,Y)$ associated with the data (1) by Theorem A.
Under this equivalence the pair $(X,Y)$ corresponds to the pair $(\AA, \rho_*V)$ in $\Perf(\AA^{op})$.

\noindent
If either $n\ge 3$ or $\cha(k)\neq 2$, then a spherical order $\AA$ is symmetric if and only if $g$ 
is a scalar multiple of the identity matrix.
}

\medskip


Note that for all examples of spherical orders over stacky curves $C$ that we were able to construct, $C$ is the usual 
nonstacky curve. However, we do not know whether this is always the case.
 
We are particularly interested in the case when $\AA$ is symmetric because in this case we get
a cyclic $A_\infty$-structure.

\medskip

\noindent
{\bf Corollary C}. {\it Let $k$ be a field. Assume that either $n\ge 3$ or $\cha(k)\neq 2$. Then
every minimal $A_\infty$-structure on the algebra $\SS(k^n,\id)$ 
can be realized by $A_\infty$-endomorphisms of the generator $\AA\oplus \rho_*V$ of $\Perf(\AA^{op})$, 
for a symmetric spherical order $\AA$ over a neat pointed stacky curve $(C,p)$
equipped with an isomorphism $\AA|_p\simeq \rho_*\End(V)^{op}$. 
If $\cha(k)=0$ then every such $A_\infty$-structure 
is gauge equivalent to a cyclic $A_\infty$-structure.}

\medskip



As was shown in \cite{LP-YB}, cyclic $A_\infty$-structures on $\SS(k^n,\id)$ correspond to formal solutions
of set-theoretical {\it Associative Yang-Baxter Equation} (AYBE).
Applying Corollary C, we get an algebro-geometric realizition of these solutions in the category of modules over
some symmetric spherical orders. 
Elsewhere we will present an explicit construction of such orders giving rise to trigonometric nondegenerate solutions of the
AYBE that were classified in \cite{P-nodal-massey}.


\subsection{Organization of the paper}

In this paper we are dealing with the following three types of structures:

\noindent
(I) $A_\infty$-structures on $\SS(g)$, or equivalently, $n$-pairs of $1$-spherical objects;

\noindent
(II) filtered algebras as in Theorem A;

\noindent
(III) orders over neat stacky curves as in Theorem B.

In Section \ref{ainf-moduli-sec} we study (I). Here the main result is the representability of the moduli functor
$\MM_\infty(\sph,n)$ of $A_\infty$-structures by an affine scheme of finite type over $\PGL_n$. 
This is based on the criterion from \cite{P-more-pts} which requires calculation of some Hochschild cohomology.
The latter calculation is performed in Sec.\ \ref{HH-sec}
using Koszulness of algebras $\SS(g)$ (with appropriate grading).

In Section \ref{sph-pairs-sec} we explain the passage from (I) to (II).
First, in Sec.\ \ref{spherical-sec} and \ref{n-pairs-1-sph-sec}, we give some background on spherical objects.  
In \ref{algebra-E-g-sec} we establish some properties of algebras $\EE(R^n,g)$.
Then in Theorem \ref{spherical-filtered-thm}
we give a construction of a filtered algebra corresponding to an $n$-pair of $1$-spherical objects $(E,F)$.
The main idea is to consider the infinite sequence of objects 
$$E_i=T_F^i(E),$$ 
where $T_F$ is the spherical twist functor
associated with $F$, and then define a graded algebra structure on $\RR=\bigoplus_i \Hom(E_0,E_i)$. It turns out that
$\RR$ has a canonical regular central element of degree $1$, so it can be identified with the Rees algebra $\RR(A)$ of
a filtered algebra $A$.
Note that this construction 
is inspired by a construction in the work of Van Roosmalen \cite{Roosmalen}.
We also describe a natural $\RR-\RR$-bimodule structure on $\bigoplus_i \Hom^1(E_i,E_0)$ in terms
of a certain filtered automorphism $\phi$ of $A$
(see Prop.\ \ref{bimodule-prop}). 
 
In Section \ref{orders-sec}, working over a field, we explain how to pass from (II) to (III) and discuss spherical orders.
Namely, in Sec.\ \ref{filtered-alg-sph-orders-sec} we construct a noncommutative order over a stacky curve
associated with a filtered algebra such that $\gr^F A\simeq \EE(R^n,g)^{op}$. 
In Sec.\ \ref{spherical-dual} we give a criterion for an order to be spherical (see Definition \ref{sph-order-def}).
In Sec.\ \ref{AS-Gor-sec} we check the AS-Gorenstein property for the Rees algebras $\RR(A)$ of our filtered algebras $A$. 
The proof uses a special order over the cuspidal cubic consructed in \ref{special-order-sec}.

In Section \ref{application-pairs-sec} we show how to go from (II) to (I).
Namely, in Sec.\ \ref{noncomm-proj-sec} we consider the noncommutative projective scheme in the sense of Artin-Zhang \cite{AZ}
associated with the graded ring $\RR(A)$ and construct a natural pair of $1$-spherical objects in
the derived category of $\qgr \RR(A)$, the category of coherent sheaves on this noncommutative projective scheme 
(see Proposition \ref{proj-coh-prop}). It is also straightforward to see that $\qgr \RR(A)$ is equivalent to the category of right modules
over the corrresponding order $\AA$ (see Prop.\ \ref{Proj-order-prop}). However, for computations it is more convenient to use the presentation
as $\qgr \RR(A)$ and to use results of \cite{AZ}.

In Section \ref{char-sec} we prove a technical result characterizing minimal $A_\infty$-structures on the subcategory of objects $(E_i)$
in terms of double products together with a single triple product (see Cor.\ \ref{m3-1-cor}). This result is a slight generalization of the main theorem of \cite{P-ext}, where a similar
problem is studied for the subcategory of powers of a very ample line bundle.
The key computation of Hochschild cohomology is Theorem \ref{HH-gen-thm}
generalizing \cite[Sec.\ 3.1]{P-ext}.
Then in \ref{triple-prod-sec} we perform a calculation of the relevant triple product in the case of the $A_\infty$-structure coming from an $n$-pair of $1$-spherical objects.

In Section \ref{proof-sec} we put everything together to prove our main results.
In Sec.\ \ref{thmA-Noeth-sec} we prove equivalence of functors corresponding to (I) and (II) over Noetherian rings. 
The main difficulty is to prove that going from (I) to (II) and back one gets an $A_\infty$-structure equivalent to the original one. For this we use the
characterization from Sec.\ \ref{min-non-form-alg-sec}. Another technical detail is that double products on the subcategory $(E_i)$ are determined by the 
filtered algebra $A$ together with a filtered automorphism $\phi$ of $A$
which we are able to characterize uniquely (see Prop.\ \ref{filt-aut-prop}, Cor.\ \ref{H1-bimodule-cor}).

In Sec.\ \ref{moduli-thmA-sec} we prove representability of the functor corresponding to (II) by an affine scheme of finite type, and as a result,
deduce the equivalence of (I) and (II) over arbitrary rings.

In \ref{thmB-sec} we prove equivalence of (II) and (III). In addition to constructions of Section \ref{orders-sec}, we use
a technical result on stacks discussed in \ref{polarizing-sec}.

Finally, in Section \ref{cyclic-ainf-sec}
we show that symmetric spherical orders over curves give rise to cyclic $A_\infty$-structures
and use this to prove Corollary C.

\medskip

\noindent
{\it Conventions}. All $A_\infty$-algebras/categories are strictly unital.
We denote by $\hom(\cdot,\cdot)$ the morphism spaces in an $A_\infty$-category and by $\Hom(\cdot,\cdot)$ their
cohomology,
For a commutative ring $S$, and a graded $S$-module $M=\bigoplus_i M_i$, where each $M_i$ is a projective
finite dimensional $S$-module, we denote by $M^*$ the restricted dual of $M$, which is a graded $S$-module
with components $(M^*)_i=(M_{-i})^\vee=\Hom_S(M_{-i},S)$.

\medskip

\noindent
{\it Acknowledgments}. I am grateful to Yanki Lekili and Riley Casper for useful discussions,
and to James Zhang for the help with proving finiteness of injective dimension in Sec.\ \ref{AS-Gor-sec}.
I also thank the anonimous referee for valuable comments.
Part of this work was done while the author was visiting the Korean Institute for Advanced Study, the ETH Zurich
and the IHES.
He would like to thank these institutions for hospitality and excellent working conditions. 

\section{A moduli space of $A_\infty$-structures}\label{ainf-moduli-sec}

\subsection{The moduli functor}\label{moduli-setup-sec}

For basics on $A_\infty$-structures we refer to \cite{keller-ainf}.
We would like to consider minimal (strictly unital) $A_\infty$-structures on a certain family of
categories with two objects (which are generalizations of the category considered in \cite{LP-YB}).

\begin{defi}\label{S-V-g-def} 
Let $R$ be a commutative ring, $V$ a finitely generated projective $R$-module, $\LL$ an invertible $R$-module,
and $g:V\to V\ot \LL$ an $R$-linear morphism.
The graded category $\SS(V,g)$ has two objects $X$ and $Y$ and the morphisms
$$\Hom(X,Y)=\Hom^0(X,Y)=V, \ \
\Hom(Y,X)=\Hom^1(Y,X)=V^\vee=\Hom_R(V,R),$$ 
$$\Hom^0(X,X)=R\id_X, \ \
\Hom^0(Y,Y)=R\id_Y, \ \ \Hom^1(X,X)=\LL, \ \ \Hom^1(Y,Y)=R.$$
The compositions are determined by 
$$v^*\circ v=\lan v^*,g(v)\ran, \ \ v\circ v^*=\lan v^*,v\ran,$$
where $v\in V$, $v^*\in V^\vee$.
\end{defi}




In this paper we will consider minimal unital $A_\infty$-structures on the algebras $\SS(k^n,g)$, where $g\in \GL_n(k)$, 
for fixed $n\ge 2$ (with $m_2$ given above), up to a gauge equivalence. 
More precisely, this family of algebras can be viewed
as a sheaf of associative algebras $\wt{\SS}_n$ on the scheme $\GL_n$ (over $\Z$). 
Note that we have the standard action of the central $\G_m$ on $\GL_n$ and the universal element
$g$ can be viewed as a morphism of bundles over $\GL_n$, $g:\OO^n\to \OO^n\ot \chi$, compatible
with the $\G_m$-action. Here $\chi$ is the identity character of $\G_m$, $\chi(\la)=\la$.
Thus, we can descend $g$ to an isomorphism 
$$\ov{g}:\OO^n\to \OO^n\ot \OO(1)$$ 
of vector bundles over $\PGL_n=\GL_n/\G_m$.
Now it is easy to see that the sheaf of algebras $\wt{\SS}_n$ over $\GL_n$ is isomorphic to the pull-back
of the sheaf of algebras $\SS_n$ over $\PGL_n$ associated with bundle $\OO^n$ and the morphism $\ov{g}$
as in Definition \ref{S-V-g-def}.

As in \cite{P-more-pts}, we consider the following functor over $\PGL_n$ on affine schemes over $\PGL_n$.

\begin{defi} The functor $\MM_\infty(\SS_n)$ associates to a pair $(R,g)$, where $R$ is a commutative ring and
$g\in\PGL_n(R)$, 
the set of gauge equivalence classes of minimal $A_\infty$-structures on $\SS(R^n,g)$.
Note that here $g$ is viewed as a morphism $R^n\to R^n\ot g^*\OO(1)$.
\end{defi}

For each $N\ge 2$ one can similarly consider the functor
$\MM_N$ of minimal $A'_N$-structures up to equivalence, where we consider $(m_2,m_3,\ldots,m_N)$ and impose
the $A_\infty$-identities up to $[m_2,m_N]+[m_3,m_{N-1}]+\ldots$
In \cite[Thm.\ 2.2.6]{P-more-pts} we gave a general criterion for the functors $\MM_N$ to be representable by an affine scheme
and for the projection $\MM_\infty\to\MM_4$ to be a closed embedding.
Note that there is a ``bookkeeping" mistake in \cite[Thm.\ 2.2.6]{P-more-pts} in that 
the results are stated for minimal $A_N$-structures (which involve $(m_2,\ldots,m_N)$ and
the identities up to $[m_2,m_{N-1}]+\ldots$), however,  these results only become correct upon replacing the $A_N$-structures
by $A'_N$-structures (the same problem occurs in \cite[Sec.\ 4.2]{P-ainf}).

In this paper we will apply general results on moduli of $A_\infty$-structures to our functor $\MM_\infty(\SS_n)$ over $\PGL_n$, for $n\ge 3$.
In the case $n=2$ we restrict to an open subscheme of $\PGL_2$ and modify the functor accordingly.
Namely, let us consider the closed embedding
$$\Spec(\Z/2)\hra \PGL_{2,\Z/2}\hra \PGL_2$$
given by the unit over $\Z/2$, and let $U\sub \PGL_2$ be the complementary open subscheme.
Let $\SS_{2,U}$ be the restriction of the sheaf of algebras $\SS_2$ to $U$.
Since $U$ is non-affine, as in \cite[Sec.\ 2.2]{P-more-pts}, 
we first consider the functor $\wt{\MM}_\infty(\SS_{2,U})$ of the set of gauge equivalence
classes of minimal $A_\infty$-structures, and then define $\MM_\infty(\SS_{2,U})$
to be the sheafification of $\wt{\MM}_\infty(\SS_{2,U})$ with respect to the Zariski topology on the base.

Note that when defining the functor $\MM_\infty(\SS_n)$ for $n\ge 3$, we do not need to pass to the
sheafification with respect to the Zariski topology on the base $\PGL_n$ since the base is affine
(see \cite[Thm.\ 2.2.6]{P-more-pts}). The same is true about the functor $\MM_\infty(\SS_2)[\tr^{-1}]$,
corresponding to the affine subset of $g$ in $\PGL_2$ with $\tr(g)\neq 0$.


\begin{thm}\label{ainf-moduli-thm}
The functor $\MM_\infty(\SS_n)$ for $n\ge 3$ (resp., $\MM_\infty(\SS_{2,U})$) 
is representable by an affine scheme of finite type over $\PGL_n$ (resp., over $U$). Furthermore,
in both cases the projection $\MM_\infty\to \MM_4$ to the moduli space of minimal $A'_4$-structures is a
closed embedding.
\end{thm}

Note that the functor $\MM_\infty(\sph,n)$ is the absolutization of $\MM_\infty(\SS_n)$: its point over an affine scheme
$\Spec(R)$ consists of a morphism $g:\Spec(R)\to \PGL_n$, together with an element in $\MM_\infty(\SS_n)(R,g)$.
Thus, we deduce the following corollary (in the case $n=2$ we restrict to the open subset $\tr\neq 0$ of $U\sub\PGL_2$).

\begin{cor}\label{ainf-moduli-cor} 
The functor $\MM_\infty(\sph,n)$ for $n\ge 3$) (resp., $\MM_\infty(\sph,2)[\tr^{-1}]$)
is representable by an affine scheme of finite type over $\Z$.
\end{cor}

We will give a proof of Theorem \ref{ainf-moduli-thm}
in Sec.\ \ref{HH-sec} after proving vanishing of some components of the Hochschild cohomology of algebras
$\SS(k^n,g)$. The computation is based on the fact that these algebras are Koszul with respect to the grading
$\deg V=\deg V^\vee=1$: we use the well known Koszul resolution for the diagonal bimodule.
The vanishing of certain components of the Hochschild cohomology translates into some statements
about the quadratic dual algebra $\SS^!$ which we establish in Lemma \ref{E!-lem}.

Note that in \cite{LP-YB} we considered minimal $A_\infty$-structures on $\SS(R^n,\id)$, cyclic with respect 
to a natural pairing. We postpone the discussion of cyclic structures until Sec.\ \ref{cyclic-ainf-sec}.


\subsection{The dual quadratic algebra to $\SS(k^n,g)$}

Let us fix a field $k$ and an element $g\in\GL_n(k)$, where $n\ge 2$. Let us set for brevity $\SS=\SS(k^n,g)$. 
Note that we can view $\SS$ as a $K$-algebra, where $K=k\oplus k$, where the idempotents $\be_X$ and $\be_Y$
correspond to the identity elements $\id_X$ and $\id_Y$. We also denote by $\xi_X$ and $\xi_Y$ the basis elements
in $\Hom^1(X,X)$ and $\Hom^1(Y,Y)$.



If $(\a_i)$ are the elements of $\Hom(X,Y)$ corresponding to the standard basis in $k^n$ and $(\b_j)$ are
the dual elements of $\Hom(Y,X)$, then the product in $\SS$ is given by
$$\b_i\a_j=a_{ij}\xi_X, \ \ \a_j\b_i=\de_{ij}\xi_Y,$$
where $g=(a_{ij})$.
Note that the $K$-algebra $\SS$ is generated by elements $(\a_i)$ and $(\b_j)$, so we can view it as a quotient of
the path algebra of the quiver with two vertices $X$ and $Y$ and $n$ arrows in each direction corresionding to $\a_i$ and
$\b_j$. We will use two different gradings on $\SS$: $\deg$ given by $\deg(\a_i)=0$, $\deg(\b_i)=1$, and $\deg_K$ given by
$\deg_K(\a_i)=\deg_K(\b_i)=1$. We denote by $\SS_j$ the graded components of $\SS$ with respect to $\deg_K$.

We are going to use the following convention about quadratic (and Koszul) duality over $K=k\cdot\id_X\oplus k\cdot\id_Y$. 
For a quadratic $K$-algebra $A$ with generators $V_{XY}$ and $V_{YX}$ of degree $1$, and quadratic relations
$R_{XX}\sub V_{YX}\ot V_{XY}$, $R_{YY}\sub V_{XY}\ot V_{YX}$, the dual quadratic algebra has generators
$V^!_{XY}=V_{YX}^\vee$ and $V^!_{YX}=V_{XY}^\vee$ and quadratic relations
$$R^!_{XX}=A_{2,XX}^\vee\sub (V_{YX}\ot V_{XY})^\vee\simeq V^!_{YX}\ot V^!_{XY},$$
and similarly for $R^!_{YY}$.

Thus, we think of the dual algebra $\SS^!$ as the quotient of the path algebra of 
the quiver with vertices $X$, $Y$, where the direction of $\a_i^*$ (resp., $\b_i^*$)
is opposite to that of $\a_i$ (resp., $\b_i$). We denote by $\SS^!_j$ the graded components of $\SS^!$ with respect
to the grading $\deg_K(\a_i^*)=\deg_K(\b_i^*)=1$.

\begin{lem}\label{Koszul-AA-lem}
With respect to the grading $\deg_K$ the algebra $\SS$ is Koszul and the dual quadratic algebra $\SS^!$
is generated by the dual generators $(\a_i^*)$ and $(\b_i^*)$ with the only relations
$$\sum_{1\le i,j\le n} a_{ij}\a_j^*\b_i^*=0, \ \ \sum_{i=1}^n \b_i^*\a_i^*=0.$$
\end{lem}

\Pf . The algebra $\SS$ is obtained by folding from the following $\Z$-algebra $\SS^\Z=\SS^\Z(g)$ (where we use the term {\it $\Z$-algebra}
in the sense of \cite[Sec.\ 4]{BP}):
$$\SS^\Z_{i,i+1}=\begin{cases} V & i \text{ even }\\ V^\vee & i \text{ odd.}\end{cases},$$
$$\SS^\Z_{i,i+2}=k, \ \ \SS^\Z_{i,j}=0 \text{ for } j>i+2.$$
Here $V$ is the $n$-dimensional space with the basis $(\a_i)$. The multiplication is given by
the pairings 
$$V\ot V^\vee\to k: v\ot v^*\mapsto \lan v^*,v\ran, \ \ V^\vee\ot V\to k: v^*\ot v\mapsto \lan v^*,g(v)\ran.$$
There is a natural isomorphism $\SS^\Z_{i+2,j+2}\simeq \SS^\Z_{i,j}$ compatible with the product,
and the algebra $\SS$ is the corresponding folding of $\SS^\Z$.

The Koszul properties for $\SS$ (with the grading $\deg_K$) and for $\SS^\Z$ are equivalent.
On the other hand, it is easy to construct the isomorphism of $\Z$-algebras between $\SS^\Z$ and
the $\Z$-algebra corresponding to the $\Z$-graded algebra 
$$B_0\oplus B_1\oplus B_2=k\oplus V\oplus k,$$ 
where the multiplication $V\ot V=B_1\ot B_1\to B_2=k$ is given by a nondegenerate symmetric bilinear form.
It is well known that the algebra $B$ is Koszul.
Hence, $\SS$ is also Koszul.
\ed

\begin{rem} The algebra $\SS^!$ for $g=\id$ is closely related to the noncommutative projective line $\P^1_n$
(see \cite{Piont}, \cite{Minamoto}). Namely, it is a folded
version of the $\Z$-algebra of a natural helix in the derived category of $\P^1_n$.
\end{rem}

Let $\be_X,\be_Y\in K$ be the idempotents corresponding to the vertices $X$ and $Y$, respectively.

\begin{lem}\label{E!-lem} 
(i) Let $m\ge 1$. If $x\in \SS^!_m\be_X$ satisfies $x\a_i^*=0$ for some $i$ then $x=0$.

\noindent (ii) Assume that either $n\ge 3$, or $g\neq \la\cdot\id$ (for any $\la\in k^*$), or the characteristic of $k$ is $\neq 2$.
If $x^0\in \SS_2^!\be_X$, $x^1\in \SS_2^!\be_Y$ satisfy
$$\b_i^*x^0=-x^1\b_i^*, \ \ \a_i^*x^1=x^0\a_i^*$$
for all $i$, then $x^0=0$ and $x^1=0$.

\noindent (iii)
If a collection of elements $(x_i)$, where $x_i\in \SS_2^!\be_Y$, satisfies 
$$\sum_{i,j} a_{ij} x_j\b_i^*=0,$$
then there exists $y\in \SS^!_1\be_X$ such that $x_j=y\a_j^*$ for each $j$.
\end{lem}

\Pf . (i) It is easy to see that the question does not depend on $g$, so we can assume $g=\id$. In this case,
we need to check that the element $x_1$ in the algebra $k\lan x_1,\ldots,x_n\ran/(x_1^2+\ldots+x_n^2)$
is not a right zero divisor. But in fact, the latter algebra is a domain by \cite[Thm.\ 0.2]{Zhang}.

\noindent
(ii) It is easy to see that we can reformulate the question as follows. Given an $n$-dimensional vector space $V$
(in our case the space spanned by $(\b_i^*)$)
and elements $x^0\in V^\vee\ot V$, $x^1\in V\ot V^\vee$, such that 
\begin{equation}\label{gv-x0-x1-eq}
g(v)\ot x^0=-x^1\ot v\mod (\id\ot V+V\ot \id) \ \text{ for any } v\in V,
\end{equation}
\begin{equation}\label{Adg-x0-x1-eq}
v^*\ot \Ad(g^{-1})x^1=x^0\ot g^*(v^*) \mod (\id\ot V^\vee+V^\vee\ot \id) \ \text{ for any } v\in V^\vee,
\end{equation}
we should deduce that $x^0$ and $x^1$ are proportional to $\id$.

Assume first that $n\ge 3$. Then we claim that \eqref{gv-x0-x1-eq} alone implies that $x^0$ and $x^1$
are proportional to $\id$. Indeed, suppose we have 
$$g(v)\ot x^0+x^1\ot v=\id\ot A(v)+B(v)\ot\id\in V\ot V^\vee\ot V$$
for all $v$, for some operators $A, B\in\End(V)$.
Taking the contraction in the third tensor component with $v^*\in V^\vee$, we get the identity
$$g(v)\ot\lan x^0,v^*\ran+\lan v,v^*\ran x^1=\lan A(v),v^*\ran \id+B(v)\ot v^*\in V\ot V^\vee.$$
Thus, whenever $\lan v,v^*\ran=0$, the operator $\lan A(v),v^*\ran \id\in \End(V)$ is the sum
of two operators of rank $1$. Since $n\ge 3$, this implies that $\lan A(v),v^*\ran=0$.
Hence, $A(v)$ is proportional to $v$ for any $v\in V$, i.e., $A=\la\cdot\id$ for some $\la\in k$. A similar argument
shows that $B=\mu\cdot g$ for $\mu\in k$. Thus, subtracting some multiples of $\id$ from $x^0$ and $x^1$
we reduce ourselves to the situation when
$$g(v)\ot x^0+x^1\ot v=0 \ \text{ for any } v\in V.$$
Using contractions as above it is easy to deduce from this that $x^0=0$ and $x^1=0$.

Next, assume that $n=2$. Then we are going to rewrite condition \eqref{gv-x0-x1-eq} by fixing a symplectic isomorphism
$s:V\to V^\vee$ and observing that $(s\ot \id_V)(\bigwedge^2 V)$ and $(\id_V\ot s)(\bigwedge^2 V)$ are precisely the lines
spanned by the identity elements in $V^\vee\ot V$ and $V\ot V^\vee$. Thus, defining $y^0,y^1\in V\ot V$ by
$$x^0=(s\ot \id)(y^0), \ \ x^1=(\id\ot s)(y^1),$$
we see that \eqref{gv-x0-x1-eq} is equivalent to the equation
$$g(v)\cdot q^0=-q^1\cdot v \ \text{ in } S^3V,$$
where $q^i$ is the image of $y^i$ in $S^2V$.
Assume first that $g$ is not proportional to $\id$. Then the relation above implies that $q^0=q^1=0$ in $S^2V$.
Indeed, let us pick $v_0\neq 0$ such that $g(v_0)$ is not a multiple of $v_0$, and assume $q^i$ are nonzero.
Then we should have the following equations in the algebra $SV$:
$$q^0=v_0\cdot v', \ \ q^1=-g(v_0)\cdot v'$$
for some $v'\in V$, $v'\neq 0$. Then for any $v$ we should have
$$g(v)\cdot v_0=g(v_0)\cdot v \ \text{ in } S^2V.$$
But this is a contratiction as soon as $v$ is not proportional to $v_0$.

Finally, let us consider the case $n=2$ and $g=\la\cdot\id$. Then the above argument gives
$$q^1=-\la q^0.$$
Similarly, we can rewrite condition \eqref{Adg-x0-x1-eq} for $g=\la\cdot\id$ as 
$$v\ot y^1=\la y^0\ot v,$$
where $v^*=s(v)$. Thus, we get
$$q^1=\la q^0.$$
Since the characteristic is $\neq 2$, we deduce that $q^0=q^1=0$.

\noindent
(iii) This follows from the exactness of the direct summand of the Koszul complex,
$$\ldots\to \be_Y \SS_1^!\ot \SS_2^*\be_X\to \be_Y\SS_2^!\ot \SS_1^*\be_X\to \be_Y\SS_3^!\be_X\to 0$$
(which holds since the algebra $\SS^!$ is Koszul).
\ed

\subsection{Calculation of Hochschild cohomology and proof of Theorem \ref{ainf-moduli-thm}}\label{HH-sec}

For a $\Z$-graded algebra $A$, we use the following convention for the bigrading on its Hochschild cohomology
$HH(A)$. We denote by $CH^{s+t}(A)\{t\}$ the space of linear maps $A^{\ot s}\to A$ of degree $t$. The corresponding
bigraded piece in the Hochschild cohomology is denoted by $HH^{s+t}(A)\{t\}$ (then the upper grading is derived Morita
invariant).

Let us consider the algebra $\SS=\SS(k^n,g)$ as above.
We denote by $\SS\{m\}$ (resp., $\SS^!\{m\}$)
the graded components of $\SS$ (resp., $\SS^!$) with respect to the grading
given by $\deg(\a_i)=0$, $\deg(\b_i)=1$ (resp., $\deg(\a_i^*)=0$, $\deg(\b_i^*)=-1$). 
We are interested in Hochschild cohomology of $\SS$ viewed as a graded
algebra with this grading (note that this affects some signs).

\begin{thm}\label{HH-thm} For $m\in\Z$, one has
$$HH^m(\SS)\{<-m\}=0.$$
Assume in addition that 
either $n\ge 3$, or $g\neq \la\cdot\id$ (for any $\la\in k^*$), or the characteristic of $k$ is $\neq 2$.
Then 
$$HH^1(\SS)\{-1\}=0.$$ 
\end{thm}

\Pf .
Recall that to compute the Hochschild cohomology of a Koszul $K$-algebra $A$ (where $K$ is commutative semisimple)
we can use the Koszul resolution (see e.g.,\cite[Sec.\ 3]{VdB-nc-hom}). 
More precisely, we have a natural embedding 
$$(A^!_m)^*\hra A_1^{\ot m}$$
(here and below all tensor products are over $K$),
so that the image consists of the intersection of kernels of the partial multiplication
maps 
$$a_1\ot\ldots\ot a_m\mapsto a_1\ot \ldots \ot a_ia_{i+1}\ot\ldots \ot a_m.$$
The corresponding subcomplex
$$A\ot (A^!_\bullet)^*\ot A\sub A\ot T^\bullet(A_+)\ot A$$
in the standard bar-resolution of $A$ by free $A-A$-bimodules is still exact.
Thus, we get a resolution of the form
$$[\ldots \to A\ot (A^!_m)^*\ot A\rTo{d_m} A\ot (A^!_{m-1})^*\ot A\to\ldots \to A\ot (A^!_1)^*\ot A\to A\ot A]\to A.$$
Let $(v_i)$ be generators in $A_1$, $(v_i^*)$ the dual generators of $A_1^!$.
Then the differential is given by
$$d_m(r\ot\phi\ot s)=\sum_i rv_i\ot v_i^*\phi\ot s+(-1)^m\sum_i r\ot \phi v_i^*\ot v_is,$$
where we use the $A^!$-bimodule structure on $(A^!)^*$ given by
the operators dual to the left and right multiplication.

Assume now that $A$ has an additional grading $\deg$, induced by some $\Z$-grading on $A_1$.
If we are interested in the Hochschild cohomology of 
$A$ as a graded algebra with respect to this grading, at this point we need to be careful in
using the appropriate signs. Namely, to compute the Hochschild cohomology $HH^*(A)$ we 
apply the functor $\Hom_{A\ot A^{op}}(?,A)$ to the above resolution and use the identification
$$A^!_m\ot A\simeq \Hom((A^!_m)^*,A)\simeq \Hom_{A\ot A^{op}}(A\ot (A^!_m)^*\ot A,A)$$
under which an element $c\in \Hom((A^!_m)^*,A)$ corresponds to the composition of
$\id\ot c\ot \id$ with the multiplication $\mu$ in $A$. Now we have to use the Koszul sign convention: 
$$\mu(\id\ot c\ot\id)(r\ot\phi\ot s)=(-1)^{\deg(c)\deg(r)}r c(\phi) s.$$ 
Thus, we get the complex computing the Hochschild cohomology of $A$,
$$A\to A_1^!\ot A\to \ldots \to A_{m}^!\ot A\rTo{\de_m} A_{m+1}^!\ot A\to\ldots$$
with the differential
$$\de_m(\psi\ot s)=(-1)^{(\deg(\psi)+\deg(s))\deg(v_i)}\sum_i \psi v_i^*\ot v_is+(-1)^{m+1}\sum_i v_i^*\psi\ot sv_i.$$
Here we 
assume that the basis $(v_i)$ is homogeneous with respect to $\deg$.

We can apply this procedure in our case since $\SS$ is Koszul, with the generators $(\a_i,\b_i)$
(see Lemma \ref{Koszul-AA-lem}). 
We are interested in the components $HH^*(\SS)\{j\}$. 
As explained above, these spaces
can be computed as cohomology of the complex $(\SS^!_\bullet\ot \SS)\{j\}$ with respect to the differential
$$\de(\psi\ot s)=\sum_i (\psi \a_i^*\ot \a_is+(-1)^j\psi \b_i^*\ot \b_is) +
(-1)^{m+1}\sum_i (\a_i^*\psi\ot s\a_i+\b_i^*\psi\ot s\b_i),$$
where $\psi\in \SS^!_m$.

Since $\SS\{j\}=0$ for $j\neq 0,1$, we have
$$(\SS^!\ot \SS)\{j\}=\SS^!\{j\}\ot \SS\{0\} \oplus \SS^!\{j-1\}\ot \SS\{1\}.$$
Note also that because $\a_i^*$ and $\b_j^*$ have to alternate in any nonzero monomial in $\SS^!_m$, we have
$\SS^!_m\{j\}=0$ unless $m\in \{-2j-1,-2j,-2j+1\}$. 
This immediately implies the vanishing
$$HH^m(\SS)\{<-m-1\}=0$$
for any integer $m$.

For $m\ge 0$ the space $HH^m(\SS)\{-m-1\}$ is identified with the kernel of the map
$$\de:\SS^!_{2m+1}\{-m-1\}\ot \SS\{0\}\to\SS^!_{2m_2}\{-m-1\}\ot \SS\{0\}.$$
But $\SS^!_{2m+1}\{-m-1\}\ot \SS\{0\}=\SS^!_{2m+1}\be_X\ot \be_X$, and
for $x\in \SS^!_{2m+1}\be_X$ we have
$$\de(x\ot \be_X)=\sum_i x\a_i^*\ot \a_i.$$
Thus, Lemma \ref{E!-lem}(i) implies that this kernel is zero.


Next,  $HH^1(\SS)\{-1\}$ is the cohomology in the middle term of 
\begin{equation}\label{HH2-1seq}
\SS^!_1\{-1\}\ot \SS\{0\}\to \SS^!_2\{-1\}\ot \SS\{0\}\to \SS^!_3\{-1\}\ot \SS\{0\}\oplus \SS^!_3\{-2\}\ot \SS\{1\}.
\end{equation}
An element of $\SS^!_2\{-1\}\ot \SS\{0\}$ has form
$$x=x^0\ot \be_X+x^1\ot \be_Y+\sum_j x_j\ot \a_j,$$
where $x^0\in \SS_2^!\be_X$, $x^1\in \SS_2^!\be_Y$ and $x_i\in \SS_2^!\be_Y$.
We have
$$\de(x)=\sum_i(x^0\a_i^*\ot \a_i-x^1\b_i^*\ot \b_i)-
\sum_{i,j} a_{ij}x_j\b_i^*\ot \xi_0-
\sum_i(\b_i^*x^0\ot \b_i+\a_i^*x^1\ot \a_i).$$
Thus, $\de(x)=0$ if and only if
\begin{align*}
&\b_i^*x^0=-x^1\b_i^*, \ \ \a_i^*x^1=x^0\a_i^* \text{  for all } i,\\ 
&\sum_{i,j} a_{ij}x_j\b_i^*=0.
\end{align*}
The coboundaries come from elements $\SS^!_1\{-1\}\ot \SS\{0\}=\SS^!_1\be_X\ot \be_X$ and have form
$$\de(y\ot \be_X)=\sum_i y\a_i^*\ot \a_i.$$
Thus, Lemma \ref{E!-lem}(ii)(iii) implies that \eqref{HH2-1seq} is exact (under our assumptions), 
and hence, $HH^1(\SS)\{-1\}=0$.
\ed

\begin{rem} In \cite[Sec.\ (2c)]{Seidel-flux} the computation similar that of Theorem \ref{HH-thm} is done
in the case $n=2$, $g=\id$, $k=\C$. In this case one also has $HH^2(\SS)\{-1\}=0$, and $HH^2(\SS)\{-2\}$
can be identified with the space of binary quartic polynomials. As explained in \cite[Sec.\ (2f)]{Seidel-flux},
this means that all minimal $A_\infty$-structures on the algebra $\SS(\C^2,\id)$ are realized by double coverings
of $\P^1$.
\end{rem}

\noindent
{\it Proof of Theorem \ref{ainf-moduli-thm}}.
We apply \cite[Thm.\ 2.2.6]{P-more-pts} to the family $\SS_n$, $n\ge 3$ (resp., $\SS_{2,U}$).
More precisely we use the following vanishing of components of Hochschild cohomology for
any algebra $\SS=\SS(k^n,g)$ (implied by Theorem \ref{HH-thm}):
$$HH^{\le 1}(\SS)\{<0\}=HH^2(\SS)\{<-2\}=0.$$
\ed

\section{Pairs of $1$-spherical objects and noncommutative algebras}\label{sph-pairs-sec}

\subsection{Spherical objects and spherical twists}\label{spherical-sec}

Recall (see \cite{ST}, \cite[I.5]{Seidel-book})
that an object $X$ of a $k$-linear $A_\infty$-category $\CC$, where $k$ be a field, 
is called {\it $n$-spherical} if $\Hom^i(X,X)=0$ for $i\neq 0,n$,
$\Hom^0(X,X)=\Hom^n(X,X)=k$, and for any object $Y$ of $\CC$ the pairing between the morphism spaces in the cohomology category $H^*\CC$,
$$\Hom^{n-i}(Y,X)\ot\Hom^i(X,Y)\to \Hom^1(X,X),$$
induced by $m_2$, is perfect.

We need the following generalization of this notion to the case of an $S$-linear $A_\infty$-category, where $S$ is a 
commutative ring (our definition is not the most general possible: we impose
rather strong assumptions on the $\hom$-complexes). 

\begin{defi} Let $X$ be an object of an $S$-linear $A_\infty$-category $\CC$. Assume that for any $Y$ the complexes
$\hom(X,Y)$ and $\hom(Y,X)$ are bounded complexes of finitely generated projective $S$-modules.
Then $X$ is called {\it $n$-spherical} if $\Hom^i(X,X)=0$ for $i\neq 0,n$, $\Hom^0(X,X)=S\cdot\id_X$,
$\LL_X:=\Hom^n(X,X)$ is a locally free $S$-module of rank $1$,
and for any $Y$ in $\CC$ the following composed chain map of complexes of $S$-modules is a quasi-isomorphism:
\begin{equation}\label{spherical-composition-map-eq}
\hom(Y,X)\to \hom(X,Y)^\vee \ot_S \hom(X,X)\to \hom(X,Y)^\vee\ot \tau_{\ge n}\hom(X,X).
\end{equation}
Here the first arrow induced by $m_2$, while the second comes from the natural map $\hom(X,X)\to \tau_{\ge n}\hom(X,X)$,
where $\tau_{\ge n}$ is the truncation functor. Also, $P^\vee=\Hom_S(P,S)$ is the dual of a finitely generated
projective module $P$.
\end{defi}

Note that the complex $\tau_{\ge n}\hom(X,X)$ is bounded, has the only cohomology at the left-most term and
all of its subsequent terms are finitely generated projective $S$-modules. 
This implies that the left-most term is also finitely generated projective and there exists a homotopy equivalence
$$\tau_{\ge n}\hom(X,X)\to \Hom^n(X,X)[-n]=\LL_X[-n].$$
Fixing such an equivalence we can view the map \eqref{spherical-composition-map-eq} as
a chain map
\begin{equation}\label{n-spherical-map-eq}
s_Y:\hom(Y,X)\to \hom(X,Y)^\vee\ot\LL_X [-n]
\end{equation}

Let $\Tw(\CC)$ denote the category of twisted complexes over $\CC$ (see e.g., \cite[Sec.\ 7.6]{keller-ainf}).

\begin{lem}\label{spherical-Tw-com-lem}
An $n$-spherical object $X$ of $\CC$ remains $n$-spherical in $\Tw(\CC)$. 
\end{lem}

\Pf . First, note that for any twisted complex $Y$ we still have that $\hom(X,Y)$ and $\hom(Y,X)$ are bounded
complexes of finitely generated projective $S$-modules.
Now assume we are given an exact triangle $Y'\to Y\to Y''\to Y'[1]$ in $\Tw(\CC)$, such that the maps
\eqref{n-spherical-map-eq} for $Y'$ and $Y''$ are quasi-isomorphisms. Then
we have a morphism of exact triangles of $S$-modules
\begin{diagram}
\hom(Y'',X)&\rTo{}&\hom(Y,X)&\rTo{}&\hom(Y',X)&\to\ldots\\
\dTo{s_{Y''}}&&\dTo{s_Y}&&\dTo{s_{Y'}}\\
\hom(X,Y'')^\vee\ot\LL_X [-n]&\rTo{}&\hom(X,Y)^\vee\ot\LL_X [-n]&\rTo{}&\hom(X,Y')^\vee\ot\LL_X [-n]&\to\ldots
\end{diagram}
so the fact that $s_{Y'}$ and $s_{Y''}$ are quasi-isomorphisms imply that $s_Y$ is also a quasi-isomorphism.
Since every object of $\Tw(\CC)$ is an iterated extension of shifts of objects of $\CC$, the assertion follows.
\ed

Given an $n$-spherical object $E$, we can define (see \cite{ST}, \cite[I.5]{Seidel-book})
the twist and the adjoint twist $A_\infty$-functors
$$T_E, T'_E:\Tw(\CC)\to \Tw(\CC)$$
by 
$$T_E(X)=\Cone(\hom(E,X)\ot E\rTo{ev} X), \ T'_E(X)=\Cone(X\rTo{\ev'}\hom(X,E)^\vee\ot E)[-1].$$
The proof of \cite[Prop.\ 2.10]{ST} extends to our situation to prove that $T'_ET_E$ and $T_ET'_E$ are
isomorphic to identity in the homotopy category of functors from $\Tw(\CC)$ to itself. 


\subsection{$n$-pairs of $1$-spherical objects}\label{n-pairs-1-sph-sec}

$A_\infty$-structures we want to consider are related to the following pairs of objects in $A_\infty$-categories.

\begin{defi} Let $\CC$ be a $S$-linear $A_\infty$-category, where $S$ is a commutative ring.

\noindent
(i) We call a pair of $1$-spherical objects $(E,F)$ in $\CC$ an
{\it $n$-pair} if $\Hom^*(E,F)$ is concentrated in degree $0$ and
is isomorphic to $S^n$. In addition we require that $\Hom^1(F,F)\simeq S$.
Note that this implies that $\Hom^*(F,E)$ is a free $S$-module of rank $n$ concentrated in degree $1$.
We say that $(E,F)$ a {\it symmetric $n$-pair} if in addition $\Hom^1(E,E)\simeq S$ and
the two perfect pairings
\begin{equation}\label{two-pair-eq}
\begin{array}{l}
\Hom^1(F,E)\ot_S \Hom^0(E,F)\to \Hom^1(E,E) \ \text{ and} \\
\Hom^0(E,F)\ot_S \Hom^1(F,E)\to \Hom^1(F,F) 
\end{array}
\end{equation}
(coming from the conditions that $E$ and $F$ are $1$-spherical)
lead to two dualities $\Hom^1(F,E)\simeq \Hom^0(E,F)^\vee$ that differ by a scalar in $S^*$.

\noindent
(ii) A {\it weak $n$-pair} in $\CC$ is a pair of objects $(E,F)$ in $\CC$ such that $F$ is $1$-spherical
with $\Hom^1(F,F)\simeq S$,
$E$ satisfies $\Hom^0(E,E)\simeq S$, $\Hom^{\neq 0,1}(E,E)=0$, $\LL_E:=\Hom^1(E,E)$ is a locally free $S$-module 
of rank $1$, and $\Hom^*(E,F)=\Hom^0(E,F)\simeq S^n$.
An {\it enhanced weak $n$-pair} is a weak $n$-pair $(E,F)$ together with chosen isomorphisms
$\Hom^0(E,F)\simeq S^n$ and
$\Hom^1(F,F)\simeq S$.
\end{defi}

Note that for an enhanced weak $n$-pair, 
the second of the pairings \eqref{two-pair-eq} is perfect, so it defines an isomorphism
$\Hom^1(F,E)\simeq V^\vee$, where $V:=\Hom^0(E,F)$. Then the first of the pairings \eqref{two-pair-eq} has form 
\begin{equation}\label{g-def-eq}
\Hom^1(F,E)\ot_S \Hom^0(E,F)\simeq V^\vee\ot V\to \LL_E: (v^*,v)\mapsto \lan v^*,gv\ran
\end{equation}
for a unique $g\in \End_S(V)\ot\LL_E$.

\begin{lem}\label{invertible-g-lem} Let $(E,F)$ be a weak $n$-pair in $\CC$, and assume that $\CC$ is split generated
by $(E,F)$. Then $(E,F)$ is an $n$-pair (i.e., $E$ is $1$-spherical) in $\CC$ if and only if the element $g\in \End_S(V)\ot\LL_E$ defined by \eqref{g-def-eq} is invertible.
\end{lem}

\Pf . Let $\CC'\sub \CC$ be the triangulated envelope of $(E,F)$. It is clear that $E$ is $1$-spherical in $\CC$ if and only if it is $1$-spherical in $\CC'$. 
But by Lemma \ref{spherical-Tw-com-lem}, $E$ is spherical in $\CC'$ if and only if the pairing
\eqref{g-def-eq} is perfect, which is equivalent to $g$ being invertible.
\ed

\begin{ex} Let $D^b(C)$ be the (enhanced) derived category of coherent sheaves on an elliptic curve $C$ over an
algebraically closed field $k$.
Then $1$-spherical objects in $D^b(C)$ are (up to shift) either simple vector bundles or the skyscraper
sheaves $\OO_p$. The group of autoequivalences of $D^b(C)$ acts transitively on them, so any $n$-pair of $1$-spherical
objects can be transformed by an autoequivalence into a pair $(E,\OO_p)$, where $E$ is an simple vector bundle of rank $n$. 
It is easy to see that any such $n$-pair is special.
\end{ex}

Given $g\in\PGL_n(S)$ and a minimal $A_\infty$-structure on $\SS=\SS(S^n,g)$ we get an $n$-pair of $1$-spherical objects
in the corresponding full subcategory $\{\be_X\cdot\SS, \be_Y\cdot\SS\}$ of right $A_\infty$-modules over $\SS$.

Conversely, starting with an enhanced $n$-pair $(E,F)$ in an $S$-linear $A_\infty$-category, let us consider the full
$A_\infty$-subcategory with the objects $E$ and $F$.
Due to our assumption that $\hom(E,E)$, $\hom(E,F)$, $\hom(F,F)$ and $\hom(F,E)$ are bounded complexes
of projective modules, we can apply the homological perturbation to get an equivalent minimal $A_\infty$-structure
on this subcategory. Furthermore, as was explained before, we
can identify $\Hom^1(F,E)$ with $V^\vee$, where $V=\Hom^0(E,F)\simeq S^n$,
so that the second of the compositions \eqref{two-pair-eq} becomes the standard pairing between $V$ and $V^\vee$, while
the first has the form \eqref{g-def-eq} for some $g\in \GL_n(S)\ot\LL_E$. 
Thus, we get a minimal $A_\infty$-structure on $\SS(S^n,g)$.


\subsection{Some properties of the algebra $\EE(S^n,g)$}
\label{algebra-E-g-sec}


Let $S$ be a commutative ring, $\LL$ a locally free $S$-module of rank $1$.
For an element $g\in \Mat_n(S)\ot \LL$ we consider the algebra $\EE(S^n,g)$ defined by \eqref{E-V-eq}.

\begin{lem}\label{E-deg-1-lem} 
Assume that $n\ge 2$ and $g\in \Mat_n(S)\ot \LL$ is such that there exists $h\in\Mat_n(S)\ot\LL^{-1}$ with $\tr(gh)=1$.

\noindent
(i) Assume that either $n\ge 3$ or $n=2$ and $\tr(g^2h')=1$ for some $h'\in\Mat_2(S)\ot\LL^{-2}$.
Then the algebra $\EE(S^n,g)$ is generated over $S$ by degree $1$ elements. 

\noindent
(ii) $\EE(S^n,g)$ is generated by degree $1$ and degree $2$ elements.

\noindent
(iii) Assume that $g$ is invertible. Then the algebra $\EE(S^n,g)$ is Koszul. 
\end{lem}

\Pf . (i) Recall that $\EE(S^n,g)_1$ is the subspace $\End_g(S^n)$ of elements $a\in \Mat_n(S)$ such that $\tr(ga)=0$. 
The existence of $h$ such that $\tr(gh)=1$ implies that $\End_g(S^n)$ is a direct summand
in $\Mat_n(S)$.

We have to prove the surjectivity of the map
\begin{equation}\label{End-g-prod-map}
\End_g(S^n)\ot_S \End_g(S^n)\to \Mat_n(S)
\end{equation}
induced by the product in $\Mat_n(S)$.

First, we claim that it is enough to prove the assertion in the case when $S$ is a field. 
We can easily reduce to the case when $S$ is local. Now let $M$
denote the cokernel of the product map \eqref{End-g-prod-map}. Note that the construction of $M$ is compatible
with any change of scalars $S\to S'$. Thus, the case of the field implies that $M/\mg M=0$, where
$\mg\sub S$ is a maximal ideal. By Nakayama lemma, this gives that $M=0$.

Thus, it is enough to consider the case when $S=k$, where $k$ is a field. In this case we will prove a more general statement
that for any pair of nonzero elements $g_1,g_2$ one has
$$\End_{g_1}(k^n)\cdot \End_{g_2}(k^n)=\Mat_n(k),$$
where in the case $n=2$ we additionally require that $g_2g_1\neq 0$.
Since the question is that certain vectors generate $\Mat_n(k)$ as a linear space, without loss of generality we
can assume $k$ to be algebraically closed.

Note that for any $a,b\in \GL_n(k)$ one has
$$a\cdot\End_{g_1}(k^n)\cdot \End_{g_2}(k^n)\cdot b=\End_{g_1a^{-1}}(k^n)\cdot \End_{b^{-1}g_2}(k^n).$$
Thus, we can replace $g_1$ by $g_1a^{-1}$ and $g_2$ by $b^{-1}g_2$.
In the case when $g_1$ and $g_2$ are invertible this reduces the statement to the case $g_1=g_2=1$, which is
easy to check.

Next, we observe that for any nonzero $g\in \Mat_n(k)$ there exists $a\in\GL_n(k)$ such that $\tr(ag)=0$.
Indeed, otherwise, the entire hyperplane $\End_g(k^n)$ would be contained in the irreducible hypersurface $\det(a)=0$.
Thus, we can assume that $\tr(g_1)=\tr(g_2)=0$.
In this case we have $1\in \End_{g_1}(k^n)$ and $1\in \End_{g_1}(k^n)$. Hence,
$$\End_{g_1}(k^n)+\End_{g_2}(k^n)\sub \End_{g_1}(k^n)\cdot \End_{g_2}(k^n).$$
Thus, the only case when this subspace is not the entire $\Mat_n(k)$ is when
$\End_{g_1}(k^n)=\End_{g_2}(k^n)$, i.e., $g_1$ and $g_2$ are proportional. 
In this case we get that $\End_{g_1}(k^n)$ is a subalgebra.
As was observed above, we can assume that $g_1$ is degenerate. Let us choose a basis in which the last row of $g_1$
vanishes. Let $g'$ denote the $(n-1)\times (n-1)$-submatrix of $g_1$ obtained by deleting the last row and last column.

Assume first that $g'=0$. Then $\End_{g_1}(k^n)$ contains the maximal parabolic subalgebra of endomorphisms preserving
the hyperplane spanned by the first $n-1$ basis vectors. In the case $n\ge 3$ this implies that $\End_{g_1}(k^n)$
cannot be a subalgebra, since it would be strictly bigger than the maximal parabolic subalgebra. In the case $n=2$
if $g'=0$ then $g_1^2=0$ which contradicts the assumption that $g_2g_1\neq 0$ (recall that $g_1$ and $g_2$ are proportional)

Next, consider the case $g'\neq 0$. Let $e_{ij}$ denote the standard matrices with $1$ as the $(i,j)$-entry. 
Then for every $i$ we have $e_{in}\in \End_{g_1}(k^n)$, and for every $j\le n-1$ there exists a matrix $A_j$ with
zero last row such that $e_{nj}+A_j\in \End_{g_1}(k^n)$. But then 
$e_{ij}=e_{in}\cdot (e_{nj}+A_j)$, so we deduce that $\End_{g_1}(k^n)\cdot\End_{g_1}(k^n)$ is $\Mat_n(k)$.

\noindent
(ii) We have to show that the product map
$$\End_g(S^n)\ot_S \Mat_n(S)\to \Mat_n(S)$$
is surjective. As before, it is enough to consider the case when $S=k$ is a field.
Furthermore, as in part (i) we reduce to the case when $\tr(g)=0$, so that $1\in \End_g(S^n)$,
in which case the assertion is clear.

\noindent
(iii) It is easy to check that the algebra $\EE(S^n,g)$ is quadratic dual to the second Veronese subalgebra of
the algebra $\SS(S^n,g)$, corresponding to the vertex $X$ (i.e., one considers paths of even length starting from $X$).
Since the latter algebra is Koszul by Lemma \ref{Koszul-AA-lem}, the result follows (see \cite[Prop.\ 2.2(i)]{PP-book}).
\ed

In the next result we consider derivations of $\EE(S^n,g)$ as an ungraded algebra 
(i.e., there is no Koszul sign in the Leibnitz rule). This result will be used later in studying
automorphisms of a filtered algebra whose associated graded algebra is isomorphic to $\EE(S^n,g)$.

\begin{prop}\label{derivations-E-prop} 
Assume that one of the following two conditions holds:
\begin{itemize}
\item $n\ge 3$ and $g$ is invertible;
\item $\tr(g)$ is a generator of $\LL$, and there exists $h_1\in \GL_n(S)$ with $\tr(gh_1)=0$. 
\end{itemize}
Then any derivation $\EE(S^n,g)\to \EE(S^n,g)$ of degree $m\le -1$ is zero.
\end{prop}

\Pf . 
The problem is local, so we can assume that $\LL=S$.

Let us set for brevity $\EE:=\EE(S^n,g)$, $\EE':=\End(V)[z,z^{-1}]$. By the assumption, we can fix an element
$h_1\in \EE_1$ such that the operators $\EE'_i\to \EE'_{i+1}$ of left and right multiplication by $h_1$ are invertible for $i\in\Z$.
Now assume we have a derivation $D:\EE\to \EE$ of degree $m\le -1$.
First, we are going to extend $D$ to a derivation $D':\EE'_{\ge 1}\to \EE'$ of degree $m$.
For this we compose $D$ with the embedding $\EE\hra \EE'$ and then set for $x\in \EE'_{\ge 1}$,
$$D'(x)=D(xh_1)h_1^{-1}-xD(h_1)h_1^{-1},$$
where we use the operation of multiplication by $h_1^{-1}$ as a degree $-1$ map 
$\EE'_{\ge 1}\to \EE'$.
Note that the expression in the right-hand side is well-defined since $h_1\in \EE_1$ and $xh_1\in \EE'_{\ge 2}=\EE_{\ge 2}$.
Also we have $D'=D$ on $\EE$. Before checking that $D'$ is a derivation we observe that for $x\in \EE'_i$, $i\ge 1$, one has
$$D'(x)=D''(x)=h_1^{-1}D(h_1x)-h_1^{-1}D(h_1)x.$$
Indeed, this can be checked by applying the Leibnitz identity to write $D(h_1xh_1)$ in two ways (note that $h_1x\in\EE$ and 
$xh_1\in \EE$).
Now we can prove the Leibnitz identity for $D'=D''$. Namely, it is enough to check that for $x_1,x_2\in \EE'_{\ge 1}$ one has
$$D(x_1x_2)=D''(x_1)x_2+x_1D'(x_2).$$
It is easy to see that this is equivalent to the identity
$$D(h_1x_1)x_2h_1+h_1x_1D(x_2h_1)=D(h_1)x_1x_2h_1+h_1D(x_1x_2)h_1+h_1x_1x_2D(h_1),$$
obtained by writing $D(h_1x_1x_2h_1)$ in two ways.

Next, we claim that there exists $a\in\End(V)$ and $s\in S$ such that 
$$D'(x)=\ad(az^m)+s\cdot z^{m+1}\frac{d}{dz}.$$ 
Indeed, as above, we can extend $D'$ to a derivation on $\End(V)[z,z^{-1}]$ of degree $m$. Using Morita equivalence
of the latter ring with $S[z,z^{-1}]$ we get the result. 

Let us first assume that $m=-1$.
Our goal is to show that $s=0$ and $a$ is proportional to the identity. 
To this end we investigate the condition 
$$D'(\End_g(V)z)\sub S\cdot\id,$$
which means that for any $x\in \End_g(V)$ one has 
$$(a+s)x-xa\in S\cdot\id,$$
Equivalently, for any $y\in \End(V)$ with $\tr(y)=0$ and any $x\in\End_g(V)$, one has
$$\tr((a+s)xy-xay)=\tr(x[y(a+s)-ay])=0.$$
Equivalently, we should have
\begin{equation}\label{a-comm-eq1}
y(a+s)-ay\in S\cdot g \ \ \text{whenever } \tr(y)=0.
\end{equation}

Assume first $n\ge 3$ and $g$ is invertible. Then substituting $y=e_{i_0j_0}$ with $i_0\neq j_0$ in \eqref{a-comm-eq1}
we get an equality of the form 
$$e_{i_0j_0}(a+s)-ae_{i_0j_0}=\la\cdot g.$$
We claim that this is possible only when $\la=0$. Indeed, for every $j\neq j_0$
we get
$$\la\cdot e_j=\mu_j\cdot g^{-1}e_{i_0}$$
for some $\mu_j\in S$. Let $g^{-1}e_{i_0}=\sum b_ie_i$.
Then we have a system of equations
$$\mu_jb_j=\la, \ \ \mu_jb_i=0 \ \text{ for } i\neq j, j\neq j_0.$$
Since $g\cdot g^{-1}=\id$, we have some $(a_i)$ in $S$ with $\sum b_ia_i=1$.
Thus, we deduce
$$\mu_j=\la a_j.$$
Plugging back in the above equation, we get that $\la\cdot I=0$, where
$I\sub S$ is the ideal generated by $(a_jb_j-1)_{j\neq j_0}$ and $(a_jb_i)_{i\neq j, j\neq j_0}$.
Now choosing a pair $i\neq j$ in $[1,n]\setminus\{j_0\}$ (this is possible since $n\ge 3$), we obtain
$$I\supset (a_ib_i-1,a_jb_j-1,a_jb_i)=(1),$$
and hence, $\la=0$.

Thus, we derive that for every $i\neq j$ one has
$$e_{ij}(a+s)-ae_{ij}=0$$
This immediately implies that $a$ is diagonal, and the diagonal entries $(a_{ii})$
satisfy $a_{jj}-a_{ii}=s$ for $i\neq j$. Using again the assumption $n\ge 3$, we obtain 
that $s=0$ and $a$ is proportional to $\id$.

Next, let us consider the case when $\tr(g)$ is invertible (but $g$ is not necessarily invertible). 
Considering traces of both sides of \eqref{a-comm-eq1} we derive that
$$y(a+s)-ay=0 \ \ \text{whenever } \tr(y)=0.$$
Now substituting $y=e_{ij}$ for $i\neq j$ one can easily derive that $a$ has to be a diagonal matrix.
As we have above, in the case $n\ge 3$ this implies that $s=0$ and $a$ is proportional to $\id$.
If $n=2$ then taking into account the equation for the diagonal $y$ with entries $1$ and $-1$, we get the same conclusion.

In the case $m\le-2$ we should have
$$D'(\End_g(V)z)=0,$$
i.e., $(a+s)x-xa=0$ for any $x\in\End_g(V)$. As we have seen above, this implies that
$a=s=0$.
\ed

\begin{rem} One can also check that there are no nonzero derivations $\EE(k^2,g)\to \EE(k^2,g)$ of degree $-1$ 
provided $k$ is a field of characteristic $\neq 2$ and $g$ is invertible.
\end{rem}

\begin{ex}
Assume that $n=2$ and $2=0$ in $S$, and let us take $g=\id$.
Then there exist nontrivial derivations of $\EE(S^2,\id)$ of degree $-1$. More precisely,
for any $2\times 2$ matrix $a$ the derivation $\ad(az^{-1})+\tr(a)\frac{d}{dz}$ of $\Mat_2(S)[z,z^{-1}]$
restricts to a derivation of $\EE(S^2,\id)$. This essentially amounts to the identity 
$$[y,a]=\tr(a)y+\tr(ay)\id$$
for $2\times 2$-matrices $a$ and $y$ such that $\tr(y)=0$ (it only holds because $2=0$ in $S$).
\end{ex}

\subsection{The algebra associated with a pair of $1$-spherical objects}\label{sph-pair-to-algebra-sec}

Let $(E,F)$ be a weak $n$-pair in a minimal $S$-linear $A_\infty$-category with $\Hom^0(E,F)=V\simeq S^n$, equipped with the trivialization 
$$\Hom^1(F,F)\simeq S.$$
We denote by $\xi_F\in\Hom^1(F,F)$ the corresponding generator.
 
We are going to associate to these data an $S$-algebra
with an increasing exhaustive filtration $(F_nA)$, together with an isomorphism of graded $S$-algebras
\begin{equation}\label{gr-E-eq}
\gr^F A=\bigoplus_{n\ge 0}F_nA/F_{n-1}A\simeq \EE(V,g)^{op},
\end{equation}
where $g\in\End(V)\ot \LL_E$ is defined as before. Namely, we use the second of the pairings \eqref{two-pair-eq}
to identify $\Hom^1(F,E)$ with $V^\vee$, and define $g$ so that the first pairing takes the form \eqref{g-def-eq}.

We will see also that (under some mild assumptions) the filtered algebra $(A,F_\bullet A)$ 
is determined by the following higher products:
\begin{equation}\label{higher-products-main-eq}
\begin{split}
m_3: & V\ot V^\vee\ot V\to V\simeq\Hom^0(E,F)\ot\Hom^1(F,E)\ot\Hom^0(E,F)\to\Hom^0(E,F)=V,\\
m_3: & V^\vee\ot V\ot V^\vee\simeq\Hom^1(F,E)\ot\Hom^0(E,F)\ot\Hom^1(F,E)\to\Hom^1(F,E)\simeq V^\vee,\\
m_3: & V^\vee\ot V\simeq\Hom^1(F,E)\ot\Hom^1(F,F)\ot\Hom^0(E,F)\to \Hom^1(E,E)=\LL_E,\\
m_4: & V^\vee\ot V\ot V^\vee\ot V\simeq\\
&\Hom^1(F,E)\ot\Hom^0(E,F)\ot\Hom^1(F,E)\ot\Hom^0(E,F)\to\Hom^0(E,E)=S.
\end{split}
\end{equation}
Let us define the maps $r, r':\End(V)\to\End(V)$ by 
\begin{equation}\label{r-r'-def-eq}
r(v\ot v^*)=\sum_i m_3(e_i,v^*,v)\ot e_i^*, \ \ r'(v\ot v^*)=e_i\ot \sum_i m_3(v^*,v,e_i^*),
\end{equation}
for $v\in V$, $v^*\in V^\vee$, where $(e_i)$ and $(e_i^*)$ are dual bases of $V$ and $V^\vee$.
Similarly, we define $s:\End(V)\to \LL_E$ by
\begin{equation}\label{s-def-eq}
s(v\ot v^*)=m_3(v^*,\xi_F,v).
\end{equation}

\begin{thm}\label{spherical-filtered-thm} 
Let $(E,F)$ be a weak $n$-pair in a minimal $S$-linear $A_\infty$-category $\CC$, such that
$\Hom^0(E,F)=V\simeq S^n$, where $n\ge 2$. Let us fix a trivialization $\Hom^1(F,F)\simeq S$, 
and let $g\in \End(V)\ot\LL_E$ be the corresponding element defined using pairings \eqref{two-pair-eq}. 
Assume in addition that there exists $h\in \End(V)\ot\LL_E^{-1}$ such that $\tr(gh)=1$.
Set 
$$E_i=T^i(E)\in \Tw(\CC),$$ 
where $T=T_F$ is the spherical
twist with respect to $F$. Let us consider the graded associative algebra 
\begin{equation}\label{R-algebra-eq}
\RR=\RR_{T,E}:=\bigoplus_{n\ge 0}\Hom(E_0,E_n),
\end{equation}
with the product $ab=T^i(a)\circ b$, where $b\in \Hom(E_0,E_i)$, $a\in \Hom(E_0,E_j)$.
Then 

\noindent
(i) $\RR$ is canonically isomorphic to the Rees algebra of a filtered algebra $(A,F_\bullet A)$ equipped
with an isomorphism $\gr^F(A)\simeq \EE(V,g)^{op}\simeq \EE(V^\vee,g^*)$.
In addition, $\Hom^{\neq 0}(E_0,E_n)=0$ for $n>0$.

\noindent
(ii) There exist embeddings $\End_g(V)\hra \RR_1$ and $\End(V)\hra \RR_2$, such that
$$\RR_1=\End_g(V)\oplus S\cdot t, \ \ \RR_2=\End(V)\oplus \End_g(V)\cdot t\oplus S\cdot t^2,$$
where $t$ is the central element of degree $1$
corresponding to the isomorphism with the Rees algebra. With respect to these decompositions, 
for $a,b\in\End_g(V)\sub\RR_1$ one has
\begin{equation}\label{R-product-formula-eq}
a\cdot b=ba+[r(b)a+br'(a)+s(ba)h]t+m_4(a\ot b)t^2.
\end{equation}
Here we use the higher products \eqref{higher-products-main-eq} and the corresponding maps
$r,r',s$ (see \eqref{r-r'-def-eq}, \eqref{s-def-eq}).
\end{thm}

\Pf . It will be notationally convenient for a while not to use the trivialization of $\Hom^1(F,F)$,
so let us set $L:=\Ext^1(F,F)$. 

For an $S$-module $M$, we set $ML^i:=M\ot L^{\ot i}$. 
Note that the second of the pairings \eqref{two-pair-eq} induces
an identification $\Hom^1(F,E)\simeq V^\vee L$. 

\medskip

\noindent
{\bf Step 1}.
We start by finding explicit twisted complexes representing $E_i$.
Namely, let us denote by $E_i$, for $i\ge 1$, the following twisted complex:
\begin{equation}\label{E-i-simple-complex-eq}
\Hom^1(F,E)L^{i-1}\ot F\rTo{\de_i} \Hom^1(F,E)L^{i-2}\ot F\rTo{\de_{i-1}}\ldots 
\rTo{\de_2}\Hom^1(F,E)\ot F\rTo{\de_1} E.
\end{equation}
Here the differentials $\de_i$ with $i>1$ are induced by the evaluation maps
$L\ot F\to F[1]$, while the differential $\de_1:\Hom^1(F,E)\ot F\to E[1]$ is
also the evaluation map. 

We are going to construct the homotopy equivalences $T(E_i)\simeq E_{i+1}$.
Note that for $i=0$ we have $T(E_0)=T(E)=E_1$ by the definition of the twist functor $T=T_F$.
The complex $\hom(F,E_i)$ has form
\begin{diagram}
\Hom^1(F,E)L^{i-1}\ot\id_F&\ \ &\Hom^1(F,E)L^{i-1}\ot\id_F&\ \ &\cdots&\ \ &\Hom^1(F,E)\ot\id_F\\
&\rdTo{}&&\rdTo{}&&\rdTo{}&&\rdTo{}\\
\Hom^1(F,E)L^{i} &\ \ &\Hom^1(F,E)L^{i-1}&\ \ &\cdots&\ \ &\Hom^1(F,E)L&\ \ &   \Hom^1(F,E)
\end{diagram}
with the first row in degree $0$ and the second row in degree $1$ (note that higher products do not appear since
we assume our $A_\infty$-structures to be strictly unital). Since all the components of the differential
are isomorphisms, the natural embedding and the projection 
give a homotopy equivalence 
$$\hom(F,E_i)\simeq \Hom^1(F,E)L^{i}[-1].$$
Hence, we deduce a homotopy equivalence 
\begin{equation}\label{T-E-E-eq}
E_{i+1}=\Cone(\Hom^1(F,E)L^{i}\ot F[-1]\rTo{\de_{i+1}} E_i)\simeq \Cone(\hom(F,E_i)\ot F\rTo{\ev} E_i)=T(E_i)
\end{equation}
as claimed.

For what follows we need to know explicitly the maps between $E_{i+1}$ and $T(E_i)$.
Note that the embedding of $\Hom^1(F,E)L^{i}\ot F[-1]$ into $\hom(F,E_i)\ot F$ commutes with the maps 
to $E_i$ used to form the above cones, however, the projection in the other direction only commutes up to homotopy. Namely,
we have a homotopy $h$ between the map $\ev:\hom(F,E_i)\ot F\to E_i$ and the composition
$$\hom(F,E_i)\ot F\to \Hom^1(F,E)L^{i}\ot F\rTo{\de_{i+1}} E_i,$$
with the nonzero components 
$$h_j: \Hom^1(F,E)L^{j}\ot F\rTo{\id} \Hom^1(F,E)L^{j}\ot F\hra E_i,$$
for $0\le j\le i-1$. Hence, the map $E_{i+1}\to T(E_i)$ is given by the obvious embedding of complexes, while the
map $T(E_i)\to E_{i+1}$ is the identity on the summands $E_i$ and all the summands $\Hom^1(F,E)L^{j}\ot F$ of 
$\hom^1(F,E_i)\ot F$.

\medskip

\noindent
{\bf Step 2}.
The complex $\hom(E_0,E_i)=\hom(E,E_i)$ has form
$$\left(\bigoplus_{j=0}^{i-1}  \Hom^1(F,E)L^{j}\ot \Hom^0(E,F)\right)\oplus \Hom^0(E,E)\to \Hom^1(E,E),$$
with the differential given by $d(\id_E)=0$, 
$$d(e\ot\xi^{\ot j}\ot x)=m_{j+2}(e,\xi,\ldots,\xi,x).$$

Recall that the map $m_2:\Hom^1(F,E)\ot\Hom(E,F)\to \Hom^1(E,E)=\LL_E$
can be identified with the map 
$$(V^\vee\ot L)\ot V=\End(V)\ot L\to \LL_E:a\mapsto \tr(ga)$$
(note that canonically $g$ is an element of $\End(V)\ot L^{-1}\ot \LL_E$).
Thus, we immediately see that for $i\ge 1$ one has $\Hom^{\neq 0}(E_0,E_i)=0$, while
$\Hom^0(E_0,E_i)$ fits into an exact sequence
$$0\to \End_g(V)L\oplus \Hom^0(E,E)\to \Hom^0(E_0,E_i)\to
\left(\bigoplus_{j=2}^i \End(V)L^{j} \right) \to 0,$$
where we use the identification  $\Hom^1(F,E)\simeq V^\vee L$.

\medskip

\noindent
{\bf Step 3}. 
For $0\le i<j$ let us consider the map of complexes
$$\hom(E_i,E_j)\rTo{T} \hom(T(E_i),T(E_j))\to \hom(E_{i+1},E_{j+1}),$$
where the second arrow is induced by the maps $E_{i+1}\to T(E_i)$ and
$T(E_j)\to E_{j+1}$ described in Step 1. Then in the case $i>0$ the following square is commutative
\begin{diagram}
\hom^0(E_i,E_j)&\rTo{}& \Hom^0(\Hom^1(F,E)L^{i-1}\ot F,\Hom^1(F,E)L^{j-1}\ot F) \\
\dTo{T}&&\dTo{\ot \id_L}\\
\hom^0(E_{i+1},E_{j+1})&\rTo{}&
\Hom^0\bigl(\Hom^1(F,E)L^{i}\ot F,\Hom^1(F,E)L^{j}\ot F\bigr)
\end{diagram}
while in the case $i=0$ the following square is commutative
\begin{diagram}
\hom^0(E_0,E_j)&\rTo{}& \Hom^0(E,\Hom^1(F,E)L^{j-1}\ot F)\\
\dTo{T}&&\dTo{}\\
\hom^0(E_{1},E_{j+1})&\rTo{}&
\Hom^0\bigl(\Hom^1(F,E)\ot F,\Hom^1(F,E)L^{j}\ot F\bigr)
\end{diagram}
where the horizontal arrows are the natural projections, while the right vertical arrow in the second diagram
sends $a\in\Hom^1(F,E)L^{j-1}\ot V\simeq \End(V)L^j$ to 
$a^*\ot\id_F\in \End(V^\vee)L^j\ot\Hom^0(F,F)$.

\medskip

\noindent
{\bf Step 4}. For each $i\ge 1$, let us consider the natural projection (see Step 2)
$$\pi_i:\Hom^0(E_0,E_i)\to V^\vee L^i \ot\Hom(E,F)\simeq \End(V)L^i.$$
Note that for $i\ge 2$ it is surjective, while for $i=1$ its image is $\End_g(V)L$.
We claim that the map $\pi=(\pi_i)$ is a homomorphism of graded algebras
$$\RR\to \EE(V,g)^{op}.$$
To prove this let us consider elements $a\in \Hom^0(E_0,E_i)$ and $b\in \Hom(E_0,E_j)$, where $i>0$, $j>0$,
and set $\ov{a}=\pi_i(a)$, $\ov{b}=\pi_j(b)$.
Iterating Step 3 we see the component of $T^j(a)\in\hom^0(E_j,E_{i+j})$ in
$\Hom^0\bigl(\Hom^1(F,E)L^{j-1}\ot F,\Hom^1(F,E)L^{i+j-1}\ot F\bigr)$ is $\ov{a}^*\ot\id_F$.
It is easy to see that the $pi_{i+j}(T^j(a)\circ b)$ is obtained by composing the above component of $T^j(a)$
with $\ov{b}\in \End(V)L^j\simeq \Hom^0(E,\Hom^1(F,E)L^{j-1}\ot F)$.
Thus, if we view $\ov{b}$ as an element of $V^\vee\ot V L^j$ then we get
$$\pi_{i+j}(T^j(a)\circ b)=(\ov{a}^*\ot\id_V)(\ov{b})=\ov{b}\ov{a}\in \End(V)L^{i+j}.$$

\medskip

\noindent
{\bf Step 5}.
Let us define $t\in\Hom^0(E_0,E_1)$ to be the element represented by the element $\id_E\in\Hom^0(E,E)$ in the
complex $\hom(E,E_1)$, so that we have a decomposition
\begin{equation}\label{Hom-E-1-ex-eq}
\Hom(E_0,E_1)\simeq \End_g(V)L\oplus S\cdot t.
\end{equation}

We claim that for $i\ge 0$ the element $T^i(t)\in\Hom(E_i,E_{i+1})$ is represented
by the closed map of twisted complexes
\begin{diagram}
&&V^\vee L^{i}\ot F&\rTo{\de}& V^\vee L^{i-1}\ot F&\rTo{\de}&\ldots& \rTo{\de}& V^\vee L\ot F&\rTo{\de}& E\\
&&\dTo{\id}&&\dTo{\id}&&&&\dTo{\id}&&\dTo{\id}\\
V^\vee L^{i+1}\ot F&\rTo{\de}& V^\vee L^{i}\ot F&\rTo{\de}&
V^\vee L^{i-1}\ot F&\rTo{\de}&\ldots &\rTo{\de}& V^\vee L\ot F&\rTo{\de}& E
\end{diagram}
Indeed, this can be easily checked by induction by computing the effect of the twist functor $T$ on the above map
and taking into account the homotopy equivalence \eqref{T-E-E-eq}.

\medskip

\noindent
{\bf Step 6}.
It follows easily from Step 5 that the map of the left multiplication by $t$ in our algebra,
$$\Hom(E_0,E_i)\rTo{t\cdot}\Hom(E_0,E_{i+1})$$ 
is injective and its image is given by the classes that have
representatives in $\hom(E_0,E_{i+1})$ with zero component in 
$$V^\vee L^{i+1} \ot\Hom(E,F)\simeq V^\vee\ot V L^{i+1}\simeq \End(V)L^{i+1}.$$
Thus, for $i\ge 1$ we get an exact sequence
\begin{equation}\label{Hom-E-i-ex-eq}
0\to \Hom(E_0,E_{i-1})\rTo{t\cdot} \Hom(E_0,E_i)\rTo{\pi_i} \EE(V,g)_i L^i\to 0.
\end{equation}
From these exact sequences we deduce that the algebra $\RR$ is generated by degree $1$ and degree $2$ elements. Indeed,
since $\pi$ is a homomorphism,
this follows from the similar property of $\EE(V,g)^{op}$ (see Lemma \ref{E-deg-1-lem}).

\medskip

\noindent
{\bf Step 7}.
Next, let us consider an element $a\in \End_g(V)L$, and let us view it as a cochain in the term 
$V^\vee L\ot \Hom(E,F)\simeq \End(V)L$ of the complex $\hom(E_0,E_1)$.
We would like to calculate $T(a)$.
By definition, $T(a)$ is obtained as the composition
$$E_1=\Cone(\hom(F,E_0)\ot F\rTo{\ev} E_0)\to \Cone(\hom(F,E_1)\ot F\rTo{\ev} E_1)\to E_2,
$$
where the arrow between the cones is induced by $a:E_0\to E_1$, and the last arrow is the projection $T(E_1)\to E_2$
that we computed before.
Now the map $\hom(F,E_0)\rTo{a}\hom(F,E_1)$ has two nonzero components:
the map $\Ext^1(F,E)\to V^\vee L^2$ induced by by the composition with $a$ and the map
\begin{equation}\label{mu-a-eq}
\mu_a:\Ext^1(F,E)\to\Ext^1(F,E):x\mapsto m_3(\de_1,a,x).
\end{equation}
It follows that $T(a)$ is represented by the map
\begin{equation}\label{T-a-eq}
\begin{diagram}
V^\vee L\ot F&\rTo{\de_1}&E\\
\dTo{a^*\ot\id_F}&\rdTo{}&\dTo{a}\\
V^\vee L^2\ot F&\rTo{}&V^\vee L\ot F&\rTo{}&E
\end{diagram}
\end{equation}
where the diagonal arrow is $\mu_a\ot\id_F$.

\medskip

\noindent
{\bf Step 8}. Now let us check that $t$ is central in $\RR$. By Step 6, it is enough to check $t$ commutes
with elements of degree $1$ and $2$, lifting arbitrary elements in $\EE(V,g)_1$ and $\EE(V,g)_2$ under the homomorphism
$\pi$. First, let us check that $at=ta$ in $\RR_2$, for $a\in \End_g(V)L\sub\RR_1$.
Note that here $at=T(a)\circ t$, $ta=T(t)\circ a$. From our description of $T(t)$, we immediately get that $ta$ is represented by
the map
\begin{equation}\label{ta-diagram}
\begin{diagram}
&&E\\
&&\dTo{a}\\
V^\vee L^2\ot F&\rTo{}&V^\vee L\ot F&\rTo{}&E
\end{diagram}
\end{equation}
On the other hand, from the description \eqref{T-a-eq} of $T(a)$ it follows immediately
that $at=T(a)\circ t$ is represented by the same chain map \eqref{ta-diagram} as $ta$.

Similarly, let us consider an element $A\in\RR_2$ represented by a closed map
\begin{diagram}
E\\
\dTo{A_0}&\rdTo{A_1}\\
V^\vee L^2\ot F&\rTo{\de}&V^\vee L\ot F&\rTo{\de}&E
\end{diagram}
where $m_3(\de,\de,A_0)+m_2(\de,A_1)=0$. Then one can easily check that $T(A)\circ t$ and $T^2(t)\circ A$
are both represented by the map
\begin{diagram}
&&E\\
&&\dTo{A_0}&\rdTo{A_1}\\
V^\vee L^3\ot F&\rTo{}&V^\vee L^2\ot F&\rTo{}&V^\vee L\ot F&\rTo{}&E
\end{diagram}

Since $t$ is a central nonzero divisor, it follows that the algebra $\RR$ is the Rees algebra of a filtered algebra $(A,F_\bullet A)$.
Furthermore, by Steps 4 and 6 we have
$$\gr^F(A)\simeq\RR/(t)\simeq \EE(V,g)^{op},$$ 
which proves (i).

\medskip

\noindent
{\bf Step 9}.
For $a,b\in \End_g(V)L\sub \RR_1$, the product $ab=T(a)\circ b$ in $\RR$ 
can be easily computed using the representation \eqref{T-a-eq} for $T(a)$: we get 
$$T(a)\circ b=(A_0,A_1,A_2)\in \hom^0(E_0,E_2)=\End(V)L^2\oplus \End(V)L \oplus \Hom^0(E,E),$$
with
$$A_0=ba, \ \ A_1=m_3(a,\de_1,b)+(\mu_a\ot\id_F)\circ b, \ \ A_2=m_4(\de_1,a,\de_1,b),$$
where $\mu_a$ is given by \eqref{mu-a-eq}.

From this point on we use the trivialization $L=S\cdot\xi_F$. Assume for simplicity of notation 
that $a=v_a\ot v^*_a$, $b=v_b\ot v^*_b$, for $v_a,v_b\in V$, $v^*_a,v^*_b\in V^\vee$ (the general case is proved similarly).

The evaluation map
$\de_1:\Hom^1(F,E)\ot F\to E[1]$ corresponds to the identity element $\sum e_i\ot e^*_i$ in 
$\Hom^1(F,E)^\vee\ot\Hom^1(F,E)\simeq V\ot V^\vee$.
Thus, the map $\mu_a:V^\vee\to V^\vee$ sends $v^*$ to 
$$\sum_i \lan v^*_a,e_i\ran\cdot m_3(e_i^*,v_a,v^*)=m_3(v^*_a,v_a,v^*).$$
Similarly,
$$m_3(a,\de_1,b)=m_3(v_a,v^*_b,v_b)\ot v^*_a,$$
$$m_4(\de_1,a,\de_1,b)=m_4(v^*_a,v_a,v^*_b,v_b).$$
Thus,
$$A_1=m_3(v_a,v^*_b,v_b)\ot v^*_a+v_b\ot m_3(v^*_a,v_a,v^*_b), \ \ A_2=m_4(v^*_a,v_a,v^*_b,v_b).$$
Using the operators $r$ and $r'$ we can rewrite the formula for $A_1$ as
$$A_1=r(b)a+br'(a).$$

Now recall that $\RR_2=\Hom^0(E_0,E_2)$ is the subspace of $\hom^0(E_0,E_2)$ consisting of $(A_0,A_1,A_2)$
such that 
$$s(A_0)+\tr(gA_1)=0,$$
where $s(v\ot v^*)=m_3(v^*,\xi_F,v)$. 
Thus, we can define the splitting $\si$ of the projection $\pi_2:\RR_2=\Hom^0(E_0,E_2)\to \End(V)$ by
setting
$$\si(A)=(A,-s(A)h,0)\in \Hom^0(E_0,E_2)\sub \hom^0(E_0,E_2),$$
for $A\in\End(V)$.
Rewriting the element $(A_0,A_1,A_2)\in\Hom^0(E_0,E_2)$ as
$\si(A_0)+(A_1+s(A)h)t+A_2t^2$ we get \eqref{R-product-formula-eq}.
\ed

\begin{cor} Under the assumptions of Theorem \ref{spherical-filtered-thm}, if in addition $g$ is invertible,
then the isomorphism class of the corresponding
filtered algebra $(A,F_\bullet A)$ 
is determined by the higher products \eqref{higher-products-main-eq}. 
\end{cor}

\Pf . Indeed, in this case by Lemma \ref{E-deg-1-lem}(iii), 
$\EE(V,g)$ is generated in degree $1$ and has quadratic defining relations.
Hence, the same is true for $\RR$. But, by Theorem \ref{spherical-filtered-thm}(ii), the product $\RR_1\ot \RR_1\to\RR_2$
is determined by the higher products \eqref{higher-products-main-eq}.
\ed

Until the end of this section we fix an $n$-pair $(E,F)$ and keep the notation of Theorem \ref{spherical-filtered-thm}.
In particular, $\RR$ is the graded $S$-algebra defined by \eqref{R-algebra-eq}. Let us define the structure of
a graded $\RR-\RR$-bimodule on
$$\bigoplus_{i\in\Z}\Hom^1(E_i,E_0)=\bigoplus_{i\ge 0}\Hom^1(E_i,E_0)$$ 
as follows. The grading of $\Hom^1(E_i,E_0)$ is set to be $-i$.
The left and right multiplication of $x\in \Hom^1(E_i,E_0)$ by $a\in \Hom^0(E_0,E_j)$, with $j\le i$, are given by
$$a\cdot x=T^{-j}(a\circ x), \ \ x\cdot a=x\circ T^{i-j}(a).
$$
We would like to calculate this $\RR-\RR$-bimodule, which together with the algebra $\RR$ completely determines the 
full subcategory with the objects $(E_i)$. This is important for proving Theorem A since, as we will show in \ref{min-non-form-alg-sec} and \ref{triple-prod-sec},
this subcategory can be equipped with a unique minimal $A_\infty$-structure up to rescaling (and gauge equivalence).

\begin{lem}\label{bimodule-lem}
Let $B$ be a non-negatively graded $S$-algebra, such that $B_0=S$ and all the graded components $B_i$ are finitely generated projective $S$-modules. Let $M$ be a graded $B-B$-bimodule such that
$M$ is isomorphic to $B^*$ (the restricted dual of $B$)
as a graded left $B$-module and as a graded right $B$-module. Then there exists an 
automorphism $\phi:B\to B$, preserving the grading, such that
$M$ is isomorphic to $(\sideset{_{\id}}{_{\phi}}{B})^*$ as a graded $B-B$-bimodule.
\end{lem}

\Pf . Note that $M_{-i}\simeq B_i^\vee$, so all the graded components of $M$ are finitely generated projective modules.
Consider the restricted dual $M^*$. Then $M^*$ is isomorphic to $B$ as a left and as a right graded $B$-module.
This easily implies the claim.
\ed

\begin{prop}\label{bimodule-prop} 
Under the assumptions of Theorem \ref{spherical-filtered-thm}, assume in addition that $g$ is invertible. 
Then there is an isomorphism of graded $\RR-\RR$-bimodules
\begin{equation}\label{restr-dual-R-isom}
\bigoplus_{i\ge 0}\Hom^1(E_i,E_0)\simeq (\sideset{_{\id}}{_{\phi}}\RR)^*\ot_S \LL_E,
\end{equation}
where $\phi$ is a graded automorphism of $\RR$, such that $\phi(t)=t$ and the automorphism
$\ov{\phi}$ of $\RR/(t)\simeq \EE(V,g)^{op}$, induced by $\phi$, is equal 
the restriction of the automorphism $\Ad(g^{-1}):x\mapsto g^{-1}xg$ of $\End(V)[z]^{op}$.
\end{prop}

\Pf . Without loss of generality we can assume that $\CC$ is generated by $(E,F)$.
Then by Lemma \ref{invertible-g-lem}, the object $E_0=E$ is $1$-spherical.
It follows that all the objects $E_i=T^i(E_0)$ are $1$-spherical in $\CC$. 
As before, we fix a trivialization $\Hom^1(F,F)=S\cdot \xi_F$, 
and consider the identification $\Hom^1(F,E)\simeq V^\vee$,
such that the second of the pairings
\eqref{two-pair-eq} gets identified with the the natural pairing $V\ot V^\vee\to S$, while the first of these pairings
is given by \eqref{g-def-eq}. 
We also use the isomorphisms $\Hom^1(E_i,E_i)\simeq\Hom^1(E,E)=\LL_E$.

\medskip

\noindent
{\bf Step 1}. There exists an isomorphism of bimodules \eqref{restr-dual-R-isom} for some $\phi$.

We can use the perfect pairing (given by the composition)
$$\Hom^1(X,E_0)\ot \Hom^0(E_0,X)\to \Hom^1(E_0,E_0)=\LL_E$$
to define an isomorphism
$$\Hom^1(X,E_0)\ot \LL_E^{-1}\rTo{\sim} \Hom^0(E_0,X)^\vee.$$
These isomorphisms are functorial in $X$, which immediately implies that summing these isomorphisms over $X=E_i$,
we get an isomorphism of right $\RR$-modules
\begin{equation}\label{restr-dual-R-isom-bis}
\bigoplus_{i\ge 0}\Hom^1(E_i,E_0)\ot \LL_E^{-1}\simeq \bigoplus_{i\ge 0}\Hom^0(E_0,E_i)^\vee=\RR^*.
\end{equation}

On the other hand, using the perfect pairings
$$\Hom^0(X,E_i)\ot \Hom^1(E_i,X)\to \Hom^1(E_i,E_i)\simeq \LL_E$$
we similarly get an isomorphism of right $\RR$-modules
$$\RR=\bigoplus_{i\ge 0}\Hom^0(E_0,E_i)\simeq\bigoplus_{i\ge 0}\Hom^1(E_i,E_0)^\vee\ot\LL_E.$$
Dualizing we get another isomorphism of the form \eqref{restr-dual-R-isom-bis}
which is compatible with the left $\RR$-module structures.
By Lemma \ref{bimodule-lem}, there exists an automorphism $\phi:\RR\to\RR$ such that \eqref{restr-dual-R-isom} holds.

It remains to check that $\phi(t)=t$ and to calculate the action of $\phi$ on $\RR/(t)$.
To this end we will calculate some compositions of morphisms between $E_0=E$ and $E_1=[V^\vee\ot F\to E]$.
Recall that we have a canonical decomposition 
$$\Hom(E_0,E_1)=\RR_1=\End_g(V)\oplus S\cdot t.$$

\medskip

\noindent
{\bf Step 2}. We construct canonical identifications 
\begin{equation}\label{Hom-E1-E0-eq}
\Hom^1(E_1,E_0)=\End(V)/(S\cdot \id) \oplus \LL_E,
\end{equation}
$$\tau:\Hom^1(E_1,E_1)\rTo{\sim} \LL_E,$$
such that the composition 
$$\Hom^1(E,E)\oplus\Hom^1(V^\vee\ot F,E)\oplus \Hom^1(V^\vee\ot F,V^\vee\ot F)=\hom^1(E_1,E_1)
\to\Hom^1(E_1,E_1)\rTo{\tau} S,$$
where the first arrow is the natural projection to cohomology, is given by
\begin{equation}\label{tau-formula-eq}
(\xi,x,y\ot y^*)\mapsto \xi+\lan y^*,gy\ran,
\end{equation}
where $\xi\in \LL_E$, $y\in V$, $y^*\in V^\vee$.

By definition, the complex $\hom(E_1,E_0)$ has the form
$$\Hom^0(E,E)\to \Hom^1(V^\vee\ot F,E)\oplus\Hom^1(E,E),$$
where the differential maps $\id_E$ to the evaluation morphism in $\Hom^1(V^\vee\ot F,E)$.
This immediately leads to the identification \eqref{Hom-E1-E0-eq}.

We have
$$\hom^0(E_1,E_1)=\Hom^0(E,E)\oplus \Hom^0(E,V^\vee\ot F)\oplus\Hom^0(V^\vee\ot F,V^\vee\ot F)$$
and the part of the differential $\hom^0(E_1,E_1)\to\hom^1(E_1,E_1)$ that maps to
the summands $\Hom^1(E,E)\oplus \Hom^1(V^\vee\ot F,V^\vee\ot F)$ is the map
$$\Hom^0(E,V^\vee\ot F)\rTo{(-(?\circ\de),\de\circ ?)} \Hom^1(E,E)\oplus \Hom^1(V^\vee\ot F,V^\vee\ot F),$$
with both components induced by the evaluation map $V^\vee\ot F\to E[1]$.
We can identify this map with the map
$$V\ot V^\vee \to \LL_E\oplus V\ot V^\vee: y\ot y^*\mapsto (-\lan y^*,gy\ran, y\ot y^*).$$
Thus, the map $\hom^1(E_1,E_1)\to S$ given by \eqref{tau-formula-eq}
descends to a map on cohomology, $\tau:\Hom^1(E_1,E_1)\to S$. It is clear from the definition that
$\tau$ is surjective. Since $E_1$ is $1$-spherical, we deduce that $\tau$ is an isomorphism.

\medskip

\noindent
{\bf Step 3}. Let $T:\LL_E=\Hom^1(E,E)\to \Hom^1(E_1,E_1)$ be the map induced by the spherical twist $T$. Then
$\tau\circ T$ is the identity map of $\LL_E$. Hence, the map
$$\Hom^0(E_0,E_1)\ot\Hom^1(E_1,E_0)\to \Hom^0(E_0,E_0)=S: a\ot x\mapsto a\cdot x=T^{-1}(a\circ x),$$
sends $a\ot x$ to $\tau(a\circ x)$.

Indeed, it is easy to see that $T(\xi)$ is given by the element $\xi\in \Hom^1(E,E)\sub\hom^1(E_1,E_1)$.

\medskip

\noindent
{\bf Step 4}. For $B\in \End(V)/(S\cdot \id)\sub \Hom^1(E_1,E_0)$ one has
$B\circ t=0$ in $\Hom^1(E_0,E_0)$ and $t\circ B=0$ in $\Hom^1(E_1,E_1)$.
Also, viewing $\xi\in\LL_E$ as an element of $\Hom^1(E_1,E_0)$ (see \eqref{Hom-E1-E0-eq}), we get
$\xi\circ t=\xi$ and $\tau(t\circ \xi)=\xi$. Hence, for any $x\in \Hom^1(E_1,E_0)$, we have
$$t\cdot x=x\cdot t$$
in the bimodule $\bigoplus_{i\ge 0}\Hom^1(E_i,E_0)$.

Indeed, recall that $t$ corresponds to the element $\id_E\in\Hom^0(E,E)\sub \hom^0(E,E_1)$.
The vanishing of the composition $B\circ t=0$ is clear from the composition rule:
\begin{diagram}
&& E\\
&&\dTo{\id_E}\\
V^\vee\ot F&\rTo{\de}&E\\
\dTo{B}\\
E
\end{diagram} 

The composition $t\circ B$ is calculated by the diagram
\begin{diagram}
&&V^\vee\ot F&\rTo{\de}&E\\
&&\dTo{B}\\
&& E\\
&&\dTo{\id_E}\\
V^\vee\ot F&\rTo{\de}&E
\end{diagram} 
Thus, it belongs to the summand $\Hom^1(V^\vee\ot F,E)\sub\hom^1(E_1,E_1)$,
which is annihilated by $\tau$, hence zero in cohomology.

The composition $\xi\circ t$ is calculated by the diagram
\begin{diagram}
&& E\\
&&\dTo{\id_E}\\
V^\vee\ot F&\rTo{\de}&E\\
&&\dTo{\xi}\\
&& E
\end{diagram} 

Finally, the composition $t\circ \xi$ is calculated by the diagram
\begin{diagram}
V^\vee\ot F&\rTo{\de}&E\\
&&\dTo{\xi}\\
&& E\\
&&\dTo{\id_E}\\
V^\vee\ot F&\rTo{\de}&E
\end{diagram} 
so it is given by the element $\xi\in\Hom^1(E,E)\sub\hom^1(E_1,E_1)$, which implies the assertion.

\medskip

\noindent
{\bf Step 5}. For $A\in\End_g(V)\sub\Hom^1(E_0,E_1)$ and $B\in \End(V)/(S\cdot \id)$, we have
$$B\circ A=\tr(BgA),$$
$$\tau(A\circ B)=\tr(AgB).$$

Indeed, the composition $B\circ A$ is calculated by the diagram 
\begin{diagram}
E\\
\dTo{A}\\
V^\vee\ot F&\rTo{\de}&E\\
\dTo{B}\\
E
\end{diagram} 
We claim that the composition of the vertical arrows is $\tr(BgA)$. Indeed, it is easy to see
using our conventions, that for $A=v^*\ot v$ and $B=w^*\ot w$, this composition will be
$\lan v^*,w\ran\cdot \lan w^*,gv\ran$, which gives our claim.

Finally, the composition $A\circ B$ is calculated by the diagram
\begin{diagram}
V^\vee\ot F&\rTo{\de}&E\\
\dTo{B}\\
E\\
\dTo{A}\\
V^\vee\ot F&\rTo{\de}&E
\end{diagram} 
Note that this composition will have a component in $\Hom^1(V^\vee\ot F,V^\vee\ot F)$ given by
the composition $m_2(A,B)$ of the vertical arrows, as well, as a component in
$\Hom^1(V^\vee\ot F,E)$ given by $m_3(\de,A,B)$. However, the latter component does not give any
contribution to the cohomology class, due to formula \eqref{tau-formula-eq}. For $A=v^*\ot v$ and
$B=w^*\ot w$, we have $m_2(A,B)=\lan w^*,v\ran\cdot w\ot v^*$, so applying $\tau$ we get
$$\tau(m_2(A,B))=\lan w^*,v\ran\cdot \lan v^*,gw\ran=\tr(AgB).$$

\medskip

\noindent
{\bf Step 6}. Since $\RR$ is generated in degree $1$, the automorphism $\phi$ is uniquely determined by its restriction
to $\RR_1$, which is in turn uniquely determined by the equation
$$x\cdot \phi(a)=a\cdot x$$
in $\LL_E=\Hom^1(E_0,E_0)$, where $a\in \RR_1$, $x\in \Hom^1(E_1,E_0)$.
Hence, by Step 4, we deduce that $\phi(t)=t$. Furthermore, still by Step 4, for $a\in\RR_1$,
one has $B\cdot a=0$ for all $B\in\End(V)/(S\cdot\id)\sub\Hom^1(E_1,E_0)$ if and only if $a\in S\cdot t$.
Therefore, for $A\in\End_g(V)$, the element
$\phi(A)\mod S\cdot t$ is determined by the products $B\cdot \phi(A)$ in $\Hom^1(E_0,E_0)$.
Now the calculation of Step 5 implies that $\phi(A)\equiv g^{-1}Ag \mod S\cdot t$.
\ed

Recall that $\RR$ is identified with the Rees algebra $\RR(A)$ (see Theorem \ref{spherical-filtered-thm}).
Thus, a graded automorphism $\phi$ of $\RR$ such that $\phi(t)=t$ is the same thing as a filtered automorphism of
$A$. In view of Proposition \ref{bimodule-prop} it is important to study such automorphisms. The next result shows
that in some cases $\phi$ is uniquely determined by the induced automorphism of $\gr^F A$.

\begin{prop}\label{filt-aut-prop} 
Let $V=S^n$, where $n\ge 2$, and let $g\in\End_S(V)\ot \LL$ be an invertible element (where $\LL$ is a locally free $S$-module of rank $1$). Assume that either $n\ge 3$ or
$\tr(g)$ is a generator of $\LL$.
Let $(A, F_\bullet A)$ be a filtered $S$-algebra, such that
\begin{equation}\label{gr-F-eq-bis}
\gr^F A\simeq \EE(V,g)^{op}.
\end{equation}
Suppose $\phi_1$ and $\phi_2$ are automorphisms of $A$, preserving the filtration (i.e., $\phi_i F_jA=F_jA$), such that 
the induced automorphisms of $\gr^F A$ are the same. Then $\phi_1=\phi_2$.
\end{prop}

\Pf . Let us consider the automorphism $\phi=\phi_1^{-1}\phi_2$ of $A$. Then $\phi$ still preserves the filtration
and induces the identity on $\gr^F A$. It is easy to check that setting for $a\in F_n A$
$$D(a)=\phi(a)-a\mod F_{n-2}A$$
we get a well defined derivation $D:\gr^F A\to \gr^F A$ of degree $-1$. Using \eqref{gr-F-eq-bis} and Proposition
\ref{derivations-E-prop}, we deduce that $D=0$. In particular, $\phi(a)=a$ for $a\in F_1 A$.
But Lemma \ref{E-deg-1-lem}(i) implies that $A$ is generated by $F_1 A$ (since this is true for $\gr^F A\simeq \EE(V,g)^{op}$).
Hence, it follows that $\phi(a)=a$ for all $a\in A$.
\ed

\section{Connection with noncommutative orders over stacky curves}\label{orders-sec}

\subsection{Filtered algebras and orders}\label{filtered-alg-sph-orders-sec}

In this section we work over a fixed ground field $k$.
Let $A$ be a filtered algebra over $k$ equipped with an isomorphism \eqref{gr-E-eq} for some $g\in\P\End(V)$, and let
 $Z\sub A$ be its center. We equip $Z$ with the induced filtration. 
 
\begin{lem}\label{center-lem} 
The algebra $A$ is finitely generated, prime and of GK-dimension $1$.
Hence, $A$ is Noetherian and finite over its center $Z$, which is a $1$-dimensional domain, finitely generated as $k$-algebra. Also, $A$ is an order in a central simple algebra over the quotient field of $Z$.
\end{lem}

\Pf . Since the algebra $\EE(V,g)$ is generated by degree $1$ and degree $2$ elements (see Lemma \ref{E-deg-1-lem}(ii)), we deduce that $A$ is generated by $F_2A$, hence, it is finitely generated.
Given a nonzero ideal $I\sub A$, let $I_0\sub \End(V)$ be the set of all elements $x$ such that $xz^n$ appears
as an initial form of an element of $I$ for some $n$. Then $I_0$ is a nonzero ideal, hence, $I_0=\End(V)$.
Hence, for a pair of nonzero ideals $I,J\sub A$ we have $I_0J_0\neq 0$, so $IJ\neq 0$, which shows that $A$ is prime.
We have $\dim F_iA/F_{i-1}A=n^2$ for $i>1$, so the GK-dimension of $A$ is one.
Now the results of \cite{Small-Warfield} and \cite{Schelter} imply that $A$ and $Z$ are Noetherian, $A$ is finite over $Z$, and
$Z$ has dimension $1$.   

Note that the center of $\gr^F(A)\simeq \EE(V,g)^{op}$ is either $k[z^2,z^3]\sub k[z]$, in the case when $\tr(g)\neq 0$,
or $k[z]$, when $\tr(g)=0$. Thus, $\gr^F(Z)$ is a graded $k$-subalgebra in $k[z]$,
i.e., a group algebra of a subsemigroup in natural numbers. This easily implies that 
the algebra $\RR(Z)$ is a domain, finitely generated as a $k$-algebra.
Next, the fact that $\gr^F(A)\simeq\EE(V,g)^{op}$ is torsion free as a module over $\gr^F(Z)\sub k[z]$ implies that
$A$ is torsion free as a $Z$-module.
Let $K$ be the quotient field of $Z$. Then $A\ot_Z K$ is a finite-dimensional prime algebra over $K$ with the center $K$,
so it is a central simple algebra over $K$.
\ed 



Next, we would like to extend $A$ to a sheaf of algebras over a projective curve compactifying $\Spec(Z)$.
The first obvious choice is to consider the Rees algebras 
$$\RR(A)=\bigoplus_{m\ge 0} F_m A \ \text{  and  } \RR(Z)=\bigoplus_{m\ge 0} F_m Z$$ 
and to consider the corresponding $\Proj$-construction.
However, the resulting structures are not always easy to analyze. Namely, the problem arises when
$\gr^F(Z)$ is contained in $k[t^d]$ for some $d\ge 2$. It turns out that a better behaved construction is
provided by the stacky version of $\Proj$, which we denote by $\Proj^{st}$.

Namely, for any commutative non-negatively graded $k$-algebra $B=\bigoplus_{n\ge 0} B_n$, where $B_0=k$, one can define
a stack
\begin{equation}\label{Proj-stack-eq}
\Proj^{st}(B):=\Spec(B)\setminus \{0\}/\G_m,
\end{equation}
where $0$ is the point corresponding to the augmentation ideal $B_+$.
Assuming in addition that $B$ is finitely generated, we have an equivalence of the category $\Coh(\Proj^{st}(B))$ with
the quotient of the category of finitely generated graded $B$-modules by the subcategory of finite-dimensional modules.
Note that we have a natural line bundle $\OO(1)$ on $\Proj^{st}(B)$ such that elements of $B_n$ can be viewed as
global sections of $\OO(n)$. The coarse moduli for $\Proj^{st}(B)$ is the usual scheme $\Proj(B)$.

Now starting with an algebra $A$ as above we define the stacky curve $C$ by 
$$C:=\Proj^{st} \RR(Z).$$
Let $d\ge 1$ be the maximal such that $\gr^F(Z)\sub k[t^d]$. We will see below that $d$ measures
the ``stackiness" of $C$ (see Lemma \ref{stacky-sheaf-lem}(i)). In particular, $d=1$ if and only if $C$ is the usual curve.

Let us denote by $t$ the element $1\in \RR_1(Z)=F_1A\cap Z$.
Note that $t$ is a non-zero-divisor, and $\RR(A)/t\RR(A)\simeq \gr^F(A)$, $\RR(Z)/t\RR(Z)\simeq \gr^F(Z)$.
Since $\gr^F(A)\simeq\EE(V,g)^{op}$ if finitely generated as a $\gr^F(Z)$-module (see the proof of Lemma \ref{center-lem}), 
we deduce that $\RR(A)$ is finitely generated as an $\RR(Z)$-module.

Thus, localizing $\RR(A)$ we get a sheaf of coherent $\OO$-algebras $\AA$ on $C$.
More precisely, $\wt{\RR(A)}$ is a $\G_m$-equivariant sheaf of coherent $\OO$-algebras on $\Spec(\RR(Z))$. The sheaf of algebras
$\AA$ is obtained by restricting $\wt{\RR(A)}$ to
$\Spec(\RR(Z))\setminus\{0\}$ and descending the result to $C$.

Note that the complement to the divisor $(t=0)$ in $C$ is identified with the affine curve $\Spec(\RR(Z)/(t-1))=\Spec(Z)$,
the restriction of $\AA$ to it gets identified with the coherent sheaf of algebras corresponding to $A$.

\begin{lem}\label{stacky-sheaf-lem} 
(i) The pointed stacky curve $C$ is neat, and the divisor $(t=0)$ is the unique
stacky point $p=\wt{p}/\mu_d$ (where $\wt{p}=\Spec(k)$).
Thus, we have
$\OO_C(1)\simeq \OO_C(p)$.
There is a natural identification of the cotangent space $T^*_{\wt{p}}C$ with $\chi$, the $1$-dimensional space
on which $\mu_d$ acts with weight $1$.

\noindent
(ii) The sheaf $\AA$ is an order on $C$, i.e., a torsion-free coherent sheaf of $\OO$-algebras, 
whose stalk at the generic point is a central simple $K(X)$-algebra. Furthermore, the center of $\AA$ is $\OO_C$.

\noindent
(iii) One has a natural isomorphism of algebras on $p$,
\begin{equation}\label{AA-res-to-p-eq}
\AA|_p\simeq \rho_*\End(V)^{op},
\end{equation}
where $\rho:\wt{p}\to p$ is the natural morphism.
Similarly, we have a natural isomorphism $\AA(mp)|_p\simeq \rho_*\End(V)$ for every $m\in\Z$, compatible with the above
isomorphism via the identification $\OO(mp)|_p\simeq\chi^{-m}$.
The rank of $\AA$ is equal to $dn^2$.

\noindent
(iv) The natural map $F_m A=\RR_m(A)\to H^0(C,\AA(mp))$ is an isomorphism for $m\in\Z$ (where $F_{-1}A=0$),
so we have an isomorphism of graded algebras
\begin{equation}\label{R-A-C-eq}
\RR(A)\simeq \bigoplus_{n\ge 0}H^0(C,\AA(n)).
\end{equation}
 One has 
$h^1(\AA)=1$ and $h^1(\AA(p))=0$. 
In particular, $\AA$ is a weakly spherical order in the sense of Definition \ref{sph-order-def}.
\end{lem}

\Pf . (i) We have seen that $\RR(Z)$ is a domain, so $C$ is integral. Also,
the coarse moduli is $\Proj \RR(Z)$ which is a projective curve.
The divisor 
$$p:=(t=0)\sub C$$ 
can be identified with $\Proj^{st}(\gr^F(Z))$.
Since $\gr^F(Z)\sub k[z^d]$, and these algebras agree in all sufficiently high degrees,
we see that $\Spec(\gr^F(Z))$ is an affine line with the pinched origin. In particular,
$$p=(\Spec(\gr^F(Z))\setminus 0)/\G_m\simeq B\mu_d.$$

Set $S=\Spec(\RR(Z))\setminus\{0\}$.
We can view $t$ as a map $S\to \A^1$ and the fiber over $0$, $D\sub S$ 
is a closed $\G_m$-orbit with the stabilizer $\mu_d$. Since $D=\Spec(\gr^F(Z))\setminus\{0\}$ is smooth, 
the surface $S$ is smooth near $D$.
By the argument of Luna's \'etale slice theorem (see \cite{Luna}), there exists a smooth 
$\mu_d$-invariant locally closed curve $\Si\sub S$ through the point $z=1$ of $D$ such that the induced map of stacks
$\Si/\mu_d\to S/\G_m$ is \'etale.

The identification of the cotangent space at $\wt{p}$ with $\chi$ comes from the fact that $p$ is given by the
equation $t=0$, where $t$ is a section of $\OO_C(1)$.


\noindent
(ii) This first assertion follows from Lemma \ref{center-lem}. Also, we know that $\RR(Z)$ is the center of $\RR(A)$,
so $\OO_C$ is the center of $\AA$.

\noindent
(iii) Note that $\AA/\AA(-p)$ is the localization of the graded module $\gr^F(A)\simeq \EE(V,g)^{op}$ over $\gr^F(Z)$.
Up to finite-dimensional pieces, this is the same as considering $\End(V)^{op}[z]$ as a $k[z^d]$-module, which easily implies
the isomorphism \eqref{AA-res-to-p-eq}. 
If we identify sheaves on $p$ with $\mu_d$-representations then the functor
$\rho_*$ is given by tensoring with the regular representation of $\mu_d$, so the image of $\rho_*$ is stable under
tensoring with $\chi$.

Finally, since $p$ is a smooth (stacky) point of $C$, the sheaf $\AA$ is locally free near $p$. Thus, by considering ranks in
the isomorphism \eqref{AA-res-to-p-eq}, we obtain that the rank of $\AA$ is $dn^2$.

\noindent
(iv) Note that the natural isomorphism 
$$\a:A\rTo{\sim} H^0(C\setminus\{p\},\AA)$$
sends $F_m A$ to $H^0(C,\AA(m))$.
It is easy to check that the induced map
of the associated graded spaces has as components the natural maps
$$F_m A/F_{m-1}A\simeq \EE(V,g)_m\hra \End(V)=H^0(p,\rho_*\End(V))\simeq H^0(p,\AA(mp)|_p),$$
where we use (iii).
It follows that $\a^{-1}H^0(C,\AA(mp))=F_m A$ for every $m\in\Z$. In particular, $H^0(C,\AA)=k$. 

Note that for $m\ge 1$ we have $\dim F_m A=mn^2$. Hence, for sufficiently large $m$ we have $\chi(\AA(mp))=h^0(\AA(mp))=mn^2$.
Hence, $\chi(\AA)=0$ and $\chi(\AA(p))=n^2$. Thus, since $h^0(\AA)=1$ and $h^0(\AA(p))=n^2$,
we get $h^1(\AA)=1$ and $h^1(\AA(p))=0$.
\ed

We will now start using the framework of noncommutative projective geometry as developed in \cite{AZ}.
Thus, for a (non-negatively) graded Noetherian algebra $B$ we consider the category $\qgr B$,
the quotient of the category of finitely generated right $B$-modules by the subcategory of torsion modules.
This category is equipped with a grading shift functor $M\mapsto M(1)$ and a special object $\OO$ (the image of $B$).
We refer to these data as $\Proj^{nc} B$, the {\it noncommutative $\Proj$ of $B$}.

\begin{prop}\label{Proj-order-prop} 
The category $\Coh(\AA^{op})$ of coherent right $\AA$-modules is equivalent to the category $\qgr \RR(A)$.
\end{prop}

\Pf . This is proved similarly to \cite[Prop.\ 2.3]{AKO}. Let $\wt{\RR(A)}$ be the coherent sheaf of $\OO$-algebras on $\Spec(\RR(Z))$
associated with $\RR(A)$, viewed as an $\RR(Z)$-algebra. The category $\Coh(\AA^{op})$ can be identified with the category of
$\G_m$-equivariant coherent sheaves of $\wt{\RR(A)}^{op}$-modules on $\Spec(\RR(Z))\setminus\{0\}$.
The latter category is the quotient of the category of $\G_m$-equivariant coherent sheaves of $\wt{\RR(A)}^{op}$-modules on $\Spec(\RR(Z))$
by the subcategory of sheaves with support in $\{0\}$. This immediately leads to the required equivalence. 
\ed

Below we refer to the condition $\chi$ introduced in \cite[Def.\ 3.7]{AZ} which is useful in the context
of noncommutative projective geometry. 
We also use the notion of the {\it cohomological dimension of $\Proj^{nc}$} of
a graded Noetherian algebra $B$ defined in terms of the cohomology functor 
$$H^i(\cdot):=\Ext^i_{\qgr B}(\OO,\cdot)$$ (see \cite[Sec.\ 7]{AZ}).
Note that for $i\ge 1$, one has an isomorphism
\begin{equation}\label{Hi-lim-formula}
H^i(M)\simeq \lim_{m\to\infty} \Ext^{i+1}_{B^{op}}(B/B_{\ge m},M)
\end{equation}
(see \cite[Prop.\ 7.2]{AZ}). Thus, finiteness of the cohomological dimension of $\Proj^{nc} B$ is equivalent to 
finiteness of the cohomological dimension of the functor 
$$\Gamma_{B_+}=\lim_{m\to\infty}\Hom_{B^{op}}(B/B_{\ge m},M).$$

\begin{cor}\label{Noeth-chi-cd-cor}
The algebra $\RR(A)$ is right Noetherian, satisfies the condition $\chi$, and $\Proj^{nc} \RR(A)$
has cohomological dimension $\le 1$. The same is true for the algebra $\RR(A)^{op}$.
\end{cor}

\Pf . We can rewrite isomorphism \eqref{R-A-C-eq} as
$$\RR(A)\simeq \bigoplus_{n\ge 0}H^0(C,\AA(n)) \simeq\bigoplus_{n\ge 0}\Hom_{\AA^{op}}(\AA,\AA(n)).$$
Hence, from Proposition \ref{Proj-order-prop}, by \cite[Thm.\ 4.5]{AZ}, 
we get that $\RR(A)$ is right Noetherian and satisfies $\chi_1$.

Next, we claim that for every coherent right $\AA$-module $M$ the spaces $\Ext^j_{\AA^{op}}(\AA,M)$ are
finite-dimensional, $\Ext^{>1}_{\AA^{op}}(\AA,M)=0$, and $\Ext^j_{\AA^{op}}(\AA,M(i))=0$ for $i\gg 0$. 
Indeed, this immediately follows from the
identification $\Ext^j_{\AA^{op}}(\AA,M)\simeq H^j(C,M)$ (note that the latter cohomology is isomorphic
to the cohomology of the push-forward of $M$ to the coarse moduli space of $C$).

By Proposition \ref{Proj-order-prop}, we deduce a similar statement for the cohomology functor $H^*$ on the
category $\qgr \RR(A)$. In particular, we see that the cohomological dimension of $\Proj^{nc} \RR(A)$ is $\le 1$.
Now the fact that $\RR(A)$ satisfies $\chi$ follows from \cite[Thm.\ 7.4(2)]{AZ}.

The last assertion follows from the fact that $\gr^F(A^{op})\simeq  \EE(V,g)\simeq \EE(V^\vee,g^*)^{op}$,
so we can repeat the argument with $A$ replaced by $A^{op}$. 
\ed


\subsection{Spherical orders and duality}\label{spherical-dual}

Let $\AA$ be an order over a proper stacky curve $C$ with a stacky point $p\simeq B\mu_d$ 
such that $\AA|_p\simeq\rho_*\End(V)^{op}$, where
$\rho:\wt{p}\to p$ is the $\mu_d$-covering of $p$ by $\wt{p}\simeq\Spec(k)$.
Then we can view $\rho_*V$ as a right $\AA$-module supported at $p$. Note that if $d=1$ then this module is $V\ot\OO_p$.

\begin{lem}\label{Perf-gen-lem} 
Let $\AA$ be an order over a neat pointed stacky curve $C$ with the unique stacky point $p\simeq B\mu_d\in C$, such that
$\AA|_p\simeq \rho_*\End(V)^{op}$, where $V$ is a finite-dimensional vector space. 
Then the pair of $\AA^{op}$-modules $(\AA,\rho_*V)$ (resp., $(\AA, \AA(-p))$)
split generates $\Perf(\AA^{op})$.
\end{lem} 

\Pf . First, we note that the $\AA$-module $\AA|_p$ is the direct sum of several copies of $\rho_*V$.
In particular, $\rho_*V$ is in $\Perf(\AA^{op})$, and the pairs $(\AA,\rho_*V)$ and $(\AA,\AA(-p))$ split generate
the same subcategory.
Next, using the exact sequences
$$0\to \AA(-(m+1)p)\to \AA(-mp)\to \AA|_p\to 0$$ 
for $m\ge 0$, we see that the subcategory $\lan \AA,\rho_*V\ran$ generated by our objects contains all $\AA(-mp)$ for 
$m\ge 0$. 

We claim that for any coherent right $\AA$-module $M$ there exists a surjection of the form 
$\bigoplus_{i=1}^N \AA(-n_i)\to M$ for some $n_i\ge 0$.
Indeed, let $p=B\mu_d$. 
This is a unique stacky point, and $\mu_d$ acts on the fiber of $\OO_C(-p)$ at $p$
by the identity character.
Hence, using \cite[Prop.\ 5.2]{OS}, we see that the bundle $\EE=\bigoplus_{i=1}^d \OO(-ip)$ over $C$ has
the property that the map
$$\pi^*\pi_*\und{\Hom}(\EE,\FF)\ot \EE\to \FF,$$
where $\pi:C\to\ov{C}$ is the coarse moduli map, is surjective for every quasicoherent sheaf $\FF$ on $C$.
Let $\ov{p}\in \ov{C}$ be the image of $p$. 
Then $\pi^*\OO_{\ov{C}}(\ov{p})\simeq \OO_C(dp)$.
Thus, viewing a coherent right $\AA$-module $M$ as a coherent sheaf of $\OO$-modules, we
get a surjection of the form
$$\EE(-dmp)^{\oplus N}\simeq
\pi^*\EE(-m\ov{p})^{\oplus N}\to \pi_*\und{\Hom}(\EE,M)\ot \EE\to M.$$
Hence, the induced map of right $\AA$-modules
$$\EE(-dmp)^{\oplus N}\ot\AA\to M$$
is also surjective, which proves our claim.

Now we can repeat the well known argument
for the category of perfect $\OO$-modules (see e.g., the proof of \cite[Thm.\ 4]{Orlov-dim}): starting with any perfect
complex of $\AA$-modules $E$, we can find bounded above complex $P^\bullet$, where each $P^i$ is a direct sum of
modules of the form $\AA(-mp)$, and a quasi-isomorphism $P^\bullet\to E$. Now we consider brutal truncation
$\si^{\ge -n}P^\bullet$ for sufficiently large $n$. The cone of the composition
$$\si^{\ge -n}P^\bullet\to P^\bullet\to E$$ 
will be isomorphic in the derived category to $F[n+1]$, where $F$ is a coherent right $\AA$-module.
Furthermore, for sufficiently large $n$, we will have $\Hom(E,F[n+1])=0$, so we deduce that
$E$ is a direct summand in $\si^{\ge -n}P^\bullet$.
\ed

We say that a pairing
$$\FF\ot \GG\to \HH,$$
where $\FF$, $\GG$ and $\HH$ are coherent sheaves on a scheme, is {\it perfect on the left} (resp. {\it on the right})
if the induced map $\FF\to\und{\Hom}(\GG,\HH)$ (resp., $\GG\to\und{\Hom}(\FF,\HH)$) is an isomorphism.
We say that such a paring is {\it perfect in the derived category} (on the left or on the right) if the similar statements
hold with $\und{\Hom}$ replaced by $R\und{\Hom}$.

\begin{prop}\label{spherical-order-prop} 
Let $\AA$ be an order over an integral proper stacky curve $C$, which is smooth near all stacky points and
satisfies $H^0(C,\OO)=k$. 

\noindent
(i) $\AA$ is spherical if and only if $h^0(C,\AA)=1$ and there is an isomorphism of left
$\AA$-modules
\begin{equation}\label{sph-order-duality-eq}
\AA\simeq\und{\Hom}(\AA,\om_C),
\end{equation}
where $\om_C$ is the dualizing sheaf on $C$
(equivalently, one can ask for an existence of an isomorphism of right $\AA$-modules above).
In particular, $\AA$ is spherical if and only if $\AA^{op}$ is spherical.

Furthermore, if $\AA$ is spherical then $h^0(C,\AA)=h^1(C,\AA)=1$ and for a nonzero morphism $\tau:\AA\to\om_C$
(which is unique up to rescaling) the pairing 
$$\AA\ot\AA\to\om_C: (x,y)\mapsto \tau(xy)$$
is perfect in the derived category (on both sides).

\noindent
(ii) Assume now that $(C,p)$ is a neat pointed stacky curve, and $\AA$ is a spherical order over it,
such that $\AA|_p\simeq \rho_*\End(V)$. 
Let $g\in \End(V)$ be the element such
that the morphism 
$$\End(V)\simeq H^0(\AA(p)|_p)\rTo{\tau|_p} H^0(\om_C(p)|_p)\simeq k$$
is of the form $x\mapsto \tr(gx)$.  
Then $g$ is invertible, and the boundary homomorphism
$$\End(V)\simeq H^0(\AA(p)|_p)\to H^1(\AA)\simeq k,$$
associated with the exact sequence $0\to \AA\to \AA(p)\to \AA(p)|_p\to 0$,
is also of the form $x\mapsto \tr(gx)$, for an appropriate choice of an isomorphism $H^1(\AA)\simeq k$.
In addition, one has $h^1(C,\AA(p))=0$.
\end{prop}

\Pf . (i) For any vector bundles $\VV,\VV'$ over $C$ we have an isomorphism
$$\und{\Hom}_{\AA}(\AA\ot \VV,\AA\ot\VV')\simeq \und{\Hom}(\VV,\AA\ot\VV'),$$
whereas $\und{\Ext}^i$ vanish for $i>0$. Hence, we have isomorphisms
$$\Ext^i_\AA(\AA,\AA\ot\VV)\simeq H^i(C,\AA\ot\VV), \ \ \Ext^i_\AA(\AA\ot\VV,\AA)\simeq \Ext^i(\VV,\AA).$$
In particular, $\Ext^i(\AA,\AA)\simeq H^i(C,\AA)$. Furthermore, the canonical pairings
$$\Ext^{1-i}_\AA(\AA\ot\VV,\AA)\ot \Ext^i_\AA(\AA,\AA\ot\VV)\to\Ext^1_{\AA}(\AA,\AA)$$
get identified with the natural composed maps
\begin{equation}\label{A-Serre-pairings-eq}
\Ext^{1-i}(\VV,\AA)\ot H^i(\AA\ot\VV)\to H^1(\AA\ot\AA)\to H^1(\AA),
\end{equation}
where the second arrow is induced by the multiplication on $\AA$.
Since the modules of the form $\AA\ot\VV$ split generate $\Perf(\AA)$,
we deduce that $\AA$ is $1$-spherical as an object of $\Perf(\AA)$ (i.e., the order $\AA^{op}$ is spherical)
if and only if $h^0(\AA)=h^1(\AA)=1$ and all the pairings
\eqref{A-Serre-pairings-eq} are perfect.

Now assume that $\AA$ is $1$-spherical in $\Perf(\AA)$. The Serre duality on $C$ gives us perfect pairings
$$\Ext^{1-i}(\AA\ot\VV,\om_C)\ot H^i(\AA\ot\VV)\to H^1(\om_C).$$
In particular, we have a nonzero generator $\tau$ in the $1$-dimensional space $\Hom(\AA,\om_C)$
such that the induced map $H^1(\AA)\rTo{H^1(\tau)} H^1(\om_C)$ is an isomorphism.
It is easy to check that the map \eqref{A-Serre-pairings-eq} for $i=1$ fits into a commutative diagram
\begin{equation}\label{Serre-pairings-diagram}
\begin{diagram}
\Hom(\VV,\AA)\ot H^1(\AA\ot\VV)&\rTo{}& H^1(\AA)\\
\dTo{}&&\dTo{H^1(\tau)}\\
\Hom(\VV,\und{\Hom}(\AA,\om_C))\ot H^1(\AA\ot\VV)&\rTo{}& H^1(\om_C)
\end{diagram}
\end{equation}
where the bottom arrow is the Serre duality pairing combined with the isomorphism 
$$\Hom(\VV,\und{\Hom}(\AA,\om_C))\simeq \Hom(\AA\ot\VV,\om_C),$$ 
and the left vertical arrow comes from the morphism of left $\AA$-modules
$$\nu=\nu_\tau:\AA\to \und{\Hom}(\AA,\om_C): a\mapsto (x\mapsto \tau(xa)).$$
Since both horizontal arrows give perfect pairing and $H^1(\tau)$ is an isomorphism, we deduce that
the map
$$\Hom(\VV,\AA)\to \Hom(\VV,\und{\Hom}(\AA,\om_C)),$$
induced by $\nu$, is an isomorphism for all vector bundles $\VV$. It follows that $\nu$ is an isomorphism.

Note for an order $\AA$ with $h^0(\AA)=1$, the biduality morphism
$$\AA\to\und{\Hom}(\und{\Hom}(\AA,\om_C),\om_C)$$
is  an isomorphism. Indeed, since $\om_C$ is a dualizing sheaf on $C$, it suffices to check that
$\und{\Ext}^{>0}(\AA,\om_C)=0$. Equivalently, we have to check that
$$\Ext^i(\AA,\om_C\ot L^m)=H^0(C,\und{\Ext}^i(\AA,\om_C)\ot L^m)=0$$ 
for $i>0$ and $m\gg 0$, where $L$ be an ample line bundle on $C$ (see \cite[Prop.\ 6.9]{Hart}). 
By Serre duality, this reduces to the vanishing of $H^0(C,\AA\ot L^{-m})$, which is clear since every
global section of $\AA$ is a scalar multiple of the unit.
 
The above biduality statement easily implies that the bilinear pairing
$$\AA\ot \AA\to \om_C: a\ot a'\mapsto \tau(aa')$$
for some $\tau:\AA\to\om_C$
induces an isomorphism $\nu$ of left $\AA$-modules as above if and only if the corresponding morphism
of right $\AA$-modules 
$$\nu':\AA\to \und{\Hom}(\AA,\om_C): a\mapsto (x\mapsto \tau(ax))$$
is an isomorphism. 

Now let us start with an order $\AA$ such that $h^0(\AA)=1$ and there exists an isomorphism of left $\AA$-modules 
$\AA\simeq  \und{\Hom}(\AA,\om_C)$. Let $\tau\in \Hom(\AA,\om_C)$ be the element corresponding 
to the unit global section of $\AA$ under this isomorphism. Then the isomorphism is equal to $\nu_\tau$.
Since $\tau$ generates the space $\Hom(\AA,\om_C)$, by Serre duality, the map 
$H^1(\AA)\rTo{H^1(\tau)} H^1(\om_C)$ is an isomorphism.
Hence, the diagram \ref{Serre-pairings-diagram} implies that the pairing 
\eqref{A-Serre-pairings-eq} for $i=0$ is perfect. Now the pairing \eqref{A-Serre-pairings-eq} for $i=1$
fits into a similar diagram
\begin{equation}\label{Serre-pairings-bis-diagram}
\begin{diagram}
\Ext^1(\VV,\AA)\ot H^0(\AA\ot\VV)&\rTo{}& H^1(\AA)\\
\dTo{}&&\dTo{H^1(\tau)}\\
\Ext^1(\VV,\und{\Hom}(\AA,\om_C))\ot H^0(\AA\ot\VV)&\rTo{}& H^1(\om_C)
\end{diagram}
\end{equation}
Now we have an isomorphism
$$\Ext^1(\VV,\und{\Hom}(\AA,\om_C))\simeq H^1(\und{\Hom}(\VV,\und{\Hom}(\AA,\om_C)))\simeq
H^1(\und{\Hom}(\AA\ot\VV,\om_C))\simeq\Ext^1(\AA\ot\VV,\om_C)$$
since $\und{\Ext}^1(\AA\ot\VV,\om_C)=0$.
Thus, the pairing given by the bottom horizontal arrow in diagram \eqref{Serre-pairings-bis-diagram} is perfect, hence, 
so is the 
pairing given by the top horizontal arrow.

For the last assertion, we use the isomorphism 
$$H^1(C,\AA)^*\simeq H^0(C,\und{\Hom}(\AA,\om_C))\simeq H^0(C,\AA),$$
together with the vanishing of $\und{\Ext}^{>0}(\AA,\om_C)$ observed before (which holds since $h^0(C,\AA)=1$).

\noindent
(ii) We can think of $\rho_*\End(V)$ as the algebra $\End(V)\ot R_{\mu_d}$ in the category of $\mu_d$-representations,
where 
$$R_{\mu_d}=\bigoplus_{i=0}^{d-1}\chi^i$$ 
is the regular representation of $\mu_d$. The restriction
$\tau|_p:\AA|_p\to \om_C|_p$ 
can be viewed as a morphism of $\mu_d$-representations,
$$\End(V)\ot R_{\mu_d}\to \chi,$$
whose unique non-trivial component, $\End(V)\ot \chi\to \chi$, corresponds to the functional $x\mapsto\tr(gx)$ on $\End(V)$.
The fact that the induced pairing $\tau|_p(xy)$ on $\End(V)\ot R_{\mu_d}$ is nondegenerate easily implies that $g$ is
invertible.

Now let us consider the morphism of exact sequences induced by $\tau$,
\begin{diagram}
0&\rTo{}& \AA &\rTo{}& \AA(p)&\rTo{}&\AA(p)|_p&\rTo{} 0\\
&&\dTo{\tau}&&\dTo{\tau}&&\dTo{\tau|_p}\\
0&\rTo{}& \om_C &\rTo{}& \om_C(p)&\rTo{}&\om_C(p)|_p&\rTo{} 0
\end{diagram}
Passing to the corresponding exact sequences of cohomology, we get a commutative square
\begin{diagram}
H^0(\AA(p)|_p)&\rTo{}&H^1(\AA)\\
\dTo{\tau|_p}&&\dTo{\tau}\\
H^0(\om_C(p)|_p)&\rTo{}&H^1(\om_C)
\end{diagram}
in which the horizontal arrows are boundary homomorphisms.
The non-degeneracy of the pairing 
$$\Hom(\AA,\om_C)\ot H^1(\AA)\to H^1(\om_C)$$
implies that the right vertical arrow is an isomorphism. Since the bottom horizontal arrow is also an isomorphism,
we deduce that the top horizontal arrow can be identified with $\tau|_p$.
\ed

\subsection{Special spherical orders over the cuspidal cubic}\label{special-order-sec}

Let $C^\cusp$ be a cuspidal curve of arithmetic genus $1$ over a field $k$, $q$ a singular point, $p$ a smooth point.
Note that the normalization map is a homeomorphism, so we can identify $C^\cusp$ with $\P^1$ as a topological space.
We assume that $p$ corresponds to $\infty\in\P^1$, while $q$ corresponds to $0\in\P^1$.
For an $n$-dimensional vector space $V$ and $g\in\GL(V)$,
let us define an order $\AA^\cusp_g$ over $C^\cusp$ as the subsheaf of algebras $\AA^\cusp_g\sub\End(V)\ot\OO_{\P^1}$,
consisting of the elements that have an expansion
$a(z)=c\cdot I+a_1 z+\ldots$ near $0\in\P^1$, with $c\in k$ and $\tr(ga_1)=0$.

Note that functions on $C^{\cusp}$ have expansion $f(z)=a_0+a_2z^2+\ldots$ near $0\in\P^1$, so 
the order $\AA^\cusp_g$ indeed contains $\OO_{C^\cusp}$. It contains $\OO_{\P^1}$ precisely when $\tr(g)=0$.

Let us denote by $\tr_g:\AA^\cusp\to \OO_{C^\cusp}$ the homomorphism induced by the map 
$$\End(V)\ot\OO_{\P^1}\to\OO_{\P^1}:A\mapsto \tr(gA).$$

\begin{lem}\label{cusp-dual-lem} 
The $\OO$-bilinear form $\tr_g(aa')$ on $\AA^\cusp$ induces an isomorphism of right $\AA^\cusp$-modules
\begin{equation}\label{cusp-dual-eq}
\AA^\cusp\rTo{\sim} \HHom_{\OO_{C^\cusp}}(\AA^\cusp,\OO_{C^\cusp}): a\mapsto (a'\mapsto \tr_g(aa')).
\end{equation}
\end{lem}

\Pf . Away from $q$ this is clear, so it is enough to consider the completions of the stalks at $q$.
Since the localization $\widehat{\AA^\cusp}_q[z^{-1}]$ is just the matrix algebra over $k((z))$, we know that any
functional $\widehat{\AA^\cusp}_q\to k[[z]]$ has form $a\mapsto \tr_g(ab)$, for some $b=b_{-n}z^{-n}+b_{-n+1}z^{-n+1}+\ldots$
Now considering the condition that $\tr_g(ab)$ has to be a formal series of the form $c_0+c_2z^2+\ldots$, we easily deduce
that $b$ has to be in $\widehat{\AA^\cusp}_q$. 
\ed

\begin{rem} We also have an isomorphism of left $\AA^\cusp$-modules like \eqref{cusp-dual-eq}, given
by $a'\mapsto (a\mapsto \tr_g(aa'))$, which is in general
different from \eqref{cusp-dual-eq}.
\end{rem}

The order $\AA^\cusp_g$ is closely related to the algebra $\EE(V,g)$ defined by \eqref{E-V-eq}.

\begin{lem}\label{chi-lem} Under the construction of Sec.\ \ref{filtered-alg-sph-orders-sec}, 
the order $(\AA^\cusp_g)^{op}$ over $C^\cusp$ comes 
from the algebra $\EE(V,g)^{op}$, viewed as a filtered algebra. In particular, we have an equivalence
$$\mod-\AA^\cusp_g\simeq \qgr \EE(V,g)[t]$$
and an isomorphism of graded algebras
$$\EE(V,g)[t]\simeq\bigoplus_{n\in\Z} H^0(C^\cusp,\AA^\cusp_g(n)).$$
Also, the algebra $\EE(V,g)^{op}[t]$ is right Noetherian and satisfies $\chi$.
\end{lem}

\Pf . First, we need to identify $\AA^\cusp_g$ with the sheafification of $\EE(V,g)[t]$, viewed as a graded module over
$k[z^2,z^3][t]$. For this we observe that $\AA^\cusp_g$ is the subsheaf in $\End(V)\ot\OO_{\P^1}$, which
coincides with $\End(V)\ot\OO_{\P^1}$ over the complement to $q$. The same is true for the sheafification of $\EE(V,g)[t]$.
Hence, it is enough to compare the restrictions of the two sheaves to the affine open subset $C^{\cusp}\setminus\{p\}$
(i.e., the open subset $t\neq 0$). It remains to observe that 
$$H^0(C^{\cusp}\setminus\{p\},\AA^\cusp_g)=\EE(V,g)\sub\End(V)[z].$$ 
This proves the first first assertion. The other assertions follow from Lemma \ref{stacky-sheaf-lem} and 
Corollary \ref{Noeth-chi-cd-cor}.
\ed


\subsection{AS-Gorenstein condition over a field}\label{AS-Gor-sec}

For graded modules $M$ and $N$ over a graded algebra $B$ we use the notation
$$\und{\Ext}_B^i(M,N):=\bigoplus_{j\in\Z}\Ext^i_{\BB-\gr}(M,N(j))$$
where $\BB-\gr$ is the category of (left) $B$-modules.

Recall that a connected graded algebra $B$ over a field $k$ 
is called left {\it Artin-Schelter Gorenstein} ({\it AS-Gorenstein}) with the parameter $(d,m)$
if $B$ has a finite left injective dimension and $\und{\Ext}_B(k,B)$ is $1$-dimensional, concentrated in
degree $d$ and internal degree $m$. Similarly one defines the notion of right AS-Gorenstein.

\begin{prop}\label{special-Gor-prop} 
For any $g\in \GL(V)$ the algebra $\EE(V,g)[x]$ (where $\deg(x)=1$) is left and right AS-Gorenstein with the parameter $(2,0)$.
\end{prop}

\Pf . 
Let us set $B=\EE(V,g)[x]$. It is easy to see that $\EE(V,g)^{op}\simeq \EE(V^\vee,g^*)$, so it is enough to check that 
$B$ is right Gorenstein.
We have a natural identification of graded algebras
$$B=\bigoplus_{n\in\Z} H^0(C^\cusp,\AA^\cusp_g(n)),$$
where $\OO(1)=\OO(p)$
(see Lemma \ref{chi-lem}). Next, by Lemma \ref{cusp-dual-lem}, for any $n\in\Z$ the map
$$\AA^\cusp(n)\rTo{\sim} \HHom_{\OO_{C^\cusp}}(\AA^\cusp(-n),\OO_{C^\cusp}),$$
$a\mapsto (a'\mapsto \tr_g(aa'))$, is an isomorphism of right $\AA^{\cusp}$-modules. 
Since the $C^\cusp$ is Gorenstein with $\om_{C^\cusp}\simeq\OO_{C^\cusp}$,
by Serre duality, we deduce an isomorphism
$$\bigoplus_{n\in\Z} H^0(C^\cusp,\AA^\cusp(n))\rTo{\sim} \bigoplus_{n\in\Z} H^1(C^\cusp,\AA^\cusp(-n))^*$$
of right $B$-modules. In other words, we have an isomorphism of left $B$-modules
\begin{equation}\label{H1-A-cusp-eq}
\bigoplus_{n\in\Z} H^1(C^\cusp,\AA^\cusp(n))\simeq B^*,
\end{equation}
where $B^*$ is the restricted dual of $B$.

Next, let us consider the bar-resolution of $k$ by the complex of free right $B$-modules,
\begin{equation}\label{bar-resolution-eq}
\ldots\to B_+\ot B_+\ot B\to B_+\ot B\to B
\end{equation}
Localizing this sequence on $C^\cusp$ and twisting, we get for each $m\in\Z$ an exact sequence
of left $\AA^\cusp$-modules,
$$\ldots\to B_+\ot B_+\ot \AA^\cusp(m)\to B_+\ot \AA^\cusp(m)\to \AA^\cusp(m)\to 0$$
Let us consider the spectral sequence computing the hypercohomology of this exact complex, i.e., abutting to zero,
with the $E_1$-term given by the cohomology of the terms of this complex. Thus, the $E_1$-term has two rows,
corresponding to $H^0$ and $H^1$. The row of $H^0$'s is the degree $m$ component of the complex
\eqref{bar-resolution-eq}, which is exact for $m\neq 0$, and has $1$-dimensional cohomology in degree $0$ for $m=0$.
On the other hand, using the isomorphism \eqref{H1-A-cusp-eq} of left $B$-modules
we can identify the row of $H^1$'s with the degree $m$ component of the complex
\begin{equation}\label{B-B*-complex-eq}
\ldots\to B_+\ot B_+\ot B^*\to B_+\ot B^*\to B^*
\end{equation}
Since the spectral sequence abuts to zero, the row of $H^1$'s should be exact for $m\neq 0$ and has one-dimensional
cohomology in the term of degree $-2$ for $m=0$. 

Now the resolution \eqref{bar-resolution-eq} for $k$ as a right $B$-module shows that
the complex \eqref{B-B*-complex-eq} computes $\Tor^B(k,B^*)$, while its
restricted dual computes $\und{\Ext}^*_{B^{op}}(k,B)$.

In other words, $\und{\Ext}^*_{B^{op}}(k,B)$ has the one-dimenisonal cohomology
concentrated in cohomological degree $2$ and internal degree $0$. 

To conclude that $B$ is right Gorenstein it remains to check that
 $B$ has finite injective dimension as a right module over itself. To this end we use 
 \cite[Thm.\ 4.5]{Jorgensen} together with \cite[Thm.\ 6.3]{VdB-dc}. More precisely,
by Lemma \ref{chi-lem}, both $B$ and $B^{op}$ are right Noetherian, satisfy $\chi$, and their $\Proj^{nc}$ has finite 
cohomological dimension.
Hence, \cite[Thm.\ 6.3]{VdB-dc} gives existence of a balanced dualizing complex over $A$. Now
the same method as in \cite[Thm.\ 4.5]{Jorgensen} can be used to prove that $B$ has finite injective dimension
(see also \cite[Thm.\ 0.3(3)]{WZ} for a similar proof in the case of local rings).
\ed


\begin{prop}\label{AS-Gor-prop}
Let $k$ be a field.
For any filtered $k$-algebra $(A,F_\bullet A)$ satisfying \eqref{gr-E-eq} for some $g\in\GL_n(k)$, the Rees algebra $\RR(A)$ is left and right AS-Gorenstein with parameters $(2,0)$.
\end{prop}

\Pf . We observe that due to the nature of the Rees algebra construction 
there is a flat family $\RR_u(A)$ of graded algebras over $\A^1$ such that specializing $u$ to a nonzero value gives
an algebra isomorphic to $\RR(A)$, while $\RR_0(A)$ is $\gr^F(A)[x]\simeq \EE(V,g)[x]$. 
Namely, $\RR_u(A)=\RR(A)[u,x]/(t-xu)$, where $t\in \RR(A)_1$ is the canonical central element (here $\deg x=1$).

Therefore, since the algebra $\EE(V,g)[x]$ is AS-Gorenstein with parameters $(2,0)$ (see Proposition \ref{special-Gor-prop}),
we get that $\und{\Ext}^i(k,\RR(A))=0$ for $i\neq 2$ and $\und{\Ext}^2(k,\RR(A))$ is at most one-dimensional, concentrated
in the internal degree $0$.

Now let us check that $\und{\Ext}^*(k,\RR(A))$ cannot be entirely zero. 
Indeed, if it were zero then using \eqref{Hi-lim-formula},
we would get that all higher cohomology of $\OO(i)$ on $\Proj^{nc}$ vanishes.
Now we use the identification of $\qgr \RR(A)$ 
with coherent modules over the corresponding order $\AA$ (see Proposition \ref{Proj-order-prop})
and get a contradiction with the fact that $H^1(C,\AA)\simeq \Ext^1_{\AA^{op}}(\AA,\AA)$ is $1$-dimensional
(see Lemma \ref{stacky-sheaf-lem}(iv)).

Finally, by Corollary \ref{Noeth-chi-cd-cor}, $\RR(A)$ and $\RR^{op}$ are both Noetherian, satisfy $\chi$ and their
$\Proj^{nc}$ has finite cohomological dimension. Thus, as in the proof of Proposition \ref{special-Gor-prop} 
we can deduce that $\RR(A)$ has finite injective dimension, so $\RR(A)$ is AS-Gorenstein.
\ed

\section{From filtered algebras to spherical $n$-pairs}\label{application-pairs-sec}

\subsection{Gorenstein condition over a base ring}

Now we are going to switch to working over an arbitrary Noetherian commutative ring $S$.
Throughout this section we 
fix a filtered $S$-algebra $(A, F_\bullet A)$ equipped with an isomorphism \eqref{gr-E-eq} for some invertible
$g\in\End_S(V)\ot\LL$ 
(where $V\simeq S^n$),
and let $\RR=\RR(A)$ be the corresponding Rees algebra. 
Note that the graded components of $\RR$ are locally free $S$-modules of finite rank.

\begin{lem}
The algebra $\RR$ is right and left Noetherian. 
\end{lem}

\Pf . First, we note that the ring $\RR/(t)\simeq \EE(V,g)^{op}$ is right and left Noetherian, since $\EE(V,g)^{op}$
is finitely generated over its central subring $S[z^2]$.
Since $t$ is a regular central element, the assertion follows.
\ed

\begin{prop}\label{Ext-R-RR-prop}
Let $\RR^*$ be the restricted dual of $\RR$. 
Then for $i\neq 2,$ one has 
$$\Tor^\RR_i(S,\RR^*)=\Tor^\RR_i(\RR^*,S)=0, \ \ \und{\Ext}^i_{\RR}(S,\RR)=\und{\Ext}^i_{\RR^{op}}(S,\RR)=0,$$
and 
$\Tor^\RR_2(S,\RR^*)$, $\Tor^\RR_2(\RR^*,S)$, $\und{\Ext}^2_\RR(S,\RR)$ and 
$\und{\Ext}^2_{\RR^{op}}(S,\RR)$  
 are locally free $S$-modules of rank $1$,
concentrated in the internal degree $0$.
\end{prop}

\Pf . {\bf Step 1}.
First, we observe that the assertion is true when $S=k$ is a field. Indeed, this follows from Proposition \ref{AS-Gor-prop}
and from the duality between $\Tor^B(k,B^*)$ and $\und{\Ext}_{B^{op}}(k,B)$.

\noindent
{\bf Step 2}.
In the general case, since $\RR$ is Noetherian, we can find a free resolution 
$$\ldots\to \PP_2\to \PP_1\to \PP_0\to S,$$ 
where $\PP_i$ are free graded $\RR$-modules of finite rank. Let us set $Q_\bullet=(\PP_\bullet)\ot_{\RR}\RR^*$,
so that $H_i(Q_\bullet)\simeq\Tor^\RR_i(S,\RR^*)$, and the graded components of $Q_i$ are free $S$-modules of finite rank.
Let us assume that $S$ is local with the maximal ideal $M$, and set $k:=S/M$, $\RR_k:=\RR\ot_S k$. 
We are going to prove $H_iQ_\bullet=0$ for $i\neq 2$, while $H_2Q_\bullet\simeq S$ is concentrated in internal degree $0$.

We have
$$Q_\bullet\ot_S k\simeq (\PP_\bullet\ot_S k)\ot_{\RR_k} \RR_k^*.$$
Note that $\PP_\bullet\ot_S k$ is a free graded resolution of $k$ over $\RR_k$,
so that 
$$H_i(Q_\bullet\ot_S k)\simeq \Tor_i^{\RR_k}(k,\RR_k^*).$$
Note that we know from Step 1 that the latter spaces are zero for $i\neq 2$ and are isomorphic to $k$ in degree $0$
for $i=2$.

Now let us consider the third quadrant spectral sequence 
$$E^2_{p,q}=\Tor_{-q}^S(H_{-p}Q,k)\implies E^\infty_n=\Tor_{-n}^S(Q_\bullet,k),$$
with the differentials $d_r:E^r_{p,q}\to E^r_{p-r+1,q+r}$.
Note that since the terms of the complex $Q_\bullet$ are free $S$-modules we have
$\Tor_i^S(Q_\bullet,k)=H_i(Q_\bullet\ot_S k)$, so $E^\infty_n=0$ for $n\neq -2$, while $E^\infty_{-2}$ is $k$ in the internal
degree $0$. It follows that the term $E^2_{0,0}$ survives in the spectral sequence, so we get
$H_0Q\ot_S k=0$. Since $H_0Q$ has finitely generated graded components, by Nakayama lemma, this implies
that $H_0Q=0$. Therefore, $E^2_{0,q}=0$, so the term $E^2_{-1,0}$ survives, and we get $H_1Q\ot_S k=0$. 
Hence, $H_1Q=0$, and so $E^2_{-1,q}=0$. Thus, the terms $E^2_{-2,0}$ and $E^2_{-2,-1}$ survive, and we deduce
$H_2Q\ot_S k\simeq k$ (in degree $0$) and $\Tor_1^S(H_2Q,k)=0$. Hence, $H_2Q\simeq S$ (sitting in degree $0$).
Now the similar argument will prove by induction in $n\ge 3$ that $H_nQ=0$ (for the base we use the vanishing of
$E^2_{-2,q}$ for $q\le -1$). 

\noindent
{\bf Step 3}. For arbitrary $\RR$, since the construction of the complex $Q_\bullet$ is compatible with localization,
we deduce that $H_iQ_\bullet=0$ for $i\neq 2$, while $H_2Q_\bullet$ is a projective $S$-module of rank $1$ sitting
in internal degree $0$. This finishes the computation of $\Tor^\RR_i(S,\RR^*)$.
Since the complex $\Hom_S(Q_\bullet,S)$ computes $\Ext^i_\RR(S,\RR)$, and since $H_iQ_\bullet$ are projective,
we deduce that 
$$\Ext^i_\RR(S,\RR)\simeq \Hom_S(\Tor^\RR_i(S,\RR^*)),$$
and the assertion about $\Ext^*_\RR(S,\RR)$ follows. It remains to apply the same argument to $\RR^{op}$.
\ed

\subsection{Noncommutative projective scheme associated with a filtered algebra}
\label{noncomm-proj-sec}

As before, we consider the noncommutative projective scheme over $S$ 
associated with $\RR=\RR(A)$, i.e., the category
$\qgr \RR$, defined as the quotient of the category of graded finitely-generated right $\RR$-modules by the subcategory
of torsion modules.
We denote by $\OO(j)$ the object of $\qgr \RR$ corresponding to the module $\RR(j)$.
Recall also that $H^i(?):=\Ext^*_{\qgr \RR}(\OO,?)$.

\begin{prop}\label{proj-coh-prop} 
(i) In the category $\qgr \RR(A)$ one has
$$H^i(\OO(j))=0 \ \text{ for } i\neq 0,1; \ \ H^1(\OO(j))=0 \ \text{ for } j>0;$$
and there is a natural isomorphism of graded algebras
$$\bigoplus_j H^0(\OO(j))\simeq\RR(A).$$

\noindent
(ii) Let $F$ be the object of $\qgr \RR(A)$ corresponding to the graded right $\RR(A)$-module
$V[z]$, with the module structure induced by the homomorphism
$$\RR(A)\to \RR(A)/t\RR(A)\simeq \EE(V,g)^{op}\hra \End(V)^{op}[z].$$
Then the multiplication by $z$ induces an isomorphism $F\simeq F(1)$.
We have a natural exact sequence in $\qgr \RR(A)$:
\begin{equation}\label{qgr-F-eq}
0\to \OO(-1)\rTo{t} \OO\rTo{} V^\vee\ot F\to 0
\end{equation}
and canonical isomorphisms $H^0(F)\simeq V$, $\Ext^1(F,\OO)\simeq V^\vee$.
Also, $\Hom(F,F)=S\cdot\id_F$, $H^{>0}(F)=0$ and $\Ext^i(F,\OO)=0$ for $i\neq 1$.

\noindent
(iii) There exist canonical isomorphisms $H^1(\OO)\simeq \LL$ and $\Ext^1(F,F)\simeq S$
such that the compositions $\Ext^1(F,\OO)\ot H^0(F)\to H^1(\OO)$ and 
$H^0(F)\ot\Ext^1(F,\OO)\to\Ext^1(F,F)$
get identified the pairings $\lan v^*,gv\ran$ and $\lan v,v^*\ran$,
where $v\in V$, $v^*\in V^\vee$. Hence, the pair $(\OO,F)$ is an $n$-pair of $1$-spherical objects with the corresponding
element $g\in\End_S(V)\ot\LL$.

\noindent
(iv) For every $n\in\Z$ we have isomorphisms $\OO(n+1)\simeq T(\OO(n))$, where $T=T_F$ is the spherical twist associated with $F$. Hence, the graded algebra $\RR_{T,\OO}$ equipped with its natural central element of degree $1$ 
(see Theorem \ref{spherical-filtered-thm}) is isomorphic to $(\RR(A),t)$.
\end{prop}

\Pf . (i) Let us set $\RR=\RR(A)$. By \cite[Prop.\ 7.2]{AZ}, we have
$$H^i(\OO(j))=\lim_{m\to\infty} \Ext^{i+1}(\RR/\RR_{\ge m}, \RR(j)) \ \text{ for } i\ge 1,$$
and there is an exact sequence
$$0\to \tau(\RR(j))_0\to \RR_j\to H^0(\OO(j))\to\lim_{m\to\infty}\Ext^1(\RR/\RR_{\ge m},\RR(j))\to 0,$$
where $\tau(M)$ denotes the torsion submodule of $M$.
Note that $\tau(\RR(j))=0$ since $t$ is a nonzero divisor. On the other hand, by Proposition \ref{Ext-R-RR-prop},
we have
$$\Ext^{\neq 2}(S(m),\RR(j))=0, \ \ \Ext^2(S(m),\RR(j))=0 \ \text{ for } m\neq j.$$
Hence, the above exact sequence implies that the natural map 
$$\RR_j\to H^0(\OO(j))$$
is an isomorphism. The compatibility of these maps with the products is well known (see \cite[Thm.\ 4.5(2)]{AZ}.
Similarly, using the above formula for $H^i(\OO(j))$ with $i>0$ we see that $H^{>1}(\OO(j))=0$ and that
$H^1(\OO(j))=0$ for $j>0$.


\noindent
(ii) The multiplication by $z$ gives an injection $V[z]\to V[z](1)$ with finite-dimensional cokernel,
hence, it induces an isomorphism $F\simeq F(1)$.
The exact sequence \eqref{qgr-F-eq} is induced by the sequence of graded $\RR$-modules
$$0\to \RR(-1)\rTo{\cdot t} \RR\to \EE(V,g)\to 0$$
since we have $\EE(V,g)_{\ge 2}\simeq \End(V)[z]_{\ge 2}$.
Twisting \eqref{qgr-F-eq} by $(2)$ and using the identification $F\simeq F(2)$, we get an exact sequence
\begin{equation}\label{qgr-F-2-eq}
0\to \OO(1)\to \OO(2)\to V^\vee\ot F\to 0
\end{equation}
Using part (i) and a long exact sequence of cohomology we immediately deduce that $H^{>0}(F)=0$.
Note that the sequence \eqref{qgr-F-2-eq} also immediately implies that $\Hom(V^\vee\ot F,\OO)=0$, hence,
$\Hom(F,\OO)=0$.

Next, we claim that the natural map $V=V[z]_0\to H^0(F)$ is an isomorphism.
Note that there is a morphism of exact sequences
\begin{diagram}
0&\rTo{}& \RR_1 &\rTo{}& \RR_2 &\rTo{}& V^\vee\ot V&\rTo{}& 0\\
&&\dTo{}&&\dTo{}&&\dTo{}\\
0&\rTo{}& H^0(\OO(1)) &\rTo{}& H^0(\OO(2))&\rTo{}& V^\vee\ot H^0(F)&\rTo{}& 0
\end{diagram}
where the bottom row is obtained from \eqref{qgr-F-2-eq} by passing to $H^0$, 
and all vertical maps are the natural maps of the form
$M_0\to H^0(\wt{M})$, where $\wt{M}$ is the object of $\qgr \RR$ associated with a graded module $M$.
Note that the exactness of the bottom row follows the vanishing of $H^1(\OO(1))$ proved in (i). 
Since the two left vertical arrows are isomorphisms (again by (i)), we deduce that rightmost vertical
arrow is also an isomorphism, which proves our claim.

Another useful observation is that we have a morphism of functors $M\rTo{\cdot t} M(1)$, which vanishes on $F$.
Hence, the natural morphisms 
$$\Hom(G,F)\to \Hom(G(-1),F), \ \ \Ext^1(F,G)\to \Ext^1(F,G(1))$$
induced by $t$, are zero. Thus, applying the functor $\Hom(?,F)$ to the exact sequence \eqref{qgr-F-eq},
we deduce the isomorphism $\Hom(V^\vee\ot F,F)\rTo{\sim} H^0(F)$, induced by the projection $\OO\to V^\vee\ot F$.
Using the isomorphism $V\to H^0(F)$, this implies that $\Hom(F,F)=S\cdot\id$.

Now let us consider another twist of \eqref{qgr-F-eq}:
\begin{equation}\label{qgr-F-1-eq}
0\to \OO\to \OO(1)\to V^\vee\ot F\to 0
\end{equation}
The corresponding extension class is an element of $\Ext^1(V^\vee\ot F,\OO)\simeq V\ot \Ext^1(F,\OO)$,
so it gives a canonical morphism $V^\vee\to\Ext^1(F,\OO)$. In other words, this is precisely the connecting homomorphism
in the long exact sequence of $\Ext^*(F,?)$ applied to \eqref{qgr-F-1-eq}:
$$0=\Hom(F,\OO(1))\to V^\vee\ot\Hom(F,F)\to \Ext^1(F,\OO)\rTo{t} \Ext^1(F,\OO(1))\to\ldots$$
Since the map on $\Ext^1$ induced by $t$ is zero, we see that the map $V^\vee\to\Ext^1(F,\OO)$ is an isomorphism.

Finally, the vanishing of $\Ext^{\ge 2}(F,\OO)$ follows from the long exact sequence of $\Ext^*(?,\OO)$ applied to 
\eqref{qgr-F-eq}.

\noindent
(iii) First, the long exact sequence of cohomology applied to \eqref{qgr-F-1-eq} has form
$$0\to \RR_0\to \RR_1\to V^\vee\ot V\to H^1(\OO)\to H^1(\OO(1))=0$$
Furthermore, the induced map $\RR_1/\RR_0\to V^\vee\ot V$ is exactly the embedding $\EE(V,g)_1\sub \End(V)$,
so its image is the $S$-submodule $\End(V)_g\sub \End(V)$. Thus, there is a unique isomorphism $H^1(\OO)\simeq \LL$,
such that the above map $V^\vee\ot V\to H^1(\OO)$ gets identified with $v^*\ot v\mapsto \lan v^*,gv\ran$.
Note that this map corresponds to the composition $\Ext^1(F,\OO)\ot H^0(F)\to H^1(\OO)$ using the identifications
$\Ext^1(F,\OO)\simeq V^\vee$, $H^0(F)\simeq V$ defined in (ii).

Next, we observe that since the map $\Ext^1(F,\OO(-1))\rTo{t}\Ext^1(F,\OO)$ is zero,
the long exact sequence of $\Ext^*(F,?)$ associated with \eqref{qgr-F-eq} gives an isomorphism
\begin{equation}\label{FO-FF-eq}
\Ext^1(F,\OO)\rTo{\sim} \Ext^1(F,V^\vee\ot F)\simeq V^\vee\ot\Ext^1(F,F).
\end{equation}
In particular, this implies that $\Ext^1(F,F)$ is a locally free $S$-module of rank $1$.
We have a split exact sequence
$$0\to \End_0(V)\ot F\to V\ot V^\vee\ot F\rTo{\tr\ot\id_F} F\to 0$$
Furthermore, the natural map $V\ot \OO\to F$ corresponding to the indentification $V=H^0(F)$, is the composition
$$V\ot \OO\rTo{\id_V\ot p} V\ot V^\vee\ot F\rTo{\tr\ot\id_F} F,$$
where $p:\OO\to V^\vee\ot F$ is the map from the sequence \eqref{qgr-F-eq}.
Thus, the map $\Ext^1(F,V\ot \OO)\to \Ext^1(F,F)$ can be identified with the composition
$$\Ext^1(F,V\ot \OO)\rTo{\sim} \Ext^1(F,V\ot V^\vee\ot F)=V\ot V^\vee\ot \Ext^1(F,F)\rTo{\tr\ot\id}\Ext^1(F,F),$$
where the first arrow is obtained from \eqref{FO-FF-eq} by tensoring with $V$.

Thus, we deduce the surjectivity of the composition map 
\begin{equation}\label{composition-map-ExtFF}
V\ot V^\vee\simeq \Hom(\OO,F)\ot \Ext^1(F,\OO)\to \Ext^1(F,F).
\end{equation}
We claim that $\End_0(V)\sub \End(V)=V\ot V^\vee$ is contained in the kernel of this map.
Indeed, it is enough to check that for any $v\in V\simeq\Hom(\OO,V)$ and any $v^*\in V^\vee\simeq\Ext^1(F,\OO)$,
such that $\lan v^*,v\ran=0$, the composition of $v^*$ and $v$ in $\Ext^1(F,F)$ vanishes.
Let us consider the push-out of \eqref{qgr-F-1-eq} by $v:\OO\to F$:
\begin{equation}\label{Ev-diagram}
\begin{diagram}
0&\rTo{}& \OO&\rTo{t}& \OO(1)&\rTo{}& V^\vee\ot F&\rTo &0\\
&&\dTo{v}&&\dTo{}&&\dTo{\id}\\
0&\rTo{}& F&\rTo{}& E_v&\rTo{}& V^\vee\ot F&\rTo &0
\end{diagram}
\end{equation}
Next, let us compute explicitly the subobject 
$$E'_v:=\ker(t:E_v\to E_v(1))\sub E_v.$$ 
We can represent $E_v$ by the graded $\RR$-module 
$$V[z]_{\ge 1}\oplus \RR(1)_{\ge 1}/\{(-v*r, tr) \ |\ r\in \RR_{\ge 1}\},$$
so that the embedding $F\to E_v$ corresponds to the embedding of the summand $V[z]_{\ge 1}$.
Here for $r\in \RR_m$ we denote by $v*r\in V\cdot z^m$ the result of the right action of the image of 
$r$ in $\RR_m/\RR_{m-1}\sub\End(V)^{op}\cdot z^m$ on $v$.
Hence, $E'_v$ corresponds to the submodule of pairs $(x,r)$ such that $v*r=0$.
This easily implies that the image of the projection $E'_v\to \OO(1)/\OO\cdot t\simeq V^\vee\ot F$ 
coincides with $\lan v\ran^\perp \ot F$. Thus, the bottom sequence in diagram \eqref{Ev-diagram} contains 
as a subsequence the exact sequence
\begin{equation}\label{E'v-eq}
0\to F\to E'_v\to \lan v\ran^\perp \ot F\to 0
\end{equation}
of objects in the subcategory $\ker(t)\sub \qgr \RR$. Note that the latter subcategory is naturally identified with 
$\qgr \End(V)[z]$ and that $F$ is a projective object
in this subcategory. Hence, the sequence \eqref{E'v-eq} splits which proves the required vanishing in $\Ext^1(F,F)$.

It follows that the composition map \eqref{composition-map-ExtFF} factors through a surjective map
$$S\simeq \End(V)/\End_0(V)\to \Ext^1(F,F).$$ 
Since $\Ext^1(F,F)$ is a locally free $S$-module of rank $1$, this map is in fact
an isomorphism.

\noindent
(iv) As we have seen above, the exact sequence \eqref{qgr-F-1-eq} induces an isomorphism $V^\vee\to \Ext^1(F,\OO)$.
Hence, it gives a canonical isomorphism 
$$\OO(1)\simeq T_F(\OO).$$
On the other hand, for any $n\in\Z$ we have an isomorphism $F(n)\simeq F$. Hence, applying the autoequivalence
$M\mapsto M(n)$ to the above isomorphism we get an isomorphism
$$\OO(n+1)\simeq T_{F(n)}(\OO(n))\simeq T_F(\OO(n)).$$
\ed

\begin{cor}\label{H1-bimodule-cor} 
There exists a graded automorphism $\phi$ of $\RR$, such that $\phi(t)=t$, the induced automorphism
of $\RR/(t)\simeq \EE(V,g)^{op}$ is $\Ad(g^{-1})$, and there is an isomorphism of graded $\RR-\RR$-bimodules
$$\und{H}^1(\OO):=\bigoplus_{i\in\Z} H^1(\OO(i))\simeq (\sideset{_{\id}}{_{\phi}}\RR)^*\ot_S \LL.$$
\end{cor}

\Pf . This follows by combining Proposition \ref{proj-coh-prop} with Proposition \ref{bimodule-prop}.
\ed

Next, we are going to construct a certain short exact sequence
in $\qgr \RR(A)$ (it will be used in Sec.\ \ref{triple-prod-sec} 
to characterize $A_\infty$-structures). Namely, let us consider an element $hz\in\End_g(V)[z]\simeq F_1A/F_0A$, where $h$
is any invertible element of $\End(V)\ot\LL^{-1}$, such that $\tr(gh)=0$. Note that we get 
take $h=g^{-1}h_0$ where $h_0$ is an element of $\GL(V)$ with $\tr(h_0)=0$ (it exists since we assume
that $V\simeq S^n$ and $n\ge 2$).
Let $\wt{h}\in F_1A\ot\LL^{-1}$ be any lifting of $hz$ to $F_1A\ot\LL^{-1}=\RR_1(A)\ot\LL^{-1}$.
The right multiplication by $hz$ induces an injective map $\EE(V,g)\to \EE(V,g)(1)\ot\LL^{-1}$ (which we can view
as a map of right $\RR(A)$-modules) with finite-dimensional cokernel,
hence, an isomorphism $V^\vee\ot F\to V^\vee\ot F(1)\ot\LL^{-1}$ in $\qgr \RR(A)$.
Since $t\in \RR_1(A)$ is central, we have a commutative diagram in $\qgr \RR(A)$ with exact rows 
\begin{diagram}
0&\rTo{}&\OO&\rTo{t\cdot}& \OO(1)&\rTo{}& V^\vee\ot F(1)&\rTo{}&0\\
&&\dTo{\wt{h}\cdot}&&\dTo{\wt{h}\cdot}&&\dTo{\cdot hz}\\
0&\rTo{}&\OO(1)\ot\LL^{-1}&\rTo{t\cdot}& \OO(2)\ot\LL^{-1}&\rTo{}& V^\vee\ot F(2)\ot\LL^{-1}&\rTo{}&0
\end{diagram}
This implies exactness of the complex
\begin{equation}\label{g-exact-seq}
0\to \OO\rTo{\a}\OO(1)\oplus \OO(1)\ot\LL^{-1}\rTo{\b}\OO(2)\ot\LL^{-1}\to 0
\end{equation}
where 
\begin{equation}\label{a-b-comp-eq}
\a=(t\cdot,\wt{h}\cdot), \ \ \b=(\wt{h}\cdot,(-t)\cdot).
\end{equation}

\begin{lem}\label{Ext1-h-lem} 
Under the isomorphism
$$\Ext^1(\OO(2),\OO)\ot\LL\rTo{\sim} \Hom_S(\RR_2,\LL)\ot\LL:c\mapsto (r\mapsto c\cdot r),$$
where we use the identification of $H^1(\OO)\simeq \LL$ from Proposition \ref{proj-coh-prop}(iii), 
the class $\ga\in \Ext^1(\OO(2),\OO)\ot\LL$ of extension \eqref{g-exact-seq} corresponds to the functional
$$\RR_2\to\LL^2: r\mapsto \tr(\frac{p_2(r)}{z^2}h^{-1}g),$$
where $p:\RR\to \RR/(t)\simeq \EE(V,g)^{op}$ is the natural projection (note that $p_2(r)$ is an
element of $\End(V)z^2$).
\end{lem}

\Pf . Let $\ga'\in\Ext^1(V^\vee\ot F,\OO(-1))\simeq \Ext^1(V^\vee\ot F(1),\OO)$ 
be the class of the extension \eqref{qgr-F-eq}, and let $\ov{p}:\OO\to V^\vee\ot F$
be the natural projection. We claim that $\ga$ is equal to the composition
$$\OO(2)\ot\LL^{-1}\rTo{\ov{p}} V^\vee\ot F(2)\ot\LL^{-1}\rTo{\cdot (hz)^{-1}} V^\vee\ot F(1)\rTo{\ga'}\OO[1].$$
Indeed, this follows immediately from the commutative diagram in which rows and columns extend to short exact sequences,
\begin{diagram}
\OO&\rTo{t}&\OO(1)&\rTo{\ov{p}}& V^\vee\ot F(1)\\
\uTo{=}&&\uTo{}&&\uTo{(\cdot (hz)^{-1})\circ\ov{p}}\\
\OO&\rTo{\a}&\OO(1)\oplus\OO(1)\ot\LL^{-1}&\rTo{\b}&\OO(2)\ot\LL^{-1}\\
\uTo{}&&\uTo{}&&\uTo{-t}\\
0&\rTo{}&\OO(1)\ot\LL^{-1}&\rTo{=}&\OO(1)\ot\LL^{-1}
\end{diagram}
Next, for any $r\in \RR_m(A)$ the composition $\OO\rTo{r} \OO(m)\rTo{\ov{p}} V^\vee\ot F(m)$
corresponds to an element $p(r)\in \EE(V)_m\simeq \End(V)z^m\simeq H^0(V^\vee\ot F(m))$.
Finally, the extension class $\ga'$ corresponds to the identity element in $V\ot V^\vee$ under an isomorphism
$\Ext^1(V^\vee\ot F(1),\OO)\simeq V\ot V^\vee$, and the composition $\Ext^1(F,\OO)\ot\Hom(\OO,F)\to H^1(\OO)$ 
is given by $v^*\ot v\mapsto \lan v^*,gv\ran$. This easily implies the assertion.
\ed

\section{Characterization of $A_\infty$-structures}\label{char-sec} 

\subsection{Minimally non-formal algebras}\label{min-non-form-alg-sec}

Let $A$ be a graded algebra over a commutative ring $S$.
We can consider structures of ($S$-linear) minimal $A_\infty$-algebras on $A$ extending the given $m_2$,
up to (strict) gauge equivalences. 
If every such structure is gauge equivalent to the one with $m_i=0$ for $i>2$, then $A$ is called {\it intrinsically formal}.
We are going to describe a class of algebras $A$ for which gauge equivalence classes of minimal $A_\infty$-structures
are classified by elements of a certain $S$-module.

Recall that the set of minimal $A_\infty$-structures on $A$ is governed by the Hochschild cohomology $HH^*(A)$ of the underlying graded associative algebra.
More precisely, if we already have products $m_i$ for $i\le n-1$, forming an $A_{n-1}$-structure, then 
the set of $m_n$ extending these to an $A_n$-structure is a torsor over $HH^2(A)_{2-n}$ (we follow the grading convention
in which $m_n$ is a Hochschild $2$-cochain of internal degree $2-n$). 
So the vanishing of $HH^2(A)_{<0}$ implies intrinsic formality. The next simplest case after that which occurs
in some situations is when $HH^2(A)_{2-n}$ is nonzero for the unique value $n=d\ge 3$.
In this case every $A_\infty$-structure is gauge equivalent to the one with $m_i=0$ for 
$2<i<d$. Furthermore, $m_d$, calculated for such a representative, gives a class in $HH^2_{2-d}$,
which uniquely determines the gauge equivalence class of the $A_\infty$-structure. 


\begin{defi}\label{ABM-def}
Let $B=\bigoplus_{n\ge 0}B_n$ be a graded $S$-algebra with $B_0=S$,
and let $M=\bigoplus_n M_n$ be a graded $B-B$-bimodule. 
We define the bigraded algebra $A(B,M,d)$ to be $B\oplus M$, where $M\cdot M=0$, with the natural internal grading
and the homological grading given by $\deg(S)=0$, $\deg(M)=d$.
\end{defi}

The next theorem is a slight generalization of the results in \cite[Sec.\ 3.1]{P-ext}.
As in \cite[Sec.\ 3.1]{P-ext}, for graded $B-B$-bimodules $M_1,\ldots,M_n$ let us consider the
bar-complex
$$\Bar^\bullet(M_1,\ldots,M_n):=M_1\ot_S T(B_+)\ot_S M_2\ot_S\ldots\ot_S T(B_+)\ot_S M_n,$$
where $T(B_+)$ is the tensor algebra of $B_+=\bigoplus_{n\ge 1} B_n$ as an $S$-module.
The grading is given by
$$\Bar^{-m}(M_1,\ldots,M_n):=\bigoplus_{m_1+\ldots+m_{n-1}=m}
M_1\ot_S T^{m_1}(B_+)\ot_S M_2\ot_S\ldots\ot_S T^{m_{n-1}}(B_+)\ot_S M_n.$$
Note that the cohomology $H^{-m}$ of the complex $\Bar^\bullet(S,M)$ (resp., $\Bar^\bullet(M,S)$) 
is isomorphic to $\Tor^B_m(S,M)$ (resp., $\Tor^B_m(M,S)$). 

\begin{thm}\label{HH-gen-thm} Assume that $M$ bounded above, i.e., $M_n=0$ for $n>n_0$, and that
$P^l:=\Tor_*^B(S,M)$ and $P^r:=\Tor^B_*(M,S)$ are both finitely generated projective $S$-modules, 
concentrated in degree $d+1$ and the internal degree $0$. Let us fix an embedding 
$$P^r\rTo{\varphi} M\ot_S T^{d+1}(B_+)$$
inducing the isomorphism of $P^r$ with the cohomology of $\Bar^\bullet(M,S)$.
Then for the algebra $A=A(B,M,d)$ one has
$HH^i_{-md}(A)=0$ for $m\ge 1$ and $i<2m$, and there is an embedding
$$HH^2_{-d}(A)\hra\Hom_S(P^r,S),$$
induced by the evaluation of a Hochschild cochain on the image of $\varphi$.
Here the lower index denotes the grading on Hochschild cohomology induced by the homological grading on $A$.
\end{thm}

\Pf . 
Below we consider graded $B-B$-bimodules $M_i$ which are equipped with a pair of isomorphisms 
$l:M_i\to M\ot_S P_i$, $r:M_i\to P'_i\ot_S M$, for some finitely generated projective $S$-modules $P_i$, $P'_i$,
where $l$ (resp., $r$) is compatible with the left (resp., right) graded 
$B$-module structures.

\noindent
{\bf Step 1}. $H^i(\Bar^\bullet(M_1,M_2))=0$ for $i\neq -d-1$ and we have an isomorphism of right 
$B$-modules
$$H^{-d-1}\Bar^\bullet(M_1,M_2)\simeq P'_1\ot_S P^r\ot_S M$$
and an isomorphism of left $B$-modules
$$H^{-d-1}\Bar^\bullet(M_1,M_2)\simeq M\ot_S P^l\ot_S P_2.$$

To prove the first assertion we consider the spectral sequence associated with the filtration on $\Bar^\bullet(M_1,M_2)$
induced by the $\Z$-grading on $M_2$. The corresponding $E_1$-term is
$$E_1\simeq H^\bullet\Bar^\bullet(M_1,S)\ot_S M_2\simeq P'_1\ot_S P^r\ot_S M_2.$$
Hence, the spectral sequence degenerates and we obtain the first assertion.
Similarly, the second assertion is obtained by considering the spectral sequence associated with the filtration induced
by the $\Z$-grading on $M_1$.

\noindent
{\bf Step 2}. $H^i\Bar^\bullet(M_1,\ldots,M_n)=0$ for $i>-(n-1)(d+1)$, $n\ge 2$.

We use induction on $n$. For $n=2$ the assertion follows from Step 1. Now for $n>2$ we can consider 
$\Bar^\bullet(M_1,\ldots,M_n)$ as the total complex associated with a bicomplex, by considering the bigrading
given by the sums of the tensor degrees in even and odd factors $T(B_+)$. This leads to a spectral sequence
abutting to $H^*\Bar^\bullet(M_1,\ldots,M_n)$ with the $E_1$-term
$$H^*\Bar^\bullet(M_1,M_2)\ot_S T(B_+)\ot_S H^*\Bar^\bullet(M_3,M_4)\ot_S T(B_+)\ot_S\ldots$$
where the last tensor factor is either $M_n$ or $H^*\Bar^\bullet(M_{n-1},M_n)$.
Thus, the $E_1$-term is isomorphic to the complex of the form
$$\Bar^\bullet(M'_1,\ldots,M'_{n'})[(n-n')(d+1)]$$
where $M'_1=H^*\Bar^\bullet(M_1,M_2)$, $M'_2=H^*\Bar^\bullet(M_3,M_4)$, etc., satisfy the same assumptions as $(M_i)$
by Step 1. Applying the induction assumption we deduce the result.

\noindent
{\bf Step 3}. $H^i\Bar^\bullet(S,M_1,M_2,\ldots,M_n,S)=0$ for $i>-n(d+1)$, and we have isomorphisms
$$H^{-d-1}\Bar^\bullet(S,M_1,S)\simeq P^l\ot_S P_1\simeq P'_1\ot_S P^r$$
of graded $S$-modules. 

Consider first the complex $\Bar^\bullet(S,M_1,S)=T(B_+)\ot_S M_1\ot_S T(B_+)$. We can view
it naturally as the total complex associated with a bicomplex (where the bigrading is induced by two tensor degrees).
Considering the corresponding two spectral sequences we immediately deduce the vanishing of
$H^i\Bar^\bullet(S,M_1,S)$ for $i>-d-1$ and get the required identifications of $H^{-d-1}\Bar^\bullet(S,M_1,S)$. 

Now we use induction on $n$. As before, we equip $\Bar^\bullet(S,M_1,\ldots,M_n,S)$ with the bigrading
using sums of tensor degrees in even and odd factors $T(B_+)$.
Thus, we get a spectral sequence of a bicomplex abutting to 
$H^*\Bar^\bullet(S,M_1,\ldots,M_n,S)$ with the $E_1$-term either of the form
\begin{align*}
&T(B_+)\ot H^*\Bar^\bullet(M_1,M_2)\ot T(B_+)\ot\ldots H^*\Bar^\bullet(M_{n-1},M_n)\ot T(B_+)=\\
&\Bar^\bullet(S,M'_1,\ldots,M'_{n/2},S)[n(d+1)/2]
\end{align*}
if $n$ is even, or of the form
\begin{align*}
&H^*\Bar^\bullet(S,M_1)\ot T(B_+)\ot H^*\Bar^\bullet(M_2,M_3)\ot T(B_+)\ot\ldots\ot 
H^*\Bar^\bullet(M_{n-1},M_n)\ot T(B_+)\\
&\simeq P^l\ot_S P_1\ot_S \Bar^\bullet(S,M'_1,\ldots,M'_{(n-1)/2},S)[(n+1)(d+1)/2]
\end{align*}
if $n$ is odd. In both cases the required cohomology vanishing follows from the induction assumption.

\noindent
{\bf Step 4.} Let us fix $m\ge 1$ and denote by 
$$C^i_{-md}\subset \Hom_S(A_+^{\ot i}, A)$$
the submodule of degree $-md$ with respect to the homological grading on $A$ and of degree $0$ with respect
to the internal grading on $A$. Note that the Hochschild differential maps $C^i_{-md}$ to $C^{i+1}_{-md}$,
and $HH^i_{-md}(A)$ is the $(i+md)$-th cohomology of the complex $C^{\bullet}_{-md}$.
We have an exact sequence of complexes
$$0\to C^\bullet_{-md}(M)\to C^\bullet_{-md}\to C^\bullet(B)\to 0$$
where $C^i_{-md}(M)\sub C^i_{-md}$ (resp., $C^i_{-md}(B)\sub C^i_{-md}$) 
consists of maps $A_+^{\ot i}\to M$ (resp., $A_+^{\ot i}\to B$). 
Note that $C^i_{-md}(B)$ (resp., $C^i_{-md}(M)$) consists of $S$-linear maps
$$[T(B_+)\ot_S M\ot_S T(B_+)\ot_S\ldots\ot_S M\ot T(B_+)]_i\to B$$
(resp.,
$$[T(B_+)\ot_S M\ot_S T(B_+)\ot_S\ldots\ot_S M\ot T(B_+)]_i\to M)$$
preserving the internal grading, where there are $m$ (resp., $m+1$) factors of $M$ in the source
and the index $i$ refers to the total number of tensor factors (of $M$ and $B_+$).

We claim that 
$$H^i(C^\bullet_{-md}(M))=H^i(C^\bullet_{-md}(B))=0 \text{  for  } i<m(d+2)$$ 
and in addition, 
$$H^{d+2}(C^\bullet_{-d}(M))=0$$
and there is an embedding 
$$H^{d+2}(C^\bullet_{-d}(B))\hra \Hom_S(P^r,S)$$
induced by the evaluation on the image of $\varphi$.

For the proof, let us consider the decomposition 
$C^\bullet_{-md}(B)=\prod_{j\ge 0}C^\bullet_{-md}(B_j)$,
where $C^\bullet_{-md}(B_j)\sub C^\bullet_{-md}(B)$ denote the maps with the image contained in $B_j$.
Let us consider the corresponding decreasing filtration on $C^\bullet_{-md}(B)$. 
Note that the corresponding associated graded complex is
$$\bigoplus_{j\ge 0} \Hom_S(\Bar^\bullet(S,(M)^m,S)_j,B_j)[-m],$$
where the lower index denotes the internal grading (and we grade the dual complex using the convention
$\Hom(K^\bullet,?)^i=\Hom(K^{-i},?)$).
By Step 3, we have the vanishing
$$H^i\Hom_S(\Bar^\bullet(S,(M)^m,S)_j,B_j)=0 \text{  for  } i<m(d+1),$$
$$H^{d+1}\Hom_S(\Bar^\bullet(S,M,S)_j,B_j)=0 \text{  for  } j>0,$$
and in addition,
$$H^{d+1}\Hom_S(\Bar^\bullet(S,M,S)_0,B_0)\simeq \Hom_S(P^r,S).$$
Hence, using \cite[Lem.\ 3.2]{P-ext},
we obtain the required statements about cohomology of $C^\bullet_{-md}(B)$.

To deal with the cohomology of $C^\bullet_{-md}(M)$ we consider the decreasing filtration on
$C^\bullet_{-md}(M)$ associated with the grading induced by the sum of the tensor degrees in the first and last
factors of the tensor product $T(B_+)\ot M\ot T(B_+)\ot\ldots\ot M\ot T(B_+)$.
The associated graded complex can be identified with
$$\Hom_{\gr-S-\mod}(T(B_+)\ot_S \Bar^\bullet((M)^{m+1})\ot_S T(B_+), M)[-m-1].$$
Thus, the required vanishing follows from Step 2.
\ed


\begin{cor}\label{ainf-min-nf-cor}
Under the conditions of Theorem \ref{HH-gen-thm}, the map $(m_\bullet)\mapsto m_{d+2}$ induces a bijection between
the set of gauge-equivalence classes of minimal ($S$-linear) $A_\infty$-structures
on $A=A(B,M,d)$ with given $m_2$ and the $S$-module $HH^2_{-d}(A)$.  
\end{cor}

\Pf . Since the grading on $A$ is concentrated in degrees $0$ and $d$, the only potentially nonzero higher products
are $m_n$ with $n\equiv 2\mod (d)$. Thus, the $A_\infty$-identities imply that $m_{d+2}$ is a Hochschild cocycle.
Now the vanishing of $HH^2_{-md}$ for $m\ge 2$ implies by the standard argument that
the gauge equivalence class of $(m_\bullet)$ is determined by the cohomology class of $m_{d+2}$.
Furthermore, since $HH^3_{-md}(A)=0$ for $m\ge 2$, by \cite[Lem.\ 3.1.2(ii)]{P-ainf}, every Hochschild cocycle
$m_{d+2}$ extends to an $A_\infty$-structure on $A$.
\ed



\begin{cor}\label{m3-1-cor} 
Under the conditions of Theorem \ref{HH-gen-thm} with $d=1$, assume that we have a finitely generated projective $S$-module
$\QQ$ and a locally free $S$-module $\LL$ of rank $1$, maps of $S$-modules
$\a:\QQ\to B_{n_1}$, $\b:\QQ^\vee\to B_{n_2}\ot\LL^{-1}$ and an element $\ga\in M_{-n_1-n_2}\ot\LL$, 
where $n_1>0$, $n_2>0$, such that
$m_2(\b\ot \a)(\id_\QQ)=0$ in $B_{n_1+n_2}\ot\LL^{-1}$ and $m_2(\ga\ot\b)=0$ in $\QQ^\vee\ot M_{-n-1}$, where $m_2$
is induced by the product in $B$ and the $B$-bimodule structure on $M$. 
Assume in addition that the $S$-module $P^r=\Tor_2(M,S)_0$ is
locally free of rank $1$.
Then up to a gauge equivalence, there exists at most one minimal
$S$-linear $A_\infty$-structure on $A=A(B,M,1)$ with the given $m_2$ and with 
\begin{equation}\label{m3-ga-b-a-eq}
m_3((\ga\ot\b\ot\a)(\id_\QQ))=1.
\end{equation}
Furthermore, if such an $A_\infty$-structure exists then the gauge equivalence classes of all minimal $A_\infty$-structures
with the given $m_2$ are classified by elements of $S$ via the map 
$$(m_\bullet)\mapsto m_3((\ga\ot\b_i\ot\a)(\id_\QQ)).$$
\end{cor}

\Pf . As in Corollary \ref{ainf-min-nf-cor}, we see that a
minimal  $A_\infty$-structure on $A$ with the given $m_2$ is determined by the class of $m_3$ in $HH^2_{-1}(A)$,
and we have an embedding 
$$HH^2_{-1}(A)\hra \Hom_S(P^r,S).$$

Next, we have a morphism $S\to P^r=\Tor_2(M,S)_0$ given by the
cycle $(\ga\ot\b_i\ot\a)(\id_\QQ)$ in the complex $\Bar^\bullet(M,S)$. 
Assume that there exists an $A_\infty$-structure on $A$ with the given $m_2$ such that
\eqref{m3-ga-b-a-eq} holds.
Then the composition
$$HH^2_{-1}(A)\to \Hom_S(P^r,S)\rTo{f\mapsto f((\ga\ot\b\ot\a)(\id_\QQ))} S$$
sends the class of $m_3$ to $1$. It follows that the morphism
$\Hom_S(P^r,S)\to S$ is surjective. Since $P^r$ is locally free of rank $1$, we deduce that this morphism is an isomorphism.
Hence, the map $HH^2_{-1}(A)\to S$, given by the evaluation on $(\ga\ot\b\ot\a)(\id_\QQ)$, is an isomorphism,
which implies our assertion.
\ed


\subsection{Triple product calculation}\label{triple-prod-sec}


Now we return to the situation of Section \ref{noncomm-proj-sec},
so $S$ is a Noetherian commutative ring,
$(A, F_\bullet A)$ a filtered $S$-algebra equipped with an isomorphism \eqref{gr-E-eq} for some invertible
$g\in\End_S(V)\ot\LL$ 
(where $V\simeq S^n$), and $\RR=\RR(A)$ be the corresponding Rees algebra. 

We would like to 
show that Corollary \ref{m3-1-cor} is applicable to 
minimal $S$-linear $A_\infty$-structures on the algebra 
$A(\RR, (\sideset{_{\id}}{_{\phi}}\RR)^*\ot\LL),1)$ (see Def.\ \ref{ABM-def}), where $\phi$ is the automorphism
of $\RR$ defined in Corollary \ref{H1-bimodule-cor}. 
Recall that we have an isomorphism of bimodules over $\RR$,
$\und{H}^1(\OO)\simeq (\sideset{_{\id}}{_{\phi}}\RR)^*\ot\LL$.

As in Lemma \ref{Ext1-h-lem}, let us fix 
 an invertible element $h\in \End(V)\ot\LL^{-1}$ such that $\tr(gh)=0$ and its lifting 
$\wt{h}\in F_1 A\ot\LL^{-1}=\RR_1\ot\LL^{-1}$, and
let us consider the extension class $\ga\in\Ext^1(\OO(2),\OO)\ot\LL=H^1(\OO(-2))\ot\LL$
of the exact sequence \eqref{g-exact-seq} in $\qgr \RR$. 
We can view the maps $\a$ and $\b$ from this exact sequence as
maps of $S$-modules 
$$\a:\QQ\to\RR_1, \ \ \b:\QQ^\vee\to \RR_1\ot\LL^{-1},$$
where $\QQ=S\oplus\LL$, satisfying $m_2(\b\ot \a)(\id_\QQ)=0$ and $m_2(\ga\ot\b)=0$.


\begin{lem}\label{qgr-enh-lem} 
Up to a gauge equivalence, there is at most one minimal $S$-linear $A_\infty$-structure on
$A(\RR,(\sideset{_{\id}}{_{\phi}}\RR)^*,1)$
with given $m_2$ and satisfying
\begin{equation}\label{m3-triangle-eq}
m_3((\ga\ot\b\ot\a)(\id_\QQ))=1.
\end{equation}
\end{lem}

\Pf . Note that
$(\sideset{_{\id}}{_{\phi}}\RR)^*\ot\LL$ is isomorphic to $\RR\ot\LL$ as a left and as a right $\RR$-module, 
so by Proposition \ref{Ext-R-RR-prop}, the conditions of
Theorem \ref{HH-gen-thm} with $d=1$ are satisfied, with $P^r$ and $P^l$ locally free over $S$ of rank $1$.
Thus, the result follows from Corollary \ref{m3-1-cor}.
\ed

\begin{rem} Using the compatibility of the higher products with exact triangles
(see e.g., \cite[Lem.\ 3.7]{Seidel-book}) one can check that
\eqref{m3-triangle-eq} holds for
the minimal $A_\infty$-structure on $A(\RR,\und{H}^1(\OO),1)$ coming from the standard 
$A_\infty$-enhancement of the derived category $D(\qgr \RR)$ (defined uniquely up to a gauge equivalence).
This explains why this is a natural condition to consider.
\end{rem}


In the case when the filtered algebra $(A,F_\bullet A)$ 
is associated with a pair of $1$-spherical objects as in Theorem \ref{spherical-filtered-thm},
we get an $A_\infty$-structure on the algebra
of the form $A(\RR,(\sideset{_{\id}}{_{\phi'}}\RR)^*\ot\LL,1)$.
Namely, we can consider the minimal $A_\infty$-structure on the subcategory of twisted complexes $(E_i)$,
obtained by homological perturbation, and use
Proposition \ref{bimodule-prop} to identify the resulting algebra with $A(\RR,(\sideset{_{\id}}{_{\phi'}}\RR)^*\ot\LL,1)$.
Note that the homological perturbation is applicable here since all the $S$-modules $\Hom(E_i,E_j)$ are finitely generated projective.
Furthermore, if we assume in addition that either $n\ge 3$ or $\tr(g)$ is a generator of $\LL$,
then we deduce that $\phi'=\phi$,
since such an automorphism is uniquely determined by its action on $\RR/(t)\simeq \EE(V,g)^{op}$
by Proposition \ref{filt-aut-prop}.

We need to check that this $A_\infty$-structure satisfies our normalization condition on $m_3$.




\begin{prop}\label{ainf-str-comparison-prop} 
Assume that $g$ is invertible, and let $\phi$ be the unique
automorphism of $\RR$, such that $\phi(t)=t$ and the induced automorphism of $\gr^F A\simeq\EE(V,g)^{op}$
is  $\Ad(g^{-1})$ (see Proposition \ref{filt-aut-prop} and Corollary \ref{H1-bimodule-cor}). Now assume that 
$\RR$ arises as the graded algebra associated with an $n$-pair $(E,F)$ as in
Theorem \ref{spherical-filtered-thm}. Then the corresponding minimal $A_\infty$-structure on 
$A(\RR,(\sideset{_{\id}}{_{\phi}}\RR)^*\ot\LL,1)$ satisfies 
$$m_3(\ga,\b,\a)=\id,$$
where $\a:E_0\to E_1\oplus E_1\ot\LL^{-1}$, $\b:E_1\oplus E_1\ot\LL^{-1}\to E_2\ot\LL^{-1}$ and
 $\ga:E_2\ot\LL^{-1}\to E_0$ correspond
to the elements \eqref{a-b-comp-eq} and to the class of the exact sequence \eqref{g-exact-seq}.
\end{prop}

\Pf . (i) It is enough to compute the relevant triple Massey product in the $A_\infty$-category of twisted complexes
over the $A_\infty$-category $\CC$ generated by our $n$-pair $(E,F)$, i.e., before applying the homological
perturbation (this follows from the functoriality of Massey products, see \cite[Prop.\ 1.1]{pol02}.

We use our presentations for $E_1$ and $E_2$ as twisted complexes from the proof of Theorem \ref{spherical-filtered-thm}
(see \eqref{E-i-simple-complex-eq}).
Thus, the complex $\hom(E_2,E)$ has form
$$\Hom^0(E,E)\to \left(\Hom^1(V^\vee L^2\ot F,E)\oplus \Hom^1(V^\vee L\ot F,E)\oplus\Hom^1(E,E)\right)$$
with the differential induced by the maps $\de_1$ and $\de_2$ (see notation in \eqref{E-i-simple-complex-eq}).

\noindent
{\bf Step 1}. We claim that the map $\ga$ is equal to the class $\ga'\in \Hom^1(E_2,E)\ot\LL_E$ of the closed element
$h^{-1}\in\End(V)\ot\LL_E\simeq \Hom^1(V^\vee L^2\ot F,E)\ot\LL_E\sub \hom^1(E_2,E)\ot\LL_E$. 
Indeed, by Lemma \ref{Ext1-h-lem},
it is enough to check that the map $\Hom^0(E,E_2)\to \Hom^1(E,E)\ot\LL_E=\LL_E^2$, given by postcomposing with $\ga'$,
coincides with the composition
$$\Hom(E,E_2)\to \Hom(E,V^\vee L^2\ot F)\simeq \End(V)\rTo{A\mapsto \tr(Ah^{-1}g)} \LL^2.$$
Indeed, this amounts to checking that the composition
$$E\to V^\vee L^2 F\rTo{h^{-1}} E\ot\LL_E[1]$$
is given by $A\mapsto \tr(Ah^{-1}g)$ which follows from the identification of the composition
$\Hom^1(F,E)\ot\Hom^0(E,F)\to \Hom^1(E,E)$ with $v^*\ot v\mapsto \lan v^*,gv\ran$.

\noindent
{\bf Step 2}. Using the computations from the proof of Theorem \ref{spherical-filtered-thm}
we see that the components of $\a$ and $\b$ are represented by the following closed maps.
The maps $E_0\rTo{t} E_1$ and $E_1\rTo{t} E_2$ are given by
\begin{diagram}
&& E\\
&&\dTo{\id_E}\\
V^\vee \ot F&\rTo{\de_1}&E
\end{diagram} 
\begin{diagram}
&& V^\vee \ot F&\rTo{\de_1}&E\\
&&\dTo{\id_E}&&\dTo{\id_E}\\
V^\vee \ot F&\rTo{\de_2}&V^\vee \ot F&\rTo{\de_1}&E
\end{diagram} 
Also, for a certain choice of $\wt{h}$, the maps 
$E_0\rTo{\wt{h}} E_1\ot\LL^{-1}_E$ and $E_1\rTo{\wt{h}} E_2\ot\LL_E^{-1}$ are given by
\begin{diagram}
E\\
\dTo{h}\\
V^\vee \ot F\ot\LL_E^{-1}&\rTo{\de_1}&E\ot\LL_E^{-1}
\end{diagram}
\begin{diagram}
V^\vee \ot F&\rTo{}&E\\
\dTo{h^*\ot\id_F}&\rdTo{}&\dTo{h}\\
V^\vee \ot F\ot\LL_E^{-1}&\rTo{\de_2}&V^\vee \ot F\ot\LL_E^{-1}&\rTo{\de_1}&E\ot\LL_E^{-1}
\end{diagram}
where the diagonal arrow is $\mu_h\ot\id_F$ (recall that $\mu_a$ is defined by \eqref{mu-a-eq}).

Note that any other choice of $\wt{h}$ is of the form $\wt{h}+c\cdot t$, for $c\in \LL^{-1}$. Hence, for a different choice
of $\wt{h}$, the sequence of maps $\a,\b$ would change by an automorphism of $E_1\oplus E_1\ot\LL^{-1}$, so the Massey
product is unaffected by such a change.

\noindent
{\bf Step 3}. Now one easily checks that $m_2(\b,\a)=0$ and that $m_3(\ga,\b,\a)=0$. However, the product
$m_2(\ga,\b)$ is not zero on the cochain level. In fact, its only nonzero component is given by the composition 
\begin{diagram}
V^\vee \ot F&\rTo{\de_1}&E\\
\dTo{h^*\ot\id_F}&\rdTo{}&\dTo{h}\\
V^\vee \ot F\ot\LL_E^{-1}&\rTo{\de_2}&V^\vee \ot F\ot\LL_E^{-1}&\rTo{\de_1}&E\ot\LL_E^{-1}\\
&\rdTo{h^{-1}}\\
&&E
\end{diagram}
It is easy to check that the composition $h^{-1}\circ (h^*\ot\id_F)$ in this diagram is equal to $\de_1$.
Hence, we obtain
$$m_2(\ga,\b)=d(\id_E),$$
where we view $\id_E$ as an element of $\hom^0(E_1,E_0)$.
Thus, by the definition of the Massey product, we get 
$$MP(\ga,\b,\a)=m_2(\id_E,t)=\id_E.$$
\ed

\section{Proofs of the main results}\label{proof-sec}

\subsection{Proof of Theorem A for Noetherian rings}\label{thmA-Noeth-sec}

For a commutative ring $S$ we can think of an $S$-point of $\PGL_n$ as isomorphism classes of
pairs $(g,\LL)$, where $\LL$ is an invertible $S$-module and
$g\in \End_S(V)\ot\LL$ (where $V=S^n$) is an invertible element.
Recall that for a graded $S$-algebra $B$ and a graded $B-B$-bimodule $M$, we define
$A(B,M,1)$ as the trivial square-zero extension $B\oplus M$ (see Def.\ \ref{ABM-def}).

Let us consider the functors on the category of Noetherian commutative rings, that associate to
$S$ the set of 

\noindent
(1) $(g,\LL)\in\PGL_n(S)$ and minimal $A_\infty$-structures on $\SS(V,g)$ up to a gauge equivalence;

\noindent
(2) $(g,\LL)\in\PGL_n(S)$ and isomorphism classes of $(A,F_\bullet A, \iota,\phi;m_\bullet)$, where $(A,F_\bullet A)$
a filtered algebra equipped with an isomorphism $\iota:\gr^F A\simeq \EE(V,g)^{op}$
and an automorphism $\phi:A\to A$ 
such that the induced automorphism $\ov{\phi}$ of $\gr^F A\simeq \EE(V,g)^{op}$ is equal to $\Ad(g^{-1})$;
and $m_\bullet$ is a minimal $A_\infty$-structure on $A(\RR(A),(\sideset{_{\id}}{_{\phi}}\RR(A))^*\ot\LL,1)$ with given $m_2$
and such that $m_3(\ga,\b,\a)=1$ (viewed up to a gauge equivalence);

\noindent
(3) $(g,\LL)\in\PGL_n(S)$ and isomorphism classes of $(A,F_\bullet A, \iota)$ as in (2).

In the case $n=2$ we always assume in addition that $\tr(g)\neq 0$.

\noindent
{\bf Step 1}. Construction of an injective map from (1) to (2).

Starting from a minimal $A_\infty$-structure on $\SS^n(V,g)$, 
we consider the corresponding $n$-pair of spherical objects $(E,F)$
(see Sec.\ \ref{n-pairs-1-sph-sec}). Now we
consider twisted objects $(E_i)$ in Theorem \ref{spherical-filtered-thm} and use this Theorem and
Proposition \ref{bimodule-prop} to identify the corresponding algebra 
$$\bigoplus_i \Hom(E_0,E_i)\oplus \bigoplus_i \Ext^1(E_i,E_0)$$ 
with $A(\RR(A),(\sideset{_{\id}}{_{\phi}}\RR(A))^*\ot\LL,1)$, for some filtered algebra $(A,F_\bullet A)$ and an automorphism
$\phi$ satisfying the conditions in (2). Applying homological perturbation (see \cite[Sec.\ 3.3]{keller-ainf}) to the $A_\infty$-structure on the subcategory $(E_i)$
we get a minimal $A_\infty$-structure $m_\bullet$ on $A(\RR(A),(\sideset{_{\id}}{_{\phi}}\RR(A))^*\ot\LL,1)$.
Furthermore, the condition on $m_3$ holds due to Proposition \ref{ainf-str-comparison-prop}.

Next, let us show injectivity of this map. For a minimal $A_\infty$-structure $m_\bullet$ on $\SS^n(V,g)$ let us denote
by $\SS(m_\bullet)$ the corresponding $A_\infty$-algebra, which can be also viewed as an $A_\infty$-category with
two objects $(E,F)$. Let $\Pi Tw \SS(m_\bullet)$ denote the $A_\infty$ split-closure of the category of twisted complexes
over $\SS(m_\bullet)$. The exact triangle
\begin{equation}\label{E0-E1-F-ex-tr}
E_0\to E_1\to V^\vee\ot F\to E_0[1]
\end{equation}
coming from the definition of $E_1=T_F(E_0)$, shows that $\Pi Tw \SS(m_\bullet)$ is split-generated by $E_0$ and $E_1$.
Thus, by \cite[Cor. 4.9]{Seidel-book}, the inclusion of the full subcategory
on objects $(E_i)$
$$\{ E_i \ |\ i\ge 0\} \hra \Pi Tw \SS(m_\bullet),$$
extends to a quasi-equivalence
$$\Pi Tw \{ E_i \ |\ i\ge 0\} \rTo{\sim} \Pi Tw \SS(m_\bullet).$$
Thus, if two minimal $A_\infty$-structures on $\AA_n$, $(m_\bullet)$ and $(m'_\bullet)$ induce
gauge-equivalent $A_\infty$-structures on $\lan E_i \ |\ i\ge 0\ran$, then there exists a quasi-equivalence
$$\Phi:\Pi Tw \SS(m_\bullet)\simeq \Pi Tw \SS(m'_\bullet)$$
such that $H^0\Phi$ is the identity on $\Hom^*(E_i,E_j)$. 
Since the functor $H^0\Phi$ is triangulated, the exact triangle \eqref{E0-E1-F-ex-tr} shows that
$$V^\vee\ot \Phi(F)\simeq \Phi(V^\vee\ot F)\simeq V^\vee\ot F.$$
Such an isomorphism is induced by a unique isomorphism
$\Phi(F)\simeq F\ot\MM$, for some locally free $S$-module of rank $1$ equipped with an isomorphism
$V^\vee\ot \MM\simeq V^\vee$. 

Localizing over an open affine covering $\Spec(S)$, we can trivialize $\LL_E$ and $\MM$.
Then $\Phi$ gives an $A_\infty$-autoequivalence of the subcategory $\{E,F\}$, identical on objects, and 
such that the induced autoequivalence
of $\{E_i \ |\ i\ge 0\}$ is isomorphic to the identity.
To prove that $(m_\bullet)$ and $(m'_\bullet)$ are
gauge equivalent, it is enough to have that $H^0\Phi|_{\{E,F\}}$ is isomorphic to the identity.
Just using the condition the the endomorphism of $\Hom^1(E,E)$ induced by $H^0\Phi$ is the identity,
it is easy to see that $H^0\Phi$ should have the following form: it is given by the maps 
\begin{align*}
&h:\Hom^0(E,F)=V\to V=\Hom^0(E,F), \ \ (h^{-1})^*:V^\vee=\Hom^1(F,E)\to \Hom^1(F, F)=V^\vee, \\
&\Hom^1(F,F)\rTo{\la\cdot ?}\Hom^1(F,F),
\end{align*}
where $h\in\GL(V)$ and $\la\in S^*$ satisfy
$$h^{-1}gh=\la g.$$
Now the condition that the induced automorphism of the graded algebra $\RR=\bigoplus_{i\ge 0}\Hom(E_0,E_i)$
is the identity, implies that the automorphism $\Ad(h^{-1})$ of $\RR/(t)\simeq\EE(V,g)^{op}$ is the identity.
But this is possible only when $h$ is a scalar matrix, $h=c\cdot \id$. In this case we can change $\Phi$
by an isomorphic $A_\infty$-equivalence (rescaling at $F$), so that $H^0\Phi$ becomes the identity.

Thus, we deduce that $(m_\bullet)$ and $(m'_\bullet)$ are gauge equivalent locally over $\Spec(S)$.
Applying \cite[Thm.\ 2.2.6(ii)]{P-more-pts}, we deduce that they are globally gauge equivalent.

\noindent
{\bf Step 2}. The forgetful map from (2) to (3) is injective. 
Indeed, the uniqueness of the automorphism $\phi$ follows from Prop.\ \ref{filt-aut-prop} (here we use the assumption that $\tr(g)\neq 0$ if $n=2$), 
while the uniqueness of the $A_\infty$-structure
follows from Lemma \ref{qgr-enh-lem}.

\noindent
{\bf Step 3}. We have a natural map from (3) to (1): starting from $(A,F_\bullet A,\iota)$ we construct an $n$-pair $(E,F)$ by
considering the derived category $D^b(\qgr \RR(A))$ and considering the objects $E=\OO$ and $F\in \qgr \RR(A)$ 
defined in Proposition \ref{proj-coh-prop}(ii). 

\noindent
{\bf Step 4}. We claim that the composition $(3)\to (1) \to (2) \to (3)$ is the identity. Together with the injectivity proved in
Steps 1 and 2, this would imply that our arrows give bijections between data (1), (2) and (3).

Thus, we start from $(A,F_\bullet A,\iota)$, consider the corresponding $n$-pair of spherical objects $(E,F)$ as in Step 3,
then look at the twisted complexes $E_i=T_F^i(E)$ and consider the corresponding $\Hom$-algebra $\RR_{T_F,E}$.
By \cite[Lem.\ 3.34]{Seidel-book}, the inclusion of the full subcategory $\{E,F\}$ extends to an $A_\infty$-functor
$$\Phi: Tw\{E,F\}\to D(\qgr \RR(A)),$$
which is a quasi-equivalence with its image.
Thus, Proposition \ref{proj-coh-prop}(iv) gives the required isomorphism of the algebra $\RR_{T_F,E}$ with $\RR(A)$,
preserving the natural central elements $t$.

Note that the last assertion of the theorem follows from the fact that the composition
$(1)\to (2)\to (3)\to (1)$ is the identity.
\ed

\begin{rem} Note that going from (3) to (1) and then to (2) equips any filtered algebra $A$ as above with a canonical filtered automorphism $\phi$ as in (2).
This is precisely the automorphism constructed in Cor.\ \ref{H1-bimodule-cor}.
\end{rem}

\subsection{Moduli spaces and the proof of Theorem A}\label{moduli-thmA-sec}

We would like to show that
the functor associating to $S$ an isomorphism class of the data
\begin{equation}\label{filtered-data}
(\LL,g,A,F_\bullet A,\iota:\gr^F(A)\simeq\EE(V,g)^{op})
\end{equation}
as before (where $\LL$ is a locally free $S$-module of rank $1$, $g\in \End(S^n)\ot\LL$ is an invertible
element, such that $\tr(g)$ is a generator of $\LL$ if $n=2$),
is representable by an affine scheme $\Spec S_{filt}$ of finite type over $\Z$.

For $n\ge 3$, let $S_0$ be the algebra of functions on $\PGL_n$ (which is the degree $0$ part in the localization
$\Z[x_{ij}][\det^{-1}]$), and let $g\in\End(S_0^n)\ot\OO_{S_0}(1)$ be the universal invertible element.
In the case $n=2$, we define $S_0$ as the algebra of functions on the open subset of $\PGL_n$
given by nonvanishing of $\tr(g)^n/\det(g)$. 

Recall that the algebra $\EE:=\EE(S_0^n,g)^{op}$ is Koszul (see Lemma \ref{E-deg-1-lem}(iii)).
This means that the filtered algebras $A$ we would like to study are given by nonhomogeneous quadratic relations
whose homogeneous quadratic parts are the quadratic relations in $\EE$.
Let $I_\EE\sub \EE_1\ot_{S_0} \EE_1$ denote the space of quadratic relations. 

\begin{lem}\label{split-emb-der-lem} 
The natural morphism 
\begin{equation}\label{nabla-eq}
\nabla:\EE_1^\vee\to \Hom_{S_0}(I_\EE,\EE_1):\xi\mapsto (e\ot e'\mapsto \xi(e)e'+\xi(e')e)
\end{equation}
is a split embedding of $S_0$-modules.
\end{lem}

\Pf . Since both $\EE_1^\vee$ and $\Hom_{S_0}(I_\EE,\EE_1)$ are finitely generated projective modules over $S_0$,
it is enough to check that for any $S_0$-algebra $S$, the morphism $\nabla_S$, obtained from $\nabla$ by extension of
scalars, is injective. Now we note that since $\EE_S=\EE\ot_{S_0}S$ is quadratic, the elements of $\ker(\nabla_S)$ are precisely
derivations of $\EE_S$ of degree $-1$. By Proposition \ref{derivations-E-prop}, any such derivation is zero.
\ed

\begin{prop}\label{filt-repr-prop} 
The functor associating to $S$ the set of isomorphism classes of triples \eqref{filtered-data} is
representable by an affine scheme $\Spec S_{filt}$ of finite type over $S_0$ (and hence, over $\Z$).
\end{prop}

\Pf . As was discussed above, our functor associates to $S$ the set of nonhomogeneous quadratic algebras deforming 
$\EE_S=\EE\ot_{S_0}S$.
For such an algebra we can always choose a splitting $s:\EE_{S,1}\to F_1 A$ of the projection 
$F_1 A\to F_1 A/F_0 A\simeq \EE_{S,1}$.
Set $U:=s(\EE_{S,1})\sub F_1 A$.
Then algebra $A$ can be given by some submodule of nonhomogeneous quadratic relations 
$$I_A\sub U\ot_S U\oplus U\oplus S,$$
such that the projection to $U\ot_S U\simeq \EE_{S,1}\ot_S \EE_{S,1}$ induces an isomorphism of $I_A$ with the submodule
$I_{\EE,S}$ of quadratic relations in $\EE_S$. Thus, $I_A$ is the graph of an $S$-linear map
$$(\phi,\th):I_{\EE,S}\to U\oplus S.$$
Conversely, starting from such data we can construct the algebra $T(U)/(I_A)$ which is equipped with a map
$\EE\to\gr^F(A)$. Since the algebra $\EE$ is Koszul, by \cite[Sec.\ V.2]{PP-book}, this correspondence gives a bijection between
the set of quadruples $(A,F_\bullet A,\iota,s)$, where $s:\EE_{S,1}\to F_1 A$ is a splitting,
and pairs of maps $(\phi,\th)$ satisfying certain quadratic equations (analogs of Jacobi identity).
\footnote{In \cite{PP-book} we work over a field, however, the argument still applies in the case when all the graded 
components $\EE_{S,i}$ are finitely generated projective modules over $S$, and $\EE_{S,0}=S$.}

Now we can change a splitting $s$ to $s+\xi$, where $\xi\in\Hom_S(\EE_{S,1}(S),S)$.
It is easy to see that this corresponds to a certain action of $\Hom_S(\EE_{S,1}(S),S)$ (viewed as an additive group)
on pairs $(\phi,\th)$. Furthermore, $\phi$ gets changed to $\phi+\nabla(\xi)$, with $\nabla$ given by \eqref{nabla-eq}.
Thus, if we fix a complementary $S_0$-submodule $K\sub \Hom_{S_0}(I_\EE,\EE_1)$ to the image of $\nabla$,
which is possible by Lemma \ref{split-emb-der-lem}, then every orbit of the above action has a unique representative
$(\phi,\th)$ with $\th\in K_S$. The set of such $(\phi,\th)$, satisfying the quadratic equations mentioned above,
is the required affine scheme of finite type over $S_0$.
\ed

\noindent
{\it End of proof of Theorem A}. 
As we have seen in Proposition \ref{filt-repr-prop} and Corollary \ref{ainf-moduli-cor}, 
we have two finitely generated $\Z$-algebras,
$S_{filt}$ and $S_{A_\infty}$, that represent the functors (3) and (1) from Sec.\ \ref{thmA-Noeth-sec}.
Since both $S_{A_\infty}$ and $S_{filt}$ are Noetherian, by the Noetherian case of Theorem A(i)
and by Yoneda lemma, we obtain an isomorphism 
$$S_{filt}\simeq S_{A_\infty}.$$
This gives the required isomorphism of functors.
\ed

\subsection{Polarizing line bundles on cyclotomic stacks}\label{polarizing-sec}

We need a technical result characterizing stacks obtained as $\Proj^{st}$ of a graded algebra in terms of polarizing line bundles
(see \eqref{Proj-stack-eq} for the definition of the stacky version of the $\Proj$-construction).
This result is based on \cite[Sec.\ 2.4.2]{AH} (in fact, all the needed arguments are already in \cite{AH}, just
not the statement itself).

Let $\XX$ be a proper algebraic stack over a field $k$.
Following \cite{AH}, we say that $\XX$ is {\it cyclotomic} if it has cyclotomic stabilizers (i.e., each geometric fiber of
$\II_{\XX}\to \XX$ is isomorphic to $\mu_d$ for some $d$). 
A line bundle $\LL$ over $\XX$ is called {\it uniformizing} if for each geometric point $p$ in $\XX$, the action of
$\Aut(p)$ on the fiber of $\LL$ is effective. Let us consider the principal $\G_m$-bundle over $\XX$ associated with $\LL$:
$$\PP_{\LL}:=\Spec_{\XX}(\bigoplus_{i\in \Z}\LL^i)\to \XX.$$
By \cite[Prop.\ 2.3.10]{AH}, $\LL$ is uniformizing if and only if the stack $\PP_{\LL}$ is representable.

Now let $X$ be the coarse moduli space of $\XX$, and let $\pi:\XX\to X$ be the projection.
A {\it polarizing} line bundle $\LL$ over $\XX$ is a uniformizing line bundle such that some positive power of $\LL$
is isomorphic to the pull-back of an ample line bundle on $X$.

\begin{prop}\label{polarizing-prop} 
If $\LL$ is a polarizing line bundle on a proper cyclotomic stack $\XX$ then there is a natural isomorphism
$$(\XX,\LL)\simeq (\Proj^{st}(R_\LL),\OO(1)),$$
where $R_\LL$ is the graded algebra given by
$$R_\LL:=\bigoplus_{j\ge 0}H^0(\XX,\LL^j).$$
\end{prop}

\Pf . Let $M$ be an ample line bundle on $X$ such that $\LL^N\simeq \pi^*M$.
Interpreting elements of $R_\LL$ as functions on $\PP_\LL$ we get a regular morphism
\begin{equation}\label{principal bundle-Spec-map}
\PP_\LL\to \Spec(R_\LL)\setminus \{0\}
\end{equation}
compatible with the $\G_m$-action.
On the other hand, the homomorphism of rings $\bigoplus M^i\to \bigoplus_j \pi_*\LL^j$ gives rise to
a finite morphism $\PP_\LL\to \PP_M$. Let us consider the commutative diagram
\begin{diagram}
\PP_\LL &\rTo{} &\Spec(R_\LL)\\
\dTo{}&&\dTo{}\\
\PP_M&\rTo{}& \Spec(R_M)
\end{diagram}
in which the lower horizontal arrow is an open embedding with a dense image (since $M$ is ample).
By Zariski's main theorem, the composed map $\PP_\LL\to \Spec(R_M)$ can be
factored as the composition 
$$\PP_\LL\rTo{j} \hat{\PP}_\LL \rTo{f} \Spec(R_M)$$
where $j$ is an open embedding with a dense image and $f$ is finite. As shown in the proof of \cite[Cor.\ 2.4.4]{AH}, in fact, $\hat{\PP}_\LL=\Spec(R_\LL)$.
Thus, the morphism \eqref{principal bundle-Spec-map} is a dense open embedding, and so is the induced map of quotients by $\G_m$,
$$\XX\to \Proj^{st}(R_\LL).$$
Since $\XX$ is proper, this map is an isomorphism.
\ed


\subsection{Proof of Theorem B}\label{thmB-sec}

We have discussed in Sec.\ \ref{filtered-alg-sph-orders-sec} the construction of the order on a neat stacky pointed curve associated with a filtered algebra
$(A,F_\bullet)$ (see especially Lemma \ref{stacky-sheaf-lem}). 

To go from a neat stacky curve $(C,p)$ with an order $\AA$ to a filtered algebra,
let us consider the filtration $F_iA=H^0(C,\AA(ip))$ on the algebra $A=H^0(C\setminus p,\AA)$.
From the trivialization of $\OO(p)|_p$ (given by the tangent vector)
we get a natural injective homomorphism of graded algebras
$$\gr^F A\to \End(V)^{op}\ot k[z].$$
We claim that its image is $\EE(V,g)^{op}$ for some $g\in \P\End(V)$ (see \eqref{E-V-eq}).
Indeed, the exact sequence
$$0\to H^0(C,\AA)\to H^0(C,\AA(p))\to H^0(p,\AA(p)|_p)\to H^1(C,\AA)\to 0$$
shows that the image of the restriction map $F_1A\to H^0(p,\AA(p)|_p)\simeq\End(V)$
is of codimension $1$. Hence, it has form $\End_g(V)$ for a unique $g\in\P\End(V)$.
On the other hand, since $H^1(C,\AA(ip))=0$ for $i\ge 1$, the restriction maps
$$F_mA=H^0(C,\AA(mp))\to H^0(p,\AA(mp)|_p)\simeq \End(V)$$
are surjective for $m\ge 2$, which proves our claim. Thus, we get an isomorphism of $\gr^F A$ with $\EE(V,g)^{op}$.

Using Lemma \ref{stacky-sheaf-lem}(iv) it is easy to check that starting from a filtered algebra $(A,F_\bullet)$ and
constructing an order $\AA$ over a stacky curve $C$, we then recover the original filtered algebra by the above construction.

Conversely, if we start with an order $\AA$ over a neat pointed stacky curve $(C,p)$ and consider the filtered algebra
$(A,F_\bullet)$ with $F_iA=H^0(C,\AA(ip))$, then we recover $(C,p,\AA)$ by the $\Proj^{st}$ construction,
described in Sec.\ \ref{filtered-alg-sph-orders-sec}.
More precisely, since $\mu_d$ acts faithfully on the fiber of $\OO_C(p)$ at $p$, this line bundle is uniformizing.
Furthermore, $\OO_C(dp)$ is a pull-back of the ample line bundle $\OO_{\ov{C}}(\ov{p})$ on the coarse moduli $\ov{C}$
(where $\ov{p}\in \ov{C}$ is the image of $p$).
Hence, $\OO_C(p)$ is polarizing, and by Proposition \ref{polarizing-prop}, we have an isomorphism
$$(C,\OO_C(p))\simeq (\Proj^{st}(\RR(Z)),\OO(1)),$$
where $Z=\cup_{i\ge 0}H^0(C,\OO_C(ip))$ is the center of $A$ (recall that $\RR(\cdot)$ denotes the
Rees construction). Furthermore, any coherent
sheaf $\FF$ on $\Proj^{st}(\RR(Z))$ is identified with the localization of the corresponding
graded $\RR(Z)$-module, $\bigoplus_i H^0(\FF(i))$. Applying this to $\AA$, we see that it corresponds
to the localization of $\RR(A)$.


Next, by Lemma \ref{Perf-gen-lem}, the pair $(\AA,\rho_*V)$ generates $\Perf(\AA^{op})$.
Thus, by Proposition \ref{Proj-order-prop}, we have an equivalence of $\Perf(\AA^{op})$ with a full subcategory in
$D\qgr \RR(A)$, which sends $(\AA,\rho_*V)$ to the pair $(\OO,F)$.

If $g$ is invertible then, by Proposition \ref{proj-coh-prop}(iii), the pair $(\OO,F)$ in $D\qgr \RR(A)$ 
is an $n$-pair of $1$-spherical objects,
so $\AA$ is a spherical order. Conversely, if $\AA$ is spherical then $g$ is invertible by Proposition
\ref{spherical-order-prop}(ii).

Finally, let us prove that $\AA$ is symmetric if and only if $g$ is scalar.
Note that for any spherical order we have a canonical {\it Nakayama automorphism} $\kappa$ defined by the equation
$$\tau(yx)=\tau(x\kappa(y)),$$
where $\tau:\AA\to\om_C$ is a nonzero morphism. Indeed, we have
two isomorphisms of coherent sheaves,
$$\nu:\AA\to\und{\Hom}(\AA,\om_C): a\mapsto (x\mapsto \tau(xa)), \ \ 
\nu':\AA\to\und{\Hom}(\AA,\om_C): a'\mapsto (x\mapsto \tau(a'x))
$$
(see Proposition \ref{spherical-order-prop}),
and we set $\kappa=\nu^{-1}\circ\nu'$. The fact that $\kappa$ is an automorphism of algebras follows from the defining identity.
 
By Proposition \ref{spherical-order-prop}(ii), the restricted functional 
$$\tau|_p:\AA|_p\simeq \End(V)\ot R_{\mu_d}\to \chi$$ 
is given by $x\mapsto \tr(gx)$ on $\End(V)\ot\chi$. Hence, we have 
$$\kappa|_p=\Ad(g)\ot\id:\End(V)\ot R_{\mu_d}\to \End(V)\ot R_{\mu_d},$$
where $R_{\mu_d}$ is the regular representation of $\mu_d$.
This immediately shows that if $\AA$ is symmetric, i.e., $\kappa=\id$, then $g$ is scalar.
Conversely, assume that $g$ is scalar. Then $\kappa|_p=\id$. Now $\kappa$ induces a filtered automorphism
of the algebra $A=H^0(C\setminus\{p\},\AA)$, and the condition that $\kappa|_p=\id$ implies that the
induced automorphism of $\gr^F A$ is equal to the identity. Hence, by Proposition \ref{filt-aut-prop}, 
this automorphism of $A$ is equal to the identity, and so $\kappa=\id$. 
\ed

\subsection{A criterion for cyclic $A_\infty$-structures}\label{cyclic-ainf-sec}

Recall that an $A_\infty$-algebra over a field $k$ is called cyclic if it is equipped with a bilinear form $\lan\cdot,\cdot\ran$
such that
$$\lan m_n(a_1,\ldots,a_n),a_{n+1}\ran=(-1)^{n(\deg(a_1)+1)}\lan a_1,m_n(a_2,\ldots,a_{n+1})\ran.$$

Kontsevich and Soibelman give a general criterion \cite[Thm.\ 10.2.2]{KS-notes} in the case when $\cha(k)=0$
stating that such a cyclic structure exists on a minimal model of an $A_\infty$-algebra $A$ with finite dimensional cohomology
$H^*(A)$, provided there is a functional $\th:HC_N(A)\to k$ (where $HC_*(A)$ is the cyclic homology of $A$) such that the
induced pairing on $H^*(A)$,
$$\lan x,y\ran=\th(\iota(xy)),$$
where $\iota:H^*(C)\to HC_*(A)$ is the natural map, is perfect.

In the case of algebras of the form $H^*(C,\AA)$, where $\AA$ is a sheaf of coherent algebras over a curve $C$,
we can provide a more direct construction of a cyclic structure, which only uses the assumption that
$\cha(k)\neq 2$, and relies instead on a cyclic version of the homological perturbation considered in \cite{Lazaroiu}.

\begin{prop}\label{cyclic-deg-[01]-prop} 
Let $B=B_0\oplus B_1$ be a dg-algebra over a field $k$, concentrated in degrees $[0,1]$, $\lan\cdot,\cdot\ran$
a pairing of degree $1$ on $B$ satisfying 
$$\lan x,y\ran=(-1)^{\deg(x)\deg(y)}\lan y,x\ran,$$
$$\lan dx,y\ran+(-1)^{\deg(x)}\lan x,dy\ran=0.$$
Assume also that $H^*(B)$ is finite-dimensional and the induced pairing on $H^*(B)$ is perfect. 
Then the data for the homological perturbation can be chosen in such a way that the resulting minimal
$A_\infty$-structure on $H^*(B)$ is cyclic with respect to the pairing induced by $\lan\cdot,\cdot\ran$.
\end{prop}

\Pf . Let $A\sub\ker(d)\sub B$ be any (graded) subspace of cohomology representatives.  We claim that
there exists a subspace $C\sub B_0$, complementary to $\ker(d)$, such that $\lan C,A_1\ran=0$.
Indeed, let us start with an arbitrary such complement $C\sub B_0$. Then the pairing $C\ot A_1\to k$ can be interpreted
as a map $C\to A_1^*\simeq A_0$ (where the latter isomorphism is given by the pairing between $A_0$ and $A_1$).
Correcting $C$ by this map, we get a new subspace in $C\oplus A_0$, which is still complementary to $\ker(d)$,
and which is orthogonal to $A_1$.

Note that we have orthogonalities $\lan C,A\ran=0$, $\lan C,C\ran=0$. Hence, the standard homotopy operator $Q:B\to B$
associated with the decomposition $B=\im(d)\oplus A\oplus C$ satsifies 
$$\lan Qx, y\ran=(-1)^{\deg(x)}\lan x,Qy\ran.$$
As was observed in \cite[Sec.\ 3.3]{Lazaroiu}, this implies that
the minimal $A_\infty$-structure on $H^*(B)$ given by the tree formula from \cite{KS} is cyclic.
\ed


We apply this general result in the following geometric setup.


\begin{prop}\label{sheaf-dg-cyclic-prop}
Let $C$ be a tame proper DM-stacky curve over a field $k$ of characteristic $\neq 2$, with a Cohen-Macaulay
coarse moduli space $\ov{C}$ such that $H^0(\ov{C},\OO)=k$.
Let $\AA$ be a coherent sheaf of $\OO_C$-algebras, equipped with a morphism $\tau:\AA\to \om_C$.
Assume that we have a morphism $\tau:\AA\to\om_C$ such that $\tau(xy)=\tau(yx)$. 
Assume that the induced pairing 
$$\AA\ot \AA\to \om_C$$
induced by $\tau(xy)$, is perfect in the derived category (either on the left or on the right). 
Then the minimal $A_\infty$-structure on $H^*(C,\AA)$ obtained by the homological perturbation can
be chosen to be cyclic with respect to the pairing $\th(xy)$, where 
$$\th:H^1(C,\AA)\to H^1(C,\om_C)\to k$$
is the functional induced by $\tau$.
\end{prop} 

\Pf . First of all, since $C$ is tame, we have an isomorphism of algebras 
$$H^*(C,\AA)\simeq H^*(\ov{C},\ov{\AA}),$$ 
where
$\ov{\AA}:=\pi_*\AA$ and $\pi:C\to \ov{C}$ to the coarse moduli map. Also, we have an isomorphism
$$\pi_*\om_C\simeq \om_{\ov{C}}$$
(see \cite[Prop.\ 2.3.1]{Nironi}). Thus, we can view the morphism $\pi_*\tau$ as a morphism 
$$\ov{\tau}:\ov{\AA}\to\om_{\ov{C}}.$$
Furthermore, the induced pairing 
$$\ov{\AA}\ot \ov{\AA}\to \om_{\ov{C}},$$
given by $\ov{\tau}(xy)$,
factors through the natural projection $\ov{\AA}\ot\ov{\AA}\to \pi_*(\AA\ot \AA)$ and hence, is still symmetric.
Finally, by the relative duality we have an isomorphism
$$\pi_*R\Hom(\AA,\om_C)\simeq \pi_*R\Hom(\AA,\pi^!\om_{\ov{C}})\simeq R\Hom(\ov{\AA},\om_{\ov{C}}),$$
which implies the $\ov{\tau}$ still satisfies the required non-degeneracy condition.
Thus, replacing $(C,\AA,\tau)$ with $(\ov{C},\ov{\AA},\ov{\tau})$, we can assume that $C$ is a usual (non-stacky) curve.

We can compute $H^*(C,\AA)$ using the Cech resolution 
$$C(\AA):\AA(U_1)\oplus \AA(U_2)\rTo{\de}\AA(U_{12}),$$
with respect to a covering $C=U_1\cap U_2$, where
$U_i$ are open affine subsets, $U_{12}=U_1\cap U_2$. 
Here $\de(f_1,f_2)=f_2-f_1$. Since $\cha(k)\neq 2$, we can equip
$C(\AA)$ with the following dg-algebra structure: the product on $\AA(U_1)\oplus \AA(U_2)$ is the one on direct
sum of algebras, while for $(f_1,f_2)\in\AA(U_1)\oplus \AA(U_2)$, $g\in \AA(U_1\cap U_2)$, we set
$$(f_1,f_2)g=\frac{(f_1+f_2)|_{U_{12}}g}{2}, \ \ g(f_1,f_2)=\frac{g(f_1+f_2)|_{U_{12}}}{2}.$$
Note that with respect to this product we have 
$$[(f_1,f_2),g]=\frac{1}{2}([f_1|_{U_{12}},g]+[f_2|_{U_{12}},g]).$$
Since the map $\tau$ vanishes on the commutators, the induced map of Cech complexes
$$C(\AA)\to C(\om_C)$$
also does, with respect to the above product. Composing this map with a map $C(\om_C)\to k[-1]$, realizing
the canonical trace morphism $H^1(C,\om_C)\to k$, we get a map $\th:C(\AA)_1\to k$, vanishing on the image of the
differential and on the commutators. Furthermore, Serre duality implies that
the map 
$$H^i(C,\AA)\ot H^{1-i}(C,\AA)\to H^1(\om_C)\to k,$$
induced by $\tau(xy)$ is a perfect pairing. Thus, using Proposition \ref{cyclic-deg-[01]-prop}  
we get a cyclic $A_\infty$-structure with 
respect to $\th(xy)$.
\ed

\begin{rem} In characteristic zero there is a generalization of Proposition \ref{sheaf-dg-cyclic-prop} 
to coherent sheaves of dg-algebras over 
higher-dimensional schemes over a field of characteristic zero. In this case one has to use Thom-Sullivan construction 
(see \cite[Sec.\ 5.2]{HS}, \cite[App.\ A,B]{VdB-def-qu}) to get a multiplicative structure on derived global sections,
and then apply the criterion of Kontsevich-Soibelman \cite[Thm.\ 10.2.2]{KS-notes}.
\end{rem}

One more observation is that the assumptions of Proposition
\ref{sheaf-dg-cyclic-prop} are preserved when passing from $\AA$ to the endomorphism sheaf of a locally projective
$\AA$-module.

\begin{lem}\label{trace-tau-lem} 
Let $\AA$ be a coherent sheaf of $\OO_C$-algebras, together with a morphism $\tau:\AA\to\om_C$, 
satisfying the assumptions of
of Proposition \ref{sheaf-dg-cyclic-prop}, and let $\PP$ be a locally projective finitely generated $\AA$-module.
Consider the sheaf of algebras $\wt{\AA}:=\und{\End}_\AA(\PP)$.
Then the assumptions of Proposition \ref{sheaf-dg-cyclic-prop} still hold for $\wt{\AA}$ and
$\wt{\tau}:\wt{\AA}\to\om_C$ defined as the composition of $\tau$ and the trace morphism
$\tr:\und{\End}_\AA(\PP)\to\AA/[\AA,\AA]$.
\end{lem}

\Pf . Note that as $\OO$-module, $\wt{\AA}$ is locally a summand in $\AA^{\oplus n}$ for some $n$.
Hence, the assumption that $\und{\Ext}^{>0}(\AA,\om_C)=0$ implies that the same holds for $\wt{\AA}$,
so we only need to check that $\wt{\tau}$ vanishes on $[\wt{\AA},\wt{\AA}]$ and that the pairing
$\wt{\tau}(xy)$ is perfect. The former is the standard fact about traces. For the latter we
can assume $\PP$ to be a direct summand of $\AA^{\oplus n}$. First, we observe
that the pairing 
$$\tau(\tr(xy)):\Mat_n(\AA)\ot_\OO \Mat_n(\AA)\to \OO$$
is a direct sum of pairings $\AA\cdot e_{ij} \ot \AA\cdot e_{ji} \to \OO$, which are perfect by assumption.
Next, consider a direct sum decomposition $\AA^{\oplus n}=\PP\oplus \QQ$, and let $e_\PP$ and $e_\QQ$ be
the corresponding idempotents in $\Mat_n(\AA)$. Then to deduce that the restriction of $\tau(\tr(xy))$
to $e_\PP\Mat_n(\AA)e_\PP$ it is enough to check that the decomposition
$$\Mat_n(\AA)=e_\PP\Mat_n(\AA)e_\PP\oplus \bigl(e_\QQ\Mat_n(\AA)\oplus e_\PP\Mat_n(\AA)e_\QQ\bigr)$$ 
is orthogonal with respect to our form. But this immediately follows from the identities $e_\PP e_\QQ=0$
and $\tr(yx)\equiv \tr(xy)\mod [\AA,\AA]$. 
\ed

\medskip

\noindent
{\it Proof of Corollary C}.
The first part follows immediately from Theorem B: we can realize every $A_\infty$-structure on $\SS(k^n,\id)$
by the one coming from a symmetric spherical order $\AA$. 

For the last assertion, we use the fact that for such $\AA$, a nonzero morphism $\tau:\AA\to\om_C$ induces a symmetric
pairing $\AA\ot\AA\to\om_C$ which is perfect in derived category (see Proposition \ref{spherical-order-prop}(ii)). 
Recall that we want to construct a cyclic minimal $A_\infty$-structure on $\Ext^*(G,G)$,
where $G=\AA\oplus \rho_*V$. Let $L$ be a sufficiently positive power of an ample line bundle on $C$.
Then twisting $G$ through the spherical object $L^{-1}\ot\AA$ gives an $\AA$-module $\PP$ fitting into an exact
sequence
$$0\to \PP\to \Hom_\AA(L^{-1}\ot\AA,G)\ot L^{-1}\ot \AA \to G\to 0.$$
Since the local projective dimension of $G$ is $1$, this immediately implies that $\PP$ is locally projective.
Furthermore, since the spherical twist can be defined on a dg-level, we can replace $G$ by $\PP$ when studying
the minimal $A_\infty$-structure on $\Ext^*_\AA(G,G)\simeq\Ext^*_\AA(\PP,\PP)$ obtained by the homological perturbation.
Now, combining Proposition \ref{sheaf-dg-cyclic-prop} with Lemma \ref{trace-tau-lem}, 
we get that the minimal $A_\infty$-structure on $H^*(C,\und{\End}_\AA(\PP,\PP))$ obtained by the homological
perturbation can be chosen to be cyclic.
\ed

\end{document}